\documentclass[a4paper, 12pt]{article}
\usepackage{enumerate, theorem}
\usepackage{amsmath, amsfonts, amssymb}
\usepackage[height=22.5cm, width=15cm]{geometry}
\usepackage[all]{xy}
\usepackage{mathrsfs, yfonts}

\DeclareMathOperator{\id}{id}

\DeclareMathOperator{\Ker}{Ker}

\DeclareMathOperator{\C}{\mathbb{C}}

\newcommand{\A}{\tilde{\mathcal{A}}}

\newcommand{\parag}[1]{\paragraph{\sc{#1.}}}

\newtheorem{thm}{Theorem}[subsection]
\newtheorem{defn}[thm]{Definition}
\newtheorem{cor}[thm]{Corollary}
\newtheorem{prop}[thm]{Proposition}
\newtheorem{lemma}[thm]{Lemma}

\begin{document}

\title{The theme of a vanishing period.}

\author{Daniel Barlet\footnote{Barlet Daniel, Institut Elie Cartan UMR 7502  \newline
Nancy-Universit\'e, CNRS, INRIA  et  Institut Universitaire de France, \newline
BP 239 - F - 54506 Vandoeuvre-l\`es-Nancy Cedex.France. r\newline
e-mail : barlet@iecn.u-nancy.fr}.}

\date{23/09/11}

\maketitle

\section*{Abstract}

Let  \ $\lambda \in \mathbb{Q}^{*+}$ \  and consider a multivalued formal  function of the type
$$ \varphi(s) : =  \sum_{j=0}^k \ c_j(s).s^{\lambda + m_j}.(Log\, s)^j $$
where \ $c_j \in \C[[s]], m_j \in \mathbb{N}$ \ for \ $j \in [0,k-1]$. The {\bf theme} associated to such a \ $\varphi$ \ is the "minimal filtered differential equation" with generator \ $\varphi$, in a sens which is made precise in this article. We study such objects and show that their isomorphism classes may be characterized by a finite set of complex numbers, when we assume the Bernstein polynomial fixed. For a given \ $\lambda$, to fix the Bernstein polynomial is equivalent to fix a finite set of integers associated to the logarithm of the monodromy in the geometric stuation described above.\\
Our purpose is to construct some analytic invariants, for instance in the following situation : Let \ $f : X \to D$ \ be a proper holomorphic function defined on a complex manifold \ $X$ \ with value in a disc \ $D$. We assume that the only critical value is \ $0 \in D$ \ and we consider this situation as a degenerating family of compact complex manifolds to a singular compact complex space \ $f^{-1}(0)$. To a smooth \ $(p+1)-$form \ $\omega$ \ on \ $X$ \ such that \ $d\omega = 0 = df \wedge \omega$ \ and to a vanishing \ $p-$cycle \ $\gamma$ \ choosen in the generic fiber \ $f^{-1}(s_0), s_0 \in D \setminus \{0\}$, we associated a vanishing period \ $\varphi(s) : = \int_{\gamma_s} \ \omega\big/df $ \ which is, when \ $\gamma$ \ is choosen in the spectral subspace of \ $H_p(f^{-1}(s_0), \C)$ \ for the eigenvalue \ $e^{2i\pi.\lambda}$ \ of the monodromy of \ $f$, of the form above. Here \ $(\gamma_s)_{s \in D^*}$ is the horizontal multivalued family of \ $p-$cycles in the fibers of \ $f$ \ obtained from the choice of  \ $\gamma$.\\
The result obtained allows, for instance, to associate "natural" holomorphic functions of the parameter space when we have a family of such degenerations depending holomorphically on a parameter.

\parag{AMS Classification} 32S05, 32S25, 32S40.

\parag{Key words} Vanishing period, Bernstein polynomial, filtered Gauss-Manin system, (a,b)-module, Brieskorn module.

\tableofcontents

\newpage

\section{Introduction}

Let begin by the definition of a vanishing period. Let \ $X$ \ be a complex connected manifold of dimension \ $n + 1$ \ and \ $f : X \to D$ \ be a holomorphic function with value in a disc \ $D \subset \C$ \ with center \ $0$. Assume\footnote{if the hypersurface \ $\{ f = 0 \}$ \ in \ $X$ \ is reduced, this is always satisfied, up to localize around \ $0$ \ in \ $D$.} that the critical set \ $ S : = \{ df = 0 \}$ \ is a closed subset in \ $\{ f = 0 \}$ \ with no interior point in \ $\{ f = 0 \}$. Let \ $\omega$ \ be a \ $(p+1)-\mathscr{C}^{\infty}-$differential form on \ $X$ \ such that \ $d \omega = 0 = df\wedge\omega$, and let \ $\gamma \in H_p(X_{s_0}, \C)$ \ be a vanishing \ $p-$cycle on the generic fiber \ $X_{s_0}$ \ of \ $f$. Then we call  {\bf a vanishing period} \ for \ $f$ \  the multivalued\footnote{so well defined and holomorphic on the universal cover \ $exp : H \to D^*$.} holomorphic function on \ $D^*$, defined as 
$$ F_{\gamma}(s) : = \int_{\gamma_s} \ \frac{\omega}{df} $$
where \ $(\gamma_s)_{s \in H}, \gamma_s \in H_p(X_s, \C)$, is the horizontal family of \ $p-$cycles in the fibers of \ $f$ \ whose value\footnote{in fact at some base point in \ $H$ \ choosen over \ $s_0$.} at \ $s_0$ \ is \ $\gamma$.

 It is known that such a function \ $F_{\gamma}$ \  is solution of regular singular meromorphic linear differential system around \ $s = 0$ \ which is the Gauss-Manin connexion of \ $f$. It admits when \ $s \to 0$ \ a convergent asymptotic expansion (see [M.74] or [A-V-G.])  in the space
 
 $$ \Xi : = \oplus_{\lambda \in \mathbb{Q} \cap]0.1]} \Xi_{\lambda} \quad {\rm} \quad \ \Xi_{\lambda} \ \oplus_{j=0}^{\infty} \ \C[[s]].s^{\lambda-1}.(Log\,s)^j .$$

  Let \ $\A : = \{\sum_{\nu = 0}^{\infty} \ P_{\nu}(a).b^{\nu} \}$, where \ $P_{\nu} $ \ are polynomials, the \ $\C-$algebra whose product is defined by the following two conditions :
 \begin{enumerate}
 \item The commutation relation \ $a.b - b.a = b^2 $.
 \item The left and right multiplication by \ $a$ \ are continuous for the \ $b-$adic topology in \ $\A$.
 \end{enumerate}
 
 Then we associate to the vanishing period \ $F_{\gamma}$ \ the left sub$-\A-$module generated by the asymptotic expansion of \ $F_{\gamma}$ \ when \ $s \to 0$. This object, called a {\bf theme}, is in some sens the minimal filtered differential equation satisfied by \ $F_{\gamma}$. Note that the regularity of the Gauss-Manin connection, which implies the convergence of the asymptotic expansion of \ $F_{\gamma}$, shows that we do not loose information by going to the formal completion.\\
 
 The aim of this article is to study themes and more precisely to try to understand the classification of \ $[\lambda]-$primitive\footnote{A \ $[\lambda]-$primitive theme is a theme contained in \ $\Xi_{\lambda}$. So it keeps only the part of the asymptotic expansion associated to a the eigenvalue \ $e^{2i\pi.\lambda}$ \ of the monodromy of \ $f$.} themes up to isomorphy. The main goal is to try to characterize an isomorphy class of a \  $[\lambda]-$primitive theme by a finite set of complex numbers, when we begin to fix some discrete invariants (call the fundamental invariants) which corresponds to fix the Bernstein polynomial of the theme. In order to obtain such a classification, we study holomorphic families of \ $[\lambda]-$primitive themes, and we construct "minimal" versal families. Then we give some sufficient conditions in order that this family is universal. In this case, it implies that we describe exactely with a finite set of complex numbers the isomorphy class of theses themes.
 
 \bigskip
 
 The main results in this article are the following.
 
 \begin{itemize}
 
 \item The theorems \ref{theme 2} , \ref{theme 5} and \ref{tw. dual} of stability of themes by quotient, twisted duality and the characterization of \ $[\lambda]-$primitive themes by the uniqueness of  Jordan-H{\"o}lder sequence.
 
 \item The existence theorem \ref{canonical form 3} for the canonical form of a \ $[\lambda]-$primitive theme, which leads to the  construction of the canonical family,  and the versality theorem \ref{versal univ. 2}.
 
 \item The theorem of uniqueness  \ref{unique 1} of the canonical form for an "invariant "\ $[\lambda]-$primitive theme, which gives many cases where the canonical family is universal.
 
  \end{itemize}

\section{Preliminaries.}

\subsection{Regular and geometric (a,b)-modules.}

Let us first recall the definition of an (a,b)-module\footnote{For more details on these basic facts, see [B.93]}.

\begin{defn} An (a,b)-module \ $E$ \ is a free finite rank \ $\C[[b]]-$module endowed with an \ $\C-$linear map \ $ a : E \to E$ \ which satisfies the following two conditions :
\begin{itemize}
\item The commutation relation \ $a.b - b.a = b^2$.
\item The map \ $a$ \ is continuous for the \ $b-$adic topology of \ $E$.
\end{itemize}
\end{defn}

Remark that these two conditions imply that for any \ $S \in \C[[b]]$ \ we have
$$ a.S(b) = S(b).a + b^2.S'(b) $$
where \ $S'$ \ is defined via the usual derivation on \ $\C[[b]]$. For a given free rank \ $k$ \ $\C[[b]]-$module with basis \ $e_1, \dots, e_k$, to define a structure of (a,b)-module it is enough to prescribe (arbitrarily) the values of \ $a$ \ on \ $e_1, \dots, e_k$.\\

An alternative way to define  (a,b)-modules is to consider the \ $\C-$algebra 
$$ \A : = \{ \sum_{\nu = 0} ^{\infty} \ P_{\nu}(a).b^{\nu} \} $$
where the \ $P_{\nu}$ \ are polynomials in \ $\C[z]$ \ and where the product by \ $a$ \ is left and right continuous for the \ $b-$adic filtration and satisfies the commutation relation \ $a.b - b.a = b^2$.\\
Then a left \ $\A-$module which is free and finite rank on the subalgebra \ $\C[[b]] \subset \A$ \ is an (a,b)-module and conversely.\\
An (a,b)-module \ $E$ \  has a {\bf simple pole} when we have \ $a.E \subset b.E$ \ and it is {\bf regular} when it is contained in a simple pole (a,b)-module. The regularity is equivalent to the finitness on \ $\C[[b]]$ \ of the saturation \ $E^{\sharp}$ \ of \ $E$ \  by \ $b^{-1}.a$ \  in \ $E \otimes_{\C[[b]]} \C[[b]][b^{-1}]$. Recall that submodules and quotients of regular (a,b)-modules are regular.\\

The {\bf Bernstein polynomial} of a regular (a,b)-module \ $E$ \ of rank \ $k$ \  is defined as the minimal polynomial of \ $-b^{-1}.a$ \ acting on the \ $k-$dimensional \ $\C-$vector space \ $E^{\sharp}\big/b.E^{\sharp}$. Of course, when \ $E$ \ is the \ $b-$completion of the Brieskorn module of a non constant  germ \ $f : (\C^{n+1}, 0) \to (\C,0)$ \ of holomorphic function with an isolated singularity, we find the "usual" (reduced) Bernstein polynomial of \ $f$ \ (see for instance [K.76] or [Bj.93]).\\

We say that a regular (a,b)-module \ $E$ \ is { \bf geometric} when all roots of its Bernstein polynomial are negative rational numbers. This condition, which corresponds to \\
 M. Kashiwara theorem [K.76],   encodes the monodromy theorem  and the positivity theorem of  B. Malgrange (see [M.75] ) extending the situation of  (a,b)-modules deduced from the Gauss-Manin connexion of an holomorphic function.\\

Another important property of regular (a,b)-modules is the existence of Jordan-H{\"o}lder sequences (J-H. sequences for short).\\
 Recall first that any  regular rank 1 (a,b)-module is characterized, up to isomorphism, by a complex number \ $\lambda$ \ and the corresponding isomorphy class is represented by the (a,b)-module  \ $E_{\lambda} : = \C[[b]].e_{\lambda}$ \ where \ $a.e_{\lambda} = \lambda.b.e_{\lambda}$, which is isomorphic to the left \ $\A-$module \ $\A\big/\A.(a - \lambda.b)$.\\
Recall also that a submodule \ $F$ \ of the (a,b)-module \ $E$ \ is called  {\bf normal}  when \ $F \cap b.E = b.F$. Normality is a necessary and sufficient condition in order that  the quotient \ $E\big/F$ \ is again an (a,b)-module.\\
A {\bf Jordan-H{\"o}lder sequence} for the rank \ $k$ \ regular (a,b)-module \ $E$ \ is a sequence of normal submodules\footnote{For \ $G \subset F \subset E$ \ submodules with  $F$ \ normal in \ $E$, the normality of \ $G$ \ in \ $F$ \ is equivalent to the normality of \ $G$ \ in \ $E$.}\ $\{0\} = F_0 \subset F_1 \subset \dots F_{k-1} \subset F_k = E$ \ such that the quotients \ $F_j\big/F_{j-1}$ \ for \ $j \in [1,k]$ \ are rank 1 (a,b)-modules. So, to each J-H. sequence of \ $E$ \ we may associate an ordered \ $k-$tuple of complex numbers \ $\lambda_1, \dots, \lambda_k$ \ such \ $F_j\big/F_{j-1} \simeq E_{\lambda_j}$ \ for each \ $j \in [1,k]$.\\
The existence of J-H. sequence for any regular (a,b)-module is given via the  following proposition  proved in [B.93].

\begin{prop}\label{J-H.2}
Let \ $E$ \ be a regular (a,b)-module of rank \ $k$. Then \ $E$ \ admits a J-H. sequence. Moreover we may choose \ $F_1 \simeq E_{\lambda_1}$ \ with \ $-\lambda_1$ \ in any prescribed class modulo \ $\mathbb{Z}$ \ of a root of the Bernstein polynomial of \ $E$.
\end{prop}

\bigskip

Recall that the tensor product  of two (a,b)-modules \ $E$ \ and \ $F$ ( see [B.06] ) is defined as the \ $\C[[b]]-$module \ $E\otimes_{\C[[b]]}F$ \ with the \ $\C-$linear endomorphism defined by the rule \ $a.(x\otimes y) = (a.x)\otimes y + x \otimes (a.y)$. The tensor product by a fix (a,b)-module preserves short exact sequences of (a,b)-modules. As we have  \ $E_{\lambda}\otimes E_{\mu} \simeq E_{\lambda+\mu}$,  the tensor product of two regular (a,b)-modules is again regular.\\

In an analoguous way, for two (a,b)-modules \ $E$ \ and \ $F$ \ we define a new (a,b)-module as the \ $\C[[b]]-$module \ $Hom_{\C[[b]]}(E,F)$ \ endowed with the \ $\C-$linear map \ $a$ \ defined by
$$ (a.\varphi)(x) =  a_F.\varphi(x) -  \varphi(a_E.x) . $$
It is an easy exercice to see that \ $(a.\varphi)$ \ is again \ $\C[[b]]-$linear and that we have again  \ $a.b - b.a = b^2$ \ on \ $Hom_{\C[[b]]}(E,F)$. We shall denote by \ $Hom_{a,b}(E,F)$ \ the (a,b)-module obtained in this way. Note that after applying the functors \ $Hom_{a,b}(E, -)$ \ or \ $Hom_{a,b}(-, F)$ \ to  any short exact sequence of (a,b)-module  we obtain  again a short exact sequence of (a,b)-modules . Moreover, as \ $Hom_{a,b}(E_{\lambda}, E_{\mu}) \simeq E_{\lambda -\mu}$, we see that for \ $E,F$ \ regular (a,b)-modules the (a,b)-module \ $Hom_{a,b}(E,F)$ \ is also regular.

\begin{defn}\label{dual}
Let \ $E$ \ be a regular (a,b)-module. The dual \ $E^*$ \ of \ $E$ \ is defined as  \ $Hom_{a,b}(E, E_0)$ \ 
where \ $E_0 : = \A\big/\A.a \simeq \C[[b]].e_0$ \ with \ $a.e_0 = 0$.
\end{defn}

 We have \ $E_{\lambda}^* \simeq E_{-\lambda}$ \ and as the duality transforms a short exact sequence of (a,b)-modules in a short exact sequence, the dual of a regular (a,b)-module is regular. It is easy also  to see that for a regular (a,b)-module \ $E$ \ the canonical map \ $E \to (E^*)^*$ \ is an isomorphism. \\
 The dual of a regular (a,b)-module is again regular, but the dual of a geometric (a,b)-module is almost never geometric. To use duality in the geometric case we shall combine it with tensor product with \ $E_N$ \ where \ $N$ \ is a big enough integer (or rational number). Then \ $E^*\otimes E_N$ \ is geometric and \ $(E^*\otimes E_N)^*\otimes E_N \simeq E$. We shall refer to this process as "twisted duality". \\

Define now the left  \ $\A-$module of "formal multivalued expansions"
$$ \Xi : = \oplus_{\lambda \in \mathbb{Q}\cap ]0,1]} \ \Xi_{\lambda} \quad {\rm with} \quad \Xi_{\lambda} : =  \oplus_{j \in \mathbb{N}} \quad \C[[b]].s^{\lambda-1}.\frac{(Log\,s)^j}{j!} $$
with the  action of \ $a$ \ given by
$$ a.(s^{\lambda-1}.\frac{(Log\,s)^j}{j!}) = \lambda.\big[b\big(s^{\lambda-1}.\frac{(Log\,s)^j}{j!}\big) + b\big(s^{\lambda-1}.\frac{(Log\,s)^{j-1}}{(j-1)!}\big)\big]$$
for \ $j \geq 1$ \ and \ $a.s^{\lambda-1} = \lambda.b(s^{\lambda-1})$, with, of course, the commutation relations \ $a.S(b) = S(b).a + b^2.S'(b)$ \ for \ $S \in \C[[b]]$.\\
For any geometric (a,b)-module of rank \ $k$,  the vector space \ $Hom_{\A}(E, \Xi)$ \ is of dimension \ $k$ \ and this functor transforms short exact sequences of geometric (a,b)-modules in short exact sequences of finite dimensional vector spaces (see [B.05] for a proof). \\
In the case of the Brieskorn module of an isolated singularity germ of an holomorphic function \ $f$ \  at the origin of \ $\C^{n+1}$ \ this vector space may be identified with the n-th homology group (with complex coefficients) of the Milnor's fiber of \ $f$ (see [B.05]). The correspondance is given by associating to a (vanishing) cycle \ $\gamma$ \ the \ $\A-$linear map
$$ [\omega] \mapsto \big[ \int_{\gamma_s} \ \omega\big/df \big] \in \Xi $$
where \ $\omega \in \Omega^{n+1}_0$,  $\gamma_s$ \ is the multivalued horizontal family of \ $n-$cycles defined by \ $\gamma$ \ in the fibers of \ $f$, and where \ $[g]$ \ denotes the formal asymptotic expansion at \ $s = 0$ \ of the multivalued holomorphic function \ $g$.\\

\subsection{Primitive and coprimitive parts.}

\begin{defn}\label{Exponents}
Let \ $E$ \ be a regular (a,b)-module. We shall denote \\ 
$Exp(E)$ \ the subset of \ $ \C\big/\mathbb{Z}$ \  of classes modulo \ $\mathbb{Z}$ \ of the numbers \ $-\alpha$ \ where \ $\alpha$ \ is a root of the Bernstein polynomial of \ $E$.
\end{defn}

\parag{Remarks}
 \begin{enumerate}
\item For any regular (a,b)-module we have \ $Exp(E) = Exp(E^{\sharp})$ \ as the Bernstein polynomial of \ $E$ \ and of its saturation \ $E^{\sharp}$ \ coincide.
\item Let \ $E$ \ be any regular (a,b)-module. We have \ $[\lambda] \in Exp(E)$ \ if and only if there exists a \ $\lambda \in [\lambda]$ \ and an  (a,b)-linear injection of \ $E_{\lambda}$ \ in \ $E$ : it is enough to prove this for \ $E^{\sharp}$ \ and in the simple pole case we may choose the \ $\lambda$ \ with the smaller  real part in \ $[\lambda]$ \ for which \ $a - \lambda.b$ \ is not injective on \ $E$, thanks to the proposition 1.3 in [B.93].
\item Let \ $E$ \ be any regular (a,b)-module. Then \ $Exp(E^*) = - Exp(E)$.
\item Using the previous  remarks and the isomorphism of \ $E$ \ with its bidual, we obtain that \ $[\lambda]$ \ is in \ $Exp(E)$ \ if and only if there exists some \ $\lambda \in [\lambda]$ \ and an (a,b)-linear surjective map \ $E \to E_{\lambda}$.
\end{enumerate}

Using the previous remarks the following lemma is an easy exercice.

\begin{lemma}\label{exo} Let \ $0 \to F \to E \to G \to 0 $ \ be an exact sequence of regular (a,b)-modules. Then we have the equality \ $Exp(E) = Exp(F) \cup Exp(G)$.\\
\end{lemma}

\begin{defn}\label{primitive 1}
A regular (a,b)-module is {\bf \ $[\Lambda]-$primitive}, where \ $[\Lambda]$ \ is a subset  in \ $\C\big/\mathbb{Z}$, if all  roots of its Bernstein polynomial are in \ $[- \Lambda]$, or equivalently when \ $Exp(E) \subset \Lambda$.
\end{defn}

\parag{Notation} When \ $\Lambda = \{[\lambda]\}$ \ we shall say that \ $E$ \ is \ $[\lambda]-$primitive. When \ $M = (\C\big/\mathbb{Z}) \setminus \Lambda$, we shall replace \ $M-$primitive by \ $[\not= \Lambda]-$primitive.

\parag{Remarks}
\begin{enumerate}
\item With our definition the nul (a,b)-module is \ $\Lambda-$primitive for any choice of \ $\Lambda$. And it is the only regular (a,b)-module which is at the same time \ $\Lambda$ \ and \ $[\not= \Lambda]-$primitive.
\item Using remark 2 above, it is easy to see that any submodule (normal or not) of a \ $\Lambda-$primitive (a,b)-module is \ $\Lambda-$primitive.

\item As a consequence,  a regular (a,b)-module is \ $[\Lambda]-$primitive if and only if it admits a J-H. sequence such that the classes modulo \ $\mathbb{Z}$ \ of  \ $\lambda_1, \dots, \lambda_k$ \ are in \ $[\Lambda]$. If it is so, this property is true for any J-H. sequence, thanks to proposition \ref{J-H.2}.

\item If we have a short exact sequence of  (a,b)-modules
$$ 0 \to F \to E \to G \to 0 $$
with \ $E$ \ regular and  $[\Lambda]-$primitive then \ $F$ \ and \ $G$ \ are regular and \ $[\Lambda]-$primitive. \\
Conversely if \ $F$ \ and \ $G$ \ are regular and \  $[\Lambda]-$primitive then \ $E$ \ is regular and \ $[\Lambda]-$primitive.

\item Let \ $f : E \to F$ \ be an \ $\A-$linear map between two regular (a,b)-modules. Then if \ $G \subset E$ \ is a \ $\Lambda-$primitive submodule, then \ $f(G)$ \ is also \ $\Lambda-$primitive: if this is not the case, by remark 3 above  we may find an submodule of \ $f(G)$ \ isomorphic to \ $E_{\mu}$ \ with \ $[\mu] \not\in \Lambda$. But then, with  \ $H : = G \cap f^{-1}(E_{\mu})$, we have an exact sequence
$$ 0 \to \Ker f \cap G \to H \to E_{\mu} \to 0  $$
of regular (a,b)-modules, contradicting remark 4  above.
\end{enumerate}

\begin{prop}\label{primitive 2} Let \ $E$ \ be a regular (a,b)-module and let \ $\Lambda$ \ be a subset of  \ $ \C\big/\mathbb{Z}$. There exists a maximal \ $\Lambda-$primitive submodule \ $E[\Lambda]$ \ of \ $E$. It is normal and the quotient \ $E\big/E[\Lambda]$ \ is \ $[\not= \Lambda]-$primitive. Moreover, any normal submodule \ $F \subset E$ \ such that \ $E\big/F$ \ is \ $[\not=\Lambda]-$primitive contains \ $E[\Lambda]$. So, for any \ $M \subset \C\big/\mathbb{Z}$ \ the quotient \ $E\big/E[\not=M]$ \ is the maximal \ $[M]-$primitive quotient of \ $E$.
\end{prop}

\parag{Proof} We shall proof the proposition by induction on the rank \ $k$ \ of \ $E$. Assume \ $k \geq 1$ \ and the proposition proved in rank \ $\leq k-1$.\\
If any \ $[\lambda] \in -\Lambda$ \ is not the class modulo \ $\mathbb{Z}$ \ of a root of the Bernstein polynomial of \ $E$, it is clear that \ $\{0\}$ \ is the biggest \ $[\Lambda]-$primitive submodule in \ $E$. 
So we may assume that there exists a root \ $-\lambda$ \ of the Bernstein polynomial of \ $E$ \ with \ $[\lambda] \in \Lambda$. Up to change the choice of \ $\lambda$ \ in \ $[\lambda]$ \ we may assume that \ $E$ \ admits a normal submodule \ $F$ \ isomorphic to \ $E_{\lambda}$. Then \ $G : = E \big/F$ \ has rank \ $k-1$ \ and the induction hypothesis gives the existence of the maximal \ $[\Lambda]-$primitive submodule \ $G[\Lambda]$ \ of \ $G$. Let \ $\pi$ \ be the quotient map \ $E \to E\big/F $. Then we shall show that \ $\pi^{-1}(G[\Lambda])$ \ is a maximal \ $[\Lambda]-$primitive submodule of \ $E$.\\
The exact sequence \ $0 \to F \to \pi^{-1}(G[\Lambda]) \to G[\Lambda] \to 0 $, thanks to the previous remark 4, shows that \ $\pi^{-1}(G[\Lambda])$ \ is \ $[\Lambda]-$primitive.\\
If \ $H \subset E$ \ is a \ $[\Lambda]-$primitive submodule of \ $E$, then \ $\pi(H)$ \ is also \ $[\Lambda]-$primitive and so contained in \ $G[\Lambda]$, giving \ $H \subset \pi^{-1}(G[\Lambda])$.\\
The other assertions in the proposition are easy  defining \ $[M] : = [\not= \Lambda]$. $\hfill \blacksquare$

\parag{Remarks}
\begin{enumerate}
\item For any regular (a,b)-module \ $E$ \ and any \ $\Lambda \subset \C\big/\mathbb{Z}$ \ we have
$$ Exp(E[\Lambda]) = Exp(E) \cap \Lambda .$$
\item Let \ $F \subset E$ \ be a submodule of a regular (a,b)-module \ $E$. Then for any \ $\Lambda$ \ we have \ $F \cap E[\Lambda] = F[\Lambda]$ : the maximality of \ $F[\Lambda]$ \ gives \ $F \cap E[\Lambda]\subset F[\Lambda]$, and the maximality of \ $E[\Lambda]$ \ implies the other inclusion.\\
\end{enumerate}

\begin{defn}\label{primitive} We shall call \ $E[\Lambda]$ \ (resp. \ $E \big/E[\not= \Lambda]$) \ the {\bf \ $[\Lambda]-$primitive part} (resp. the {\bf $[\Lambda]-$coprimitive part}) of \ $E$. \\
\end{defn}

An easy consequence of the proposition \ref{primitive 2} is the following.

\begin{cor}\label{primitive 3} Let \ $E$ \ be a regular (a,b)-module and let \ $[\lambda_1], \dots, [\lambda_d]$ \ with \ $[\lambda_i] \not= [\lambda_j]$ \ for \ $i \not= j$ \ in \ $[1,d]$ \ be an order on the set \ $Exp(E)$. Then there exists an unique sequence of normal submodules of \ $E$ :
$$ 0 = F_0 \subset F_1 \subset \dots \subset F_d = E $$
such that \ $F_j\big/F_{j-1}$ \ is \ $[\lambda_j]-$primitive for each \ $j \in [1,d]$.
\end{cor}

Note that we have \ $F_j = E[\{[\lambda_1], \dots, [\lambda_j]\}]$.\\
But, in general for \ $j \geq 2$, the quotient \ $F_j\big/F_{j-1}$ \ is not isomorphic to \ $E[\lambda_j]$, as it may be seen on the rank 2  (a,b)-module \ $E_{\lambda,\mu} : = \A\big/\A.(a - \lambda.b)(a - \mu.b)$ \ when \ $[\lambda] \not= [\mu]$.

\begin{lemma}\label{primitive 4}  Let \ $E$ \ be a regular (a,b)-module and let \ $\Lambda$ \ be a subset of  \ $ \C\big/\mathbb{Z}$. The duality functor of (a,b)-modules transforms the exact sequence
$$ 0 \to E[\Lambda] \to E \to E\big/E[\Lambda] \to 0  $$
in the exact sequence
$$ 0 \to E^*[\not= -\Lambda] \to E^* \to (E[\Lambda])^* \to 0 .$$
So there is a canonical isomorphism \ $(E[\Lambda])^* \simeq E^*\big/[-\Lambda]$.
\end{lemma}

\parag{Proof} As the dual of a \ $[\Lambda]-$primitive (a,b)-module is \ $[-\Lambda]-$primitive, the universal property of the inclusion \ $E[\Lambda] \hookrightarrow E$ \ relative to \ $\A-$linear maps from \ $[\Lambda]-$primitive (a,b)-modules to \ $E$ \ gives that the surjection map \ $E^* \to (E[\Lambda])^*$ \ factorizes any \ $\A-$linear map from \ $E^*$ \ to a \ $[-\Lambda]-$primitive (a,b)-module. So \ $(E[\Lambda])^*$ \  is the \ $[-\Lambda]-$coprimitive part of \ $E^*$. This implies that the kernel of this quotient is the \ $[\not= -\Lambda]-$primitive part of \ $E^*$.\\

\section{Themes.}

\subsection{Stability by quotient and twisted duality.}

As it is explained in the introduction, a theme is, in the precise sens given below, a minimal filtered differential equation associated to a formal "standard multivalued asymptotic expansion". 

\begin{defn}\label{theme 1}
A {\bf theme} is a left \ $\A-$module isomorphic as \ $\A-$module to \ $\A.\varphi \subset \Xi$ \ for some \ $\varphi \in \Xi$.
\end{defn}

\parag{Remarks}
\begin{enumerate}
\item A theme is  a monogenic geometric (a,b)-module. Note that the finitness on \ $\C[[b]]$ \ comes immediately from our definition of \ $\Xi$ : there are only finitely many exponents  \ $\lambda \in ]0,1] \cap \mathbb{Q}$ \ involved in \ $\varphi$ \ and also the degree of logarithms are bounded. Moreover there is no \ $b-$torsion because \ $\Xi$ \ itself  has no \ $b-$torsion. So such an \ $\A.\varphi$ \ is an (a,b)-module, and the monogenic assertion is obvious. The regularity is an easy consequence of the fact that \ $\Xi$ \ has a simple pole and so there exists a simple pole sub-(a,b)-module of \ $\Xi$ \ which contains \ $\A.\varphi$, which is of finite type on \ $\C[[b]]$. The geometric aspect is then obvious.
\item We have proved that any  monogenic geometric (a,b)-module  may be embeded in \ $\Xi\otimes_{\C} V$ \ where \ $V$ \ is a finite dimensional complex vector space in [B.09], but nevertheless not every monogenic geometric (a,b)-module is a theme. \\
For instance, if \ $\lambda > 1$ \ is rational and \ $n$ \ is an integer the (a,b)-module \ $E : = \A\big/\A.(a - \lambda.b).(a - (\lambda+n).b)$ \ is a rank \ $2$ \ (a,b)-module which is not a theme : if \ $f : E \to \Xi$ \ is an \ $\A-$linear map, its \ $b-$rank is less or equal to \ $1$, for the following reason:\\
Let \ $e_1$ \ denote the image of \ $1$ \ in \ $E$ \ and put \ $e_2 : = (a-(\lambda+n).b).e_1$. Then we have \ $(a - \lambda.b).e_2 = 0 $ \ in \ $E$. Then \ $f(e_2) = \rho.s^{\lambda-1}$ \ for some \ $\rho \in \C$. Now it is easy to see that \ $\psi : = f(e_1)$ \ has to be a solution of the equation \ $(a - (\lambda+n).b).\psi = \rho.s^{\lambda-1}$. But if \ $\rho \not= 0$ \ this solution is \ $\psi = \sigma.s^{\lambda-2} + \tau.b^n.s^{\lambda-1}$ \ for some well choosen \ $\sigma \in \C$. So \ $f(e_1)$ \ is in \ $\C[[b]].f(e_2)$ \ and the rank on \ $\C[[b]]$ \ is at most \ $1$ \ for any \ $f$.
\item Any geometric (a,b)-module of rank 1 is a theme.
\item A theme is \ $[\lambda]-$primitive if and only if it may be embedded in \ $\Xi_{\lambda}$.
\item It is an easy exercice to show that if \ $\varphi \in \Xi_{\lambda}$ \ has degree \ $k-1$ \ in \ $Log\,s$ \ then \ $\A.\varphi$ \ is a ($[\lambda]-$primitive) theme of rank \ $k$.\\
\end{enumerate}

The next lemma contains the classification of \ $[\lambda]-$primitive themes of rank \ $2$.

\begin{lemma}\label{rang 2} Let \ $\lambda \in  ]0,1] $ \ be a rational number.
A rank \ $2$ \ \ $[\lambda]-$primitive theme is isomorphic to one and only one rank \ $2$  \ $[\lambda]-$primitive theme in the following list:
\begin{enumerate}
\item Let \ $\lambda_1 > 1$ \ in \ $[\lambda]$ \ and define \ $E_{\lambda_1,\lambda_1}$ \ as the quotient  $$ E_{\lambda_1,\lambda_1} : =  \A\big/\A.(a - \lambda_1.b).(a - (\lambda_1-1).b) .$$
\item Let \ $\lambda_1 > 1$ \ in \ $[\lambda]$, \ $n \in \mathbb{N}^*$ \ and \ $\alpha \in \C^*$ \ and define \ $E_{\lambda_1, \lambda_1+n}(\alpha)$ \ as the quotient 
 $$E_{\lambda_1, \lambda_1+n}(\alpha) : = \A\big/\A.(a - \lambda_1.b).(1 + \alpha.b^n)^{-1}.(a - (\lambda_1+n-1).b).$$
\end{enumerate}
In all cases, there exists an unique normal submodule of rank 1 and it is  isomorphic to \ $E_{\lambda_1}$.
\end{lemma}

\parag{Proof} Recall the classification of regular rank 2 (a,b)-module given in [B.93] proposition 2.4 :
\begin{enumerate}
\item  \ $E = E_{\lambda} \oplus E_{\mu}$.
\item For \ $n \in \mathbb{N}$ \ let   \ $E_{\lambda}(n) $ \ be the \ $\C[[b]]-$module with basis \ $(x,y)$ \ and where \ $a$ \ is defined by
$$ a.x = (\lambda+n).b.x + b^{n+1}.y \quad {\rm and} \quad  a.y = \lambda.b.y .$$
\item For \ $(\lambda,\mu) \in \C^2\big/\frak{s}_2$ \ let \ $E_{\lambda,\mu}$ \ be the \ $\C[[b]]-$module with basis \ $(x,t)$ \ and where \ $a$ \ is defined by
$$ a.y = \mu.b.t \quad {\rm and} \quad  a.t = y + (\lambda-1).b.t $$
\item For \ $n \in \mathbb{N}^*$ \ and \ $\alpha \in \C^*$ \ let  \ $E_{\lambda, \lambda-n}(\alpha)$ \  be  the \ $\C[[b]]-$module with basis \ $(y,t)$ \ and where \ $a$ \ is defined by
$$ a.y = (\lambda-n).b.y \quad {\rm and} \quad a.t = y + (\lambda-1).b.t + \alpha.b^n.y .$$
\end{enumerate}
 The cases 1 and 2 are not monogenic. The cases 4 corresponds to the case 2 in the lemma when \ $\lambda = \lambda_1 > 1$. So it is enough to prove that they may be embedded in \ $\Xi$. But it is an easy exercice to verify that it is isomorphic to \ $\A.\psi \subset \Xi$ \ with
 $$ \psi : =  s^{\lambda_1+n-2}.Log\, s + \gamma.s^{\lambda_1-2} \quad 
 {\rm with} \quad  \gamma = - \frac{(\lambda_1-1).\lambda_1\dots (\lambda_1+n-2)}{n} .$$
 We shall show that in case 3 we get a \ $[\lambda]-$primitive theme if and only if we are in  the case 1 of the lemma. Of course we may assume \ $\mu \in [\lambda]$ \ because we are looking for \ $[\lambda]-$primitive themes. Then it is enough (up to exchange \ $\lambda$ \ and \ $\mu$) to prove that \ $E_{\lambda,\lambda+p}$ \ is not a theme for \ $p \in \mathbb{N}^*$. But this has been proved in the remark 2 following the definition \ref{theme 1}. \\
 To check uniqueness of the normal rank 1 submodule is  an easy exercice left to the reader. $\hfill \blacksquare$\\

\parag{Remark} An easy consequence of the uniqueness of the normal rank 1 normal submodule is the uniqueness  of the J-H. sequence for a \ $[\lambda]-$primitive rank \ $2$ \  theme. Moreover in every case we have the inequatity \ $\lambda_1 \leq \lambda_2 - 1$, where \ $\lambda_1, \lambda_2$ \ are numbers associated to the rank 1  quotients of the J.H. sequence.\\

If it is not surprising that a monogenic submodule of a theme is again a theme, it is less obvious  that the quotient of a theme is again a theme, as it is shown in the next theorem.

\begin{thm}\label{theme 2}
Let \ $E$ \ be a theme and let  \ $F$ \ be a monogenic submodule of \ $E$. Then \ $F$ \ is a theme. Moreover, if \ $F$ \ is normal in \ $E$, the quotient \ $E\big/F$ \ is also a theme.
\end{thm}

\parag{Remark} A normal submodule of a monogenic (a,b)-module is monogenic : the normality of \ $F$ \ in \ $E$ \ implies the injectivity of the map \ $F\big/b.F \to E\big/b.E$. As the action of \ $a$ \ on \ $E\big/b.E$ \ is a principal nilpotent \ $\C-$linear endomorphism, the \ $a-$stable subspace \ $F\big/b.F$ \ is equal to \ $Im(a^h)$ \ for some integer \ $h$. So \ $F\big/a.F + b.F$ \ is 1-dimensional. This implies that \ $F$ \ is monogenic.\\

The proof of the theorem  will use the following proposition.

\begin{prop}\label{theme 3}
Let \ $E$ \ be a non zero  \ $[\lambda]-$primitive theme. Then \ $E$ \ has an unique normal rank 1  submodule. Moreover, if \ $F \simeq E_{\lambda_1}$ \ is this normal rank 1 submodule the quotient \ $E\big/F$ \ is again a \ $[\lambda]-$primitive theme.
\end{prop}

\parag{Proof} The existence of J-H. sequences for regular (a,b)-modules shows that there exists at least one normal rank 1 submodule of \ $E \not= \{0\}$. As \ $E$ \ is \ $[\lambda]-$primitive, this normal rank 1 submodule \ $F_1$ \ is isomorphic to \ $E_{\lambda_1}$ \ for some \ $\lambda_1 \in [\lambda]$.\\
Assume that we have another normal rank 1 submodule \ $G$ \ of \ $E$. Put \ $H : = F_1 + G$. Then, assuming \ $F_1 \not= G$ \ the rank of \ $H$ \ is two : if not, we have \ $F_1 = b^p.H$ \ and \ $G = b^q.H$. But the normality of \ $F_1$ \ and \ $G$ \ imply \ $p = q = 0$ \ contradicting \ $F_1 \not= G$.\\
So \ $H$ \ is a \ $[\lambda]-$primitive theme of rank \ $2$. But then it has a unique rank 1 normal submodule, contradicting again \ $F_1 \not= G$. This proves the uniqueness of \ $F_1$.\\
To show that the quotient \ $E\big/F_1$ \ is again a theme, put for \ $N \in \mathbb{N}$
$$ \Xi^{(N)}_{\lambda} : = \oplus_{j=0}^N \ \C[[b]].s^{\lambda-1}.(Log\,s)^j \subset \Xi_{\lambda} .$$
This is a \ $[\lambda]-$primitive geometric submodule of \ $\Xi_{\lambda}$ \ for \ $\lambda \in ]0,1] \cap \mathbb{Q}$. And we have the exact sequence of \ $\A-$modules, for any \ $N \geq 1$ :
$$ 0 \to \C[[b]].s^{\lambda-1} \to \Xi^{(N)}_{\lambda} \overset{f_{\lambda}}{\longrightarrow}  \Xi^{(N-1)}_{\lambda} \to 0 $$
where \ $f_{\lambda}$ \ is \ $\C[[b]]-$linear and defined by 
 $$f_{\lambda}(s^{\lambda-1}) = 0 \quad {\rm and} \quad  f_{\lambda}(s^{\lambda-1}.(Log\,s)^j) : = j.s^{\lambda-1}.(Log\,s)^{j-1}.$$
 It is easy to check  (but not obvious)  that \ $f_{\lambda}$ \ is \ $a-$linear, surjective and with kernel \ $\C[[b]].s^{\lambda-1}$.\\
Consider now a non zero  \ $[\lambda]-$primitive theme \ $E \hookrightarrow \Xi_{\lambda}^{(N)}$ \ for some \ $N$. Its unique normal submodule \ $F_1$ \ is send bijectively on  \ $\C[[b]].s^{\lambda+p-1}$ \ for some \ $p \in \mathbb{N}$ \ (in fact \ $\lambda+p-1 = \lambda_1$). It coincides with \ $E \cap \Ker f_{\lambda}$\ because  \ $E \cap \Ker f_{\lambda}$ \ is normal and rank 1 in \ $E$. So the restriction of \ $f_{\lambda}$ \ to \ $E$ \ induces an injection of \ $E\big/F_1$ \ in \ $\Xi_{\lambda}^{(N-1)}$. $\hfill \blacksquare$\\

\parag{Proof of theorem \ref{theme 2}} The first assertion is obvious. As the quotient of a monogenic geometric (a,b)-module is again monogenic\footnote{See the remark following the statement of theorem \ref{theme 2}.} and geometric, the point is to construct an embedding of \ $E$ \ in \ $\Xi$. We shall do it by induction on the rank of \ $F$. The case of rank 1 is an easy consequence of the proposition \ref{theme 3} : assume that \ $F \simeq E_{\lambda_1}$ \ with \ $\lambda_1 \in [\lambda]$ and define  \ $ f : \Xi^{(N)} \to \Xi^{(N)} $ \  as the direct sum of the identity on \ $\Xi_{\mu}^{(N)}$ \ for \ $\mu \in ]0,1] \cap \mathbb{Q}, \mu \not= \lambda$ \ and the composition of \ $f_{\lambda}$ \ with the obvious inclusion \ $\Xi_{\lambda}^{(N-1)} \hookrightarrow \Xi_{\lambda}^{(N)}$. Then if \ $g : E \hookrightarrow \Xi^{(N)}$ \ is an embedding, then \ $f\circ g : E\big/F \to \Xi^{(N)}$ \ is also an embedding.\\
Assume now that the theorem is proved for rank \ $\leq k-1$ \ and that \ $F$ \ has rank \ $k \geq 2$. Then choose \ $G$ \ a normal rank 1 submodule in \ $F$. As we have \ $E\big/F \simeq (E\big/G)\Big/(F\big/G)$ \ the induction hypothesis allows to conclude.$\hfill \blacksquare$

\begin{cor}\label{theme 3bis}
Let \ $E$ \ be a monogenic geometric (a,b)-module. Then \ $E$ \ is a theme if and only if for any \ $[\lambda] \in Exp(E)$ \ its \ $[\lambda]-$coprimitive part \ $E\big/E[\not=\lambda]$ \ is a theme.
\end{cor}

We shall show that the analoguous statement with the \ $[\lambda]-$primitive parts of \ $E$ \ is also true, using the twisted duality theorem \ref{tw. dual}.

\parag{Proof}The previous theorem implies that the condition is necessary. To show that it is sufficient, let \ $Exp(E) = \{[\lambda_1], \dots, [\lambda_d]\}$. For \ $N \in \mathbb{N}$ \ large enough let \ $\theta_i : E \to \Xi_{\lambda_i}^{(N)}$ \ the composition of the quotient map \ $E \to E\big/E[\not=\lambda_i]$ \ with an injection of this quotient, which is, by assumption a \ $[\lambda_i]-$primitive theme, in \ $\Xi_{\lambda_i}^{(N)}$. Now \ $\theta : = \oplus_{i=1}^d \ \theta_i : E \to \Xi^{(N)}$ \ is an embedding : by construction we have \ $Ker \theta_i = E[\not=[\lambda_i]]$ \ and as \ $E[\cap_i P_i] = \cap_i E[P_i]$ \ the kernel of \ $\theta$ \ is \ $\{0\}$. $\hfill \blacksquare$\\

\begin{prop}\label{theme 4}
\begin{enumerate}
\item Let \ $E$ \ be a rank \ $k$ \ $[\lambda]-$primitive theme. Then there exists\footnote{with the convention \ $\Xi_{\lambda}^{(-1)} = \{0\}$.} \ $\varphi \in \Xi_{\lambda}^{(k-1)} \setminus \Xi_{\lambda}^{(k-2)}$ \ such that \ $E$ \ is isomorphic to \ $\A.\varphi \subset \Xi_{\lambda}^{(k-1)}$.
\item Conversely, for any such \ $\varphi$, $\A.\varphi$ \ is a rank \ $k$ \ $[\lambda]-$primitive theme. 
\item In this situation, for any \ $j \in [1,k]$ \ $F_j : = \A.\varphi \cap \Xi^{(j-1)}_{\lambda}$ \ is a rank \ $j$ \ normal submodule of \ $\A.\varphi$ \ contained in \ $b^{k-j}.\Xi^{(j-1)}_{\lambda}$.
\end{enumerate}
\end{prop}

\parag{Proof} We shall let as an exercie for the reader  the following  easy fact :\\
Fix \ $\lambda \in ]0,1] \cap \mathbb{Q}, k \in \mathbb{N}^*$ \ and \ $\mu \in \mathbb{Q}^{+*}$. 
Then we have:
\ $$(a - \mu.b).\Xi_{\lambda}^{(k-1)} \cap \C.s^{\mu}.(Log\,s)^{k-1} = \{0\}.$$

By definition of a \ $[\lambda]-$primitive theme there exist an integer \ $N $ \ large enough and \ $\varphi \in \Xi_{\lambda}^{(N)}$ \ such that \ $\A.\varphi$ \ is isomorphic to \ $E$. So the first two points in the proposition reduce to show that the rank of \ $\A.\varphi$ \ is (for \ $\varphi \not= 0$)  the maximal power of \ $Log\,s$ \ which appears in \ $\varphi$ \ plus 1. We shall show this fact by induction on the maximal power \ $d$ \ of \ $Log\,s$ \ appearing in \ $\varphi$.\\
The result is clear for \ $d = 0$. So assume the result proved for \ $d-1\geq 0$ \ and consider a \ $\varphi$ \ for which the maximal power is \ $d$. Then we may write
$$ \varphi = T.s^{\mu-1}.(Log\,s)^d + \psi $$
where \ $T\in \C[[b]]$ \ is invertible, \ $\mu \in \lambda + \mathbb{N}$ \ and \ $\psi $ \ is in \ $\Xi_{\lambda}^{(d-1)}$. To replace \ $\varphi$ \ by \ $T^{-1}.\varphi$ \ does not change \ $\A.\varphi$ \ so we may assume that \ $T = 1$. Then we have
$$\theta : =  (a - \mu.b).\varphi \in \Xi_{\lambda}^{(d-1)}$$
and is not in \ $\Xi^{(d-2)}_{\lambda}$ \ because of the  exercice above. So the induction hypothesis gives that \ $\A.\theta$ \ has rank \ $d$. \\
We have a surjective \ $\A-$linear map, where \ $E : = \A.\varphi$ \ and \ $F : = \A.\theta$
$$ E\big/F \to E_{\mu} \to 0 $$
obtained by sending \ $\varphi$ \ onto the standard generator of \ $E_{\mu}$, because, by definition, \ $(a - \mu.b).\varphi \in F$.
We want to prove now that \ $\A.\theta$ \ is exactely the kernel of this map. Let \ $x \in \A$.

The euclidian division gives \ $q \in \A$ \ and \ $r \in \C[[b]]$ \ such that
$$ x = q.(a - \mu.b) + r(b) $$
and this leads to \ $x.\varphi = q.\theta + r(b).\varphi$. So \ $x.\varphi$ \ is in the kernel if and only if \ $r(b).\varphi$ \ is in the kernel. But its image is \ $r(b).e_{\mu}$ \ which is zero only when \ $r = 0$.

This implies that the rank of \ $E$ \ is \ $d+1$.\\
But this proof also gives that \ $\A.\varphi \cap \Xi_{\lambda}^{(d-1)}$ \ is a normal submodule of \ $\A.\varphi$ \ of rank \ $d$. An easy iteration implies then the first part of assertion 3. \\
We have also the second part of \ 3 \ for \ $j = k-1 = d$ \ because \ $\theta = (a - \mu.b).\varphi$ \ is in \  $b.\Xi_{\lambda}^{(d)} \,\cap\, \Xi_{\lambda}^{(d-1)}$ \ as \ $\Xi_{\lambda}^{(N)}$ \ has a  simple pole for any \ $N$.\\
 Now we conclude that \ $\theta \in b.\Xi_{\lambda}^{(d-1)}$ \ because \ $\Xi_{\lambda}^{(N-1)}$ \ is normal in \ $\Xi_{\lambda}^{(N)}$ \ for any \ $N$.\\
By an iteration of this result we obtain 3. $\hfill \blacksquare$\\

\parag{Remark} With the notations of the previous proposition, if \ $s^{\mu-1}$ \ is in \ $\A.\varphi$ \ we have \ $\mu > k-1$, where \ $k$ \ is the rank of \ $\A.\varphi$.\\

\begin{thm}\label{theme 5}
Let \ $E$ \ be a monogenic geometric (a,b)-module having an unique normal rank 1  submodule. Then \ $E$ \ is a \ $[\lambda]-$primitive for some \ $\lambda \in ]0,1] \cap \mathbb{Q}$.
\end{thm}

\parag{Proof} We make an induction on the rank of \ $E$. The rank 1 case is obvious. Assume that \ $E$ \ has  rank \ $k +1 \geq 2$ \ and that the result is proved for rank \ $\leq k$. Let \ $F$ \ be a normal rank \ $k$ \ submodule of \ $E$. Then we have an exact sequence of monogenic geometric (a,b)-modules
\begin{equation*}
 0 \to F \to E \to E_{\lambda'} \to 0  \tag{@}
 \end{equation*}
where \ $\lambda' > 0 $ \ is a rational number. The induction hypothesis shows that \ $F$ \ is a \ $[\lambda]-$primitive theme. If \ $\lambda' \not\in [\lambda]$ \ then \ $E$ \ has a rank 1 normal submodule\footnote{see proposition \ref{J-H.2}.} isomorphic to \ $E_{\lambda''}$ \ with \ $\lambda'' \in [\lambda']$ \ and as \ $F$ \ has a rank 1 normal submodule isomorphic to \ $E_{\lambda_1}$ \ with \ $\lambda_1 \in [\lambda]$, we contradict our assumption.\\
Fix an \ $\A-$linear  injection \ $ j : F \to \Xi^{(N)}_{\lambda}$ \ and consider the exact sequence of finite dimensional vector spaces deduced from \ $(@)$, thanks to [B.05]  thm 2.2.1 p.24:
$$ 0 \to Hom_{\A}(E_{\lambda'}, \Xi) \to Hom_{\A}(E, \Xi) \to Hom_{\A}(F, \Xi) \to 0 . $$
Let \ $j' \in Hom_{\A}(E, \Xi)$ \ extending the image of \ $j$ \ in \ $Hom_{\A}(F, \Xi)$. It is now enough to prove that \ $j'$ \ is injective. Let \ $G : = Ker\,j'$, and assume that \ $G$ \ is not \ $\{0\}$. As \ $G \cap F = \{0\}$ \ the rank of \ $G$ \ is 1 and it is normal in \ $E$. But, as \ $k \geq 1$ \ $F$ \ has a normal rank 1 submodule, which cannot be \ $G$. So we contradict our assumption. So \ $G = \{0\}$ \ and \ $j'$ \ is injective. $\hfill \blacksquare$

\begin{cor}\label{theme 6}
Let \ $E$ \ be a \ $[\lambda]-$primitive theme of rank \ $k$. Then \ $E$ \ has exactly one normal submodule of rank \ $j$ \ for each \ $j \in [0,k]$. This implies the following properties:
\begin{enumerate}[i)]
\item The Jordan-H{\"ol}der sequence of \ $E$ \ is unique.
\item If we choose for each \ $j \in [1,k]$ \ an injection \ $\A-$linear \ $\theta_j : E\big/F_{k-j} \to \Xi_{\lambda}$ \ of the ($[\lambda]-$primitive) theme \ $E\big/F_{k-j}$, then \ $\theta_1, \dots, \theta_k$ \ is a basis of the complex vector space \ $Hom_{\A}(E, \Xi)$.
\end{enumerate}
\end{cor}

\parag{Proof} We shall prove  by induction on \ $j \in [0,k]$ \   the uniqueness of the normal rank \ $j$ \ submodule of a   \ $[\lambda]-$primitve theme of rank \ $k$.\\
The case \ $j = 1$ \ is proved in the previous theorem. Assume \ $j \geq 2$ \ and the assertion proved for rank \ $\leq j-1$. Consider  two normal submodules \ $F$ \ and \ $G$ \ of rank \ $j$ \ in \ $E$. As \ $j \geq 2$ \ they both contain the unique rank 1 normal submodule \ $H$ \ of \ $E$. But then \ $E\big/H$ \ is  a rank \ $k-1$ \ $[\lambda]-$primitive theme and as  \ $F\big/H$ \ and \ $G\big/ H$ \ are two rank \ $j-1$ \ normal sumodules of \ $E\big/H$, the inductive assumption implies there equality. So \ $F = G$ \ and the first assertion of the corollary is proved.\\
Then the properties i) and ii) are immediate. $\hfill \blacksquare$

\parag{Remarks}
\begin{enumerate}
\item The uniqueness of the J-H. sequence implies that there is a canonical order fixed on the roots of the Bernstein polynomial of a \ $[\lambda]-$primitive theme.
\item We proved in [B.09] proposition 3.5.2 that any monogenic geometric (a,b)-module has an J-H. sequence such that the corresponding  numbers \ $\lambda_1, \dots, \lambda_k$ \ satisfy the following property
\begin{itemize}
\item The sequence \ $\lambda_j + j$ \ is increasing (large).
\end{itemize}
So, when \ $E$ \ is a \ $[\lambda]-$primitive theme, the numbers \ $\lambda_1, \dots, \lambda_k$ \ associated to the unique J-H. sequence of \ $E$ \ satisfy this property and we may write
$$ \lambda_{j+1} = \lambda_j + p_j - 1 \quad \forall j \in [1,k-1] $$
where \ $p_j $ \ is in \ $\mathbb{N}$.
\item As we have seen (see the  remark following proposition \ref{theme 4}) that for a rank \ $k$ \  $[\lambda]-$primitive theme we have   \ $\lambda_1 >  k-1$ \ the inequatity  \ $\lambda_1 + 1 \leq \lambda_j + j$ \ implies
$$ \lambda_j > k-j \quad \forall j \in [1,k] .$$
\item The property ii) of the previous corollary shows that the  stratification by the rank defines a full flag in \ $Hom_{\A}(E, \Xi)$.
\end{enumerate}

\begin{defn}\label{Fund. inv.}
Let \ $E$ \ be a \  $[\lambda]-$primitive theme. We shall call {\bf fundamental invariants} of \ $E$ \ the (ordered) sequence \ $\lambda_1, \dots, \lambda_k$ \ associated to its unique J-H. sequence. We shall often fix the fundamental invariants by giving \ $\lambda_1, p_1, \dots, p_{k-1}$.
\end{defn}

\bigskip

\begin{thm}[Twisted duality]\label{tw. dual}
Let \ $E$ \ be a \ $[\lambda]-$primitive theme with fundamental invariants \ $\lambda_1, \dots, \lambda_k$. For any rational number \ $\delta$ \ satisfying \ $\delta > \lambda_k + k - 1 $ \ the (a,b)-module \ $E^* \otimes_{a,b}E_{\delta} $ \ is a \ $[\delta-\lambda]-$primitive theme with fundamental invariants \ $\delta -\lambda_k, \dots, \delta - \lambda_1$, where \ $E^*$ \ denotes the dual\footnote{see section 2.1.} of \ $E$.
\end{thm}

\parag{Proof} First the dual of a  regular (a,b)-module is  regular thanks to the fact that by duality a J-H. sequence gives a J-H. sequence and the equality \ $(E_{\lambda})^* \simeq E_{-\lambda}$ \ for any \ $\lambda \in \C$. The fact that \ $E^*$ \ is again monogenic when \ $E$ \ is monogenic regular is obtained  as follows :\\
It is proved in [B.08] section 3.3 (and also in [K.09] ) that there is, for an (a,b)-module \ $E$ \  a canonical isomorphism
$$ E^* \simeq Ext^1_{\A}(E, \A) $$
where the \ $\C-$vector space \ $Ext^1_{\A}(E, \A)$ \ is endowed with a left \ $\A-$module structure deduced from the right module structure on \ $\A$ \ using the anti-automorphism \\
 $\tau : \A \to \A$ \ defined by \ $\tau(1) = 1, \tau(a) = -a , \tau(b) = b $ \ and \ $\tau(x.y) = \tau(y).\tau(x)$. \\ 
Now the fact that \ $E$ \ is an monogenic (a,b)-module correspond to a resolution\footnote{The structure theorem \ref{struct. th.} precises are  the \ $P \in \A$ \ which may  appears here.}
$$ 0 \to \A \overset{.P}{\longrightarrow} \A \to E \to 0$$
which gives the resolution of \ $E^*$ \ as a right \ $\A-$module,
$$ 0 \to \A \overset{P.}{\longrightarrow} \A \to E^* \to 0 $$
 using the fact that \ $Ext^0_{\A}(E, \A) = 0$. This implies that \ $E^*$ \ is a monogenic (a,b)-module.\\
Then \ $E^*$ \ and also \ $E^*\otimes_{a,b} E_{\delta}$ \ for any \ $\delta$ \ has only one normal submodule for each rank \ $j \in [0,k]$ \ where \ $k : = rk(E)$. Dualizing and twisting  the J-H. sequence of \ $E$ \ gives the unique  J-H. sequence of \ $E^*\otimes_{a,b} E_{\delta}$. Now to have a \ $[\delta - \lambda]-$primitive, it is enough to have a geometric (a,b)-module. This is acheive if and only if the rational number \ $\delta$ \ satisfies  \ $\delta > \lambda_k + k - 1 $. $\hfill \blacksquare$

\parag{Remark} For a general theme, it is easy to deduce from the previous theorem that, for \ $\delta$ \ large enough in \ $\mathbb{Q}$,  $E^*\otimes_{a,b}E_{\delta}$ \ is again a theme because the \ $[\lambda]-$coprimitive part of \ $E^*\otimes_{a,b}E_{\delta}$ \ is \ $(E[-\lambda])^*\otimes_{a,b} E_{\delta}$ \ for any \ $[\lambda] \in \C\big/\mathbb{Z}$. We conclude using corollary \ref{theme 3bis}.

Then the twisted duality theorem combined with corollary \ref{theme 3bis} gives the following characterization of a theme.

\begin{cor}\label{theme 6}
Let \ $E$ \ be a monogenic geometric (a,b)-module. It is a theme if and only if for each \ $[\lambda] \in Exp(E)$ \ its \ $[\lambda]-$primitive part \ $E[\lambda]$ \ is a theme.
\end{cor}

\subsection{Standard and canonical forms for a \ $[\lambda]-$primitive theme.}

The structure theorem of [B.09] (thm 3.4.1)  for monogenic regular (a,b)-modules gives the following structure theorem for \ $[\lambda]-$primitive themes.

\begin{thm}\label{struct. th.}
Let \ $E$ \ be a \ $[\lambda]-$primitive theme with fundamental invariants \ $\lambda_1, p_1, \dots, p_{k-1}$. Then there exists \ $S_1, \dots, S_{k-1}$ \ in \ $\C[b]$ \ such that \ $S_j(0) = 1$ \ and \ $deg(S_j) \leq p_j+ \dots + p_{k-1}$ \ for all \ $j \in [1,k-1]$ \ and with \ $E \simeq \A\big/\A.P$ \ where
$$ P : = (a - \lambda_1.b).S_1^{-1}.(a - \lambda_2.b).S_2^{-1} \dots S_{k-1}^{-1}.(a - \lambda_k.b) .$$
Moreover, for each \ $j$ \ the coefficient of \ $b^{p_j}$ \ in \ $S_j$ \ is not zero.\\
Conversely, for any \ $\lambda_1 > k-1$ \ rational and any natural integers \ $p_1, \dots, p_{k-1}$, any \ $S_1, \dots, S_{k-1}$ \ invertible elements in \ $\C[[b]]$ \  such that the coefficient of \ $b^{p_j}$ \ in \ $S_j$ \ is not zero, the quotient \ $\A\big/\A.P$ \ is a rank \ $k$ \ \ $[\lambda]-$primitive theme.
\end{thm}

\parag{Proof} The direct part is an immediate consequence of the theorem 3.4.1 and of the lemma 3.5.1 of [B.09], using the uniqueness of the J-H. sequence of a \ $[\lambda]-$primitive theme.\\
Let  us show the converse. The quotient \ $\A\big/\A.P$ \ is clearly a rank \ $k$ \  monogenic geometric (a,b)-module. We shall show that it is a theme by induction on \ $k$. As the case \ $k = 1$ \ is obvious, assume the result proved for \ $k-1 \geq 1$. Put \ $Q : = (a - \lambda_1.b).S_1^{-1} \dots S_{k-2}^{-1}.(a - \lambda_{k-1}.b)$. The induction hypothesis gives that \ $F : = \A\big/\A.Q$ \ is a \ $[\lambda]-$primitive theme. Let \ $\varphi \in \Xi^{(N)}_{\lambda}$ \ such that \ $F$ \ is isomorphic to \ $\A.\varphi$. To construct \ $\psi \in \Xi^{(N+1)}_{\lambda}$ \ such that \ $(a - \lambda_k.b).\psi = S_{k-1}.\varphi$ \ it is enough to solve an elementary differential equation. Explicitely \ $s.f'(s) - \lambda_k.f(s) = S_{k-1}.\varphi$, where we put \ $b.\psi : = f$. Le main point is that, as \ $F$ has rank \ $k-1$, we may, up to change \ $\varphi$ \ by an invertible element in \ $\C[[b]]$, assume that \ $\varphi$ \ is a polynomial in \ $Log\, s$ \ of degree \ $k-2$ \ with leading coefficient \ $s^{\lambda_{k-1}-1}$. Then we see that the degree in \ $Log\, s$ \ of \ $\psi$ \ is exactely \ $k-1$ \ because the coefficient of \ $b^{p_{k-1}}$ \ in \ $S_{k-1}$ \ does not vanish  and we have \ $\lambda_k = \lambda_{k-1} + p_{k-1} -1$.\\
Then the \ $\A-$linear map \ $\A\big/\A.P \to \Xi_{\lambda}^{(N+1)}$ \ defined by sending \ $1$ \ to \ $\psi$ \ will have an image \ $\A.\psi$ \ which is a \ $[\lambda]-$primitive theme of rank \ $k$. As it is surjective between free \ $\C[[b]]-$modules of the same rank, it is an isomorphism. So \ $E = \A\big/\A.P $ \ is a theme. $\hfill \blacksquare$

\parag{Notation} When we consider a \ $[\lambda]-$primitive theme with fundamental invariants \ $\lambda_1, p_1, \dots, p_{k-1}$, we shall say that an isomorphism of \ $E \simeq \A\big/\A.P$ \ gives a {\bf standard form} for \ $E$ \ when \ $P$ \ is as in the statement above, but omitting the condition on the degree of the \ $S_j$.\\

\begin{prop}\label{base standard 1}
Let \ $E \simeq \A\big/\A.P$ \ a standard form for a \ $[\lambda]-$primitive theme. Then for each \ $j \in [0,k]$ \ define
$$ P_j : = (a- \lambda_{j+1}.b)S_{j+1}^{-1} \dots S_{k-1}^{-1}.(a -\lambda_k.b) $$
and put \ $e_j : = P_j.[1]$ \ where\footnote{As \ $P_0 = P$ \ we have \ $e_0 = 0$ \ and \ $P_k = 1$ \ shows that \ $e_k = [1]$.} \ $[1]$ \ is the image of the class of \ $1 \in \A$ \ via the isomorphism \ $E \simeq \A\big/\A.P$. Then \ $e_1, \dots, e_k$ \ is a \ $\C[[b]]-$base of \ $E$ \ and \ $a$ \ is defined on \ $E$ \ by the relations 
\begin{equation*}
a.e_j = \lambda_j.b.e_j + S_j(b).e_{j-1} \quad \forall j \in [1,k]. \tag{1}
\end{equation*}
The rank \ $j$ \ normal submodule \ $F_j$ \ of \ $E$ \ admits \ $e_1, \dots, e_j$ \ as a \ $\C[[b]]-$basis.\\
Any \ $x \in F_{k-j}$ \ which is annihilated by \ $P_0$ \ satisfies \ $P_j.x = 0$ \ in \ $E$. Moreover there exists \ $\rho \in \C$ \ such that \ $x - \rho.b^{\lambda_k-\lambda_{k-j}}.e_{k-j}$ \ belongs to \ $F_{k-j-1}$.
\end{prop}

Remark that if \ $\lambda_k - \lambda_{k-j} < 0$ \ this means that \ $\rho = 0 $ \ and \ $x \in F_{k-j-1}$.

\parag{Proof} Note first that the class of \ $e_j$ \ in \ $E\big/b.E$ \ is equal to \ $[a^{k-j}.1]$, so \ $e_1, \dots, e_k$ \   induces a basis of \ $E\big/b.E$.  Then \ $e_1, \dots, e_k$ \ is a \ $\C[[b]]-$basis of \ $E$. Define \ $F_j$ \ as the \ $\C[[b]]-$submodule of \ $E$ \ generated by \ $e_1, \dots, e_j$. The relations \ $(1)$ \ shows that \ $F_j$ \ is stable by \ $a$ \ so it is a \ $\A-$submodule. As \ $E\big/F_j$ \ has no \ $b-$torsion, \ $F_j$ \ is the unique normal rank \ $j$ \ submodule of \ $E$.\\
Let \ $x \in F_{k-j}$ \ such that \ $P_0.x = 0$ \ in \ $E$. Then sending \ $e_k$ \ to \ $x$ \ defines a surjective \ $\A-$linear map from \ $E = \A\big/\A.P_0$ \ to \ $F_{k-j}$. The kernel is normal and has rank \ $j$, so it is \ $F_j$. But as \ $P_j.e_k = S_j.e_j$ \ lies in \ $F_j$, we must have \ $P_j.x = 0$ \ in \ $E_{k-j}$.\\
Then the image of \ $x$ \ in \ $F_{k-j}\big/F_{k-j-1} \simeq E_{\lambda_{k-j}}$ \ is in the kernel of \ $P_j$ \ acting on \ $E_{\lambda_{k-j}}$. This kernel is equal to \ $Hom_{\A}(\A\big/\A.P_j, E_{\lambda_{k-j}})$ \ which is at most 1-dimensional because  \ $\A\big/\A.P_j \simeq E\big/F_j$ \ is a \ $[\lambda]-$primitive theme (so has an unique corank 1 submodule which is \ $F_{k-1}\big/F_j$.). This kernel is zero if \ $\lambda_k-\lambda_{k-j} < 0$ \ and is given by \ $\C.b^{\lambda_k-\lambda_{k-j}}.e_{\lambda_{k-j}}$ \ in the other case. So we conclude that there always exists \ $\rho \in \C$ \ such that \ $x - \rho.b^{\lambda_k-\lambda_{k-j}}.e_{k-j}$ \ lies in \ $F_{k-j-1}$. $\hfill \blacksquare$

\begin{defn}\label{base standard 2}
In the situation of the previous proposition \ref{base standard 1} we shall call the \ $\C[[b]]-$
base \ $e_1, \dots, e_k$ \ of \ $E$ \ a  {\bf standard base} associated to the standard form \ $E \simeq \A\big/\A.P$. The element \ $e_k$ \  will be call a {\bf standard generator} for \ $E$.
\end{defn}

\parag{Remarks} 
\begin{enumerate}
\item In this situation \ $e_k$ \ generates \ $E$ \ as a left  \ $\A-$module and \ $\A.P$ \ is its annihilator ideal. The standard generator \ $e_k$ \  is uniquely defined modulo \ $F_{k-1}$, up to a non zero multiplicative constant, because its image in the quotient \ $E\big/F_{k-1} \simeq E_{\lambda_k}$ \ has to be a generator \ $e_{\lambda_k}$ \ of \ $E_{\lambda_k}$ \ such that \ $a.e_{\lambda_k} = \lambda_k.b.e_{\lambda_k}$.\\
 More generally, \ $e_j$ \ generates \ $F_j$ \ as a left  \ $\A-$module and \ $\A.P_{k-j}$ \ is its annihilator ideal. 
\item Let fix \ $\varphi \in \Xi^{(k-1)}_{\lambda} \setminus \Xi^{(k-2)}_{\lambda} $ \ and put \ $E : = \A.\varphi$. Then the unique normal rank \ $j$ \ submodule of the  \ $[\lambda]-$primitive theme  \ $E$ \ is given by the intersection \ $E \cap \Xi^{(j-1)}_{\lambda}$ \ for any \ $j \in [0,k]$;  see proposition \ref{theme 4}.\\
Note that any  \ $[\lambda]-$primitive rank \ $k$ \ theme  \ $E$ \ is isomorphic to such an "example", and in such a realization, a standard generator \ $e_k$ \ will be of the form
$$ e_k = c.s^{\lambda_k-1}.(Log\,s)^{k-1} \quad {\rm modulo} \  \Xi_{\lambda}^{(k-2)} $$
where \ $c \in \C^*$.
\end{enumerate}

The standard form for a given isomorphism class of a \ $[\lambda]-$primitive theme is not enough precise to have a chance to be unique : as it is shown in the following example the condition on the degrees of the polynomials \ $S_1, \dots, S_{k-1}$ \ are compatible with many different choices in order to give isomorphic themes.

\parag{Example} Let \ $\lambda_1 > 1$ \ be a rational number and \ $p_1 \geq 2$ \ and integer. Put \ $S_1 : = 1 + \sum_{i=1}^{p_1-1} \ x_i.b^i + \alpha.b^{p_1}$ \ where \ $\alpha \in \C^*$ \ and where \ $x_1, \dots, x_{p_1-1}$ \ are complex numbers. We shall show that \ $E : = \A\big/\A.P$ \ where 
 $$P : = (a - \lambda_1.b).S_1^{-1}.(a - (\lambda_1+p_1-1).b) $$
 is isomorphic to \ $\A\big/\A.(a - \lambda_1.b).(1 + \alpha.b^{p_1})^{-1}.(a - (\lambda_1+p_1-1).b) $.\\
 This clearly shows that the coefficients \ $x_1, \dots,x_{p_1-1}$ \ are irrelevant in the determination of the isomorphism class of \ $E : = \A\big/\A.P$.\\
 Remark first that \ $E \simeq \C[[b]].e_2 \oplus \C[[b]].e_1$ \ as a \ $\C[[b]]-$module, with \ $a$ define by the relations
 $$ a.e_2 = (\lambda_1 +p_1-1).b.e_2 + S_1.e_1 \quad {\rm and} \quad a.e_1 = \lambda_1.b.e_1 .$$
 We look for a \ $V \in \C[[b]]$ \ in order that \ $\tilde{e}_2 : = e_2 + V(b).e_1$ \ satisfies the relation
 $$ a.\tilde{e}_2 = (\lambda_1 +p_1-1).b.\tilde{e}_2 + (1 +\alpha.b^{p_1}).e_1 .$$
 It is clear that such an \ $\tilde{e}_2$ \ will give an isomorphism from \ $E$ \ to  
  $$\A\big/\A.(a - \lambda_1.b).(1 + \alpha.b^{p_1})^{-1}.(a - \lambda_1+p_1-1) .$$
 Let us compute \ $a.\tilde{e}_2$ :
 \begin{align*}
 & a.\tilde{e}_2 = (\lambda_1 +p_1-1).b.e_2 + S_1.e_1 + V.a.e_1 + b^2.V'.e_1 \\
 & \qquad  = (\lambda_1 +p_1-1).b.(\tilde{e}_2 - V.e_1) + S_1.e_1 + \lambda_1.b.V.e_1 + b^2.V'.e_1
 \end{align*}
 so we want to solve the differential equation
 $$ b^2.V' - (p_1-1).b.V = 1 + \alpha.b^{p_1} - S_1 = b.(x_1 + x_2.b + \dots+ x_{p_1-1}.b^{p_1-2}) $$
 with \ $V \in \C[[b]]$. This is easy because after simplification by \ $b$ \ we see that there is no term in\ $b^{p-1}$ \ in the right handside. $\hfill \blacksquare$\\

\bigskip

In what follows we shall fix \ $k \in \mathbb{N}$ \ and \ $\lambda_1, p_1, \dots, p_{k-1}$ \ assuming that \ $\lambda_1 > k-1$ \ is rational and that \ $p_1, \dots, p_{k-1}$ \ are natural integers. Then we define for each \ $j \in [1,k-1]$ \ the complex vector space \ $V_j \subset \C[[b]]$ \ in the following way :
\begin{enumerate}
\item If \ $p_j+\dots +p_{k-1} < k-j $ \ put \ $V_j : = \oplus_{i=0}^{k-j-1} \C.b^i $.
\item If  \ \ $p_j+\dots +p_{k-1} \geq  k-j $ \ put \ $q_j : = p_j + \dots + p_{j+h}$ \ where \ $h \in \mathbb{N}$ \ is minimal with the property that \ $q_j \geq k-j$ \ and put \ $ V_j : = \oplus_{i=0}^{k-j-1} \C.b^i \oplus \C.b^{q_j}$.
\end{enumerate}

\parag{Remarks}
\begin{enumerate}
\item The definition  of the vector spaces  \ $V_j$ \ does not depend on \ $\lambda_1$.
\item The definition of  the vector spaces \ $V_j$ \ depend only on \ $p_j, \dots, p_{k-1}$.
\item The vector space \ $V_j$ \ contains always  \ $\C.b^{p_j}$.
\end{enumerate} 

\begin{prop}\label{canonical form 1}
Let \ $E \simeq \A\big/\A.P$ \ a standard form for a \ $[\lambda]-$primitive theme. Then we have for each \ $j \in [1,k-1]$ \ the equality
\begin{equation*}
E_{\lambda_j} = P_j.E_{\lambda_j} \oplus V_j.e_{\lambda_j} \tag{@@@}
\end{equation*}
where \ $e_{\lambda_j}$ \ is a generator of \ $E_{\lambda_j}$ \ such that \ $a.e_{\lambda_j} = \lambda_j.b.e_{\lambda_j}$.
\end{prop}

\parag{proof} First notice that \ $E\big/F_j \simeq \A\big/\A.P_j$ \ is a \ $[\lambda]-$primitive theme of rank \ $k-j$. So, thanks to the lemma \ref{Ext} in the Appendix, we have
$$ \dim_{\C}\big(Ext^1_{\A}(E\big/F_j, E_{\lambda_j})\big) - \dim_{\C}\big(Ext^0_{\A}(E\big/F_j, E_{\lambda_j})\big) = k-j .$$
But \ $Ext^0_{\A}(E\big/F_j, E_{\lambda_j})$ \ has dimension \ $\leq 1$ \ and is \ $0$ \ if \ $\lambda_j > \lambda_k$ \ because such a map must factorize through an injection of the unique rank 1 quotient \ $E_{\lambda_k}$ \ of \ the theme \ $E\big/F_j$. So we have \ $ \dim_{\C}\big(Ext^1_{\A}(E\big/F_j, E_{\lambda_j})\big) = k-j+1 \quad {\rm or} \quad = k-j$ \ depending on the fact that  \ $\lambda_j \leq \lambda_k$ \ or \ $\lambda_j > \lambda_k$. Note that this corresponds to the inequalities \ $p_j + \dots+ p_{k-1} \geq k-j $ \ and \ $p_j + \dots + p_{k-1} < k-j $.\\
The exact sequence of left  \ $\A-$modules
$$ 0 \to \A \overset{.P_j}{\longrightarrow} \A \to E\big/F_j \to 0 $$
shows that the vector space \ $Ext^1_{\A}(E\big/F_j, E_{\lambda_j})$ \ is isomorphic to the cokernel of \ $P_j$ \ acting on \ $E_{\lambda_j}$. So the codimension of \ $P_j.E_{\lambda_j}$ \ in \ $E_{\lambda_j}$ \ is \ $k-j+1$ \ when \ $p_j+ \dots + p_{k-1} \geq k-j $ \ and \ $k-j$ \ when \ $p_j+ \dots + p_{k-1} <  k-j $.
But remark that in any case \ $P_j.E_{\lambda_j} \subset b^{k-j-1}.E_{\lambda_j}$. This implies  the result in the case 1. To prove the result in the case 2 it is enough to show that any linear combination
$$ \sum_{i=0}^{k-j-1} \ c_i.b^i.e_{\lambda_j} + \gamma.b^{q_j}.e_{\lambda_j} $$
which is in \ $P_j.E_{\lambda_j}$ \ is zero. The inclusion  \ $P_j.E_{\lambda_j} \subset b^{k-j-1}.E_{\lambda_j}$ \ implies that it is enough to prove that \ $b^{q_j}.e_{\lambda_j}$ \ is not in \ $P_j.E_{\lambda_j}$.\\
For this purpose remark first that if \ $x \in E_{\lambda_j} $ \ has a \ $b-$adic valuation equal to \ $q$ \ the \ $b-$adic valuation of \ $P_j.x$ \ will be exactely \ $q + k-j$ \ when \ $q$ \ is not one of the numbers \ $p_j + \dots+ p_{j+h} - (k-j)$ \ for an integer  \ $h \in [0,k-j-1]$. This comes from the fact that we may ignore the invertible elements \ $S_{j+1} \dots S_{k-1}$ \ in the computation of the \ $b-$adic valuation. Moreover after  the action of \ $(a - \lambda_{j+h+1}.b)\dots (a - \lambda_k)$ \ either we have the valuation equal to \ $q+k-(j+h)$ \ or the final valuation cannot be \ $q + k - j$. Now the action of \ $(a - \lambda_{j+h}.b) $ \ on \ $b^{q+k-(j+h)}.e_{\lambda_j}$ \ gives
$$ (q + \lambda_j + k -(j+h) - \lambda_{j+h}).b^{q+k-(j+h)+1}.e_{\lambda_j} $$
and we have
\begin{align*}
& (q + \lambda_j + k -(j+h) - \lambda_{j+h}) = q - [p_j+ \dots +p_{j+h-1} - h] + k - (j+h) \\
& \qquad   \quad = q - [p_j+ \dots +p_{j+h-1} - h] + k  - j
\end{align*}
which is not zero as long as \ $q$ \ is different from \ $p_j+ \dots +p_{j+h-1} - (k-j)$.\\
Assume now that \ $P_j.x = b^{q_j}.e_{\lambda_j}$ \ for some \ $x \in E_{\lambda_j}$. Let \ $q$ \ the \ $b-$adic valuation of \ $x$. If \ $q$ \ is not one of the numbers \ $p_j + \dots+ p_{j+h} - (k-j)$ \ for any integer  \ $h \in [0,k-j-1]$, then we must have \ $q_j = q + k - j$, which contradicts the definition of \ $q_j$.\\
If we have an integer \ $h \in [0,k-j-1]$ \ such that \ $q = p_j + \dots+ p_{j+h} - (k-j)$ \ then \ $q \geq q_j -(k-j)$. If we have \ $q = q_j -(k-j)$ \ what we say above implies that the valuation of \ $P_j.x$ \ is not \ $q_j$ \ and if \ $q > q_j -(k-j)$ \ the valuation of \ $P_j.x$ \ is strictely bigger than \ $q_j$. So such an \ $x$ \ cannot exists. $\hfill \blacksquare$.

\parag{Remark}
Note that \ $P_j$ \ depends only on \ $\lambda_j, p_j, \dots, p_{k-1}$ \ and \ $S_{j+1}, \dots, S_{k-1}$.

\bigskip

\begin{defn}\label{canonical form 2}
Let \ $E \simeq \A\big/\A.P$ \ a standard form for a \ $[\lambda]-$primitive theme. We shall say that it is a {\bf canonical form}  for  \ $E$ \ when we have \ $S_j \in V_j$ \ for each \ $j \in [1,k-1]$.
\end{defn}

\begin{thm}\label{canonical form 3}
Any  \ $[\lambda]-$primitive theme admits a canonical form.
\end{thm}

\parag{proof} We shall prove by descending induction on \ $j \in [0,k]$ \ that \ $E\big/F_j$ \ admits a canonical form. The case \ $j = k$ \ and \ $j = k-1$ \ are obvious, so assume that \ $E\big/F_{j+1}$ \ has a canonical form given by the isomorphism \ $\A\big/\A.P_{j+1} \simeq E\big/F_{j+1}$. Let \ $e \in E$ \ be an element inducing the image of \ $1$ \ in the quotient \ $E\big/F_{j+1}$. Then \ $P_{j+1}.e$ \ lies in \ $F_{j+1}$. Denote by \ $\pi : F_{j+1} \to F_{j+1}\big/F_j \simeq E_{\lambda_{j+1}}$ \ the quotient map and put \ $\pi(P_{j+1}.e) = T.e_{\lambda_{j+1}}$ \ where \ $T \in \C[[b]]$. As \ $P_{j+1}.e$ \ generates \ $F_{j+1}$, \ $T$ \ is invertible in \ $E_{\lambda_{j+1}}$ \ and, up to change the choice of \ $e_{\lambda_{j+1}}$ \ by a non zero constant factor, we may assume that \ $T(0) = 1$.\\
Write \ $T.e_{\lambda_{j+1}} = S_{j+1}.e_{\lambda_{j+1}} + P_{j+1}.z$ \ where \ $S_{j+1} $ \ lies in \ $V_{j+1}$. This possible thanks to the proposition \ref{canonical form 1}. Note that we have \ $S_{j+1}(0) = T(0) = 1$. Choose now \ $\varepsilon \in F_{j+1}$ \ satisfying \ $\pi(\varepsilon) = z$. Now \ $\tilde{e} : = e - \varepsilon$ \ is a generator of \ $E$\footnote{it induces the same class than \ $e$ \ in \ $E\big/(a.E+b.E)$.} and we have \ $S_{j+1}^{-1}.P_{j+1}.\tilde{e}$ \ which is send to \ $e_{\lambda_{j+1}}$ \ via \ $\pi$. So \ $(a - \lambda_{j+1}.b).S_{j+1}^{-1}.P_{j+1}.\tilde{e} = P_j.\tilde{e}$ \ lies in \ $F_j$. So \ $E\big/F_j$ \ admits a canonical form. $\hfill \blacksquare$

\section{The uniqueness of a canonical form}

To understand the classification of \ $[\lambda]-$primitive themes up to isomorphism, a key point is the uniqueness of the canonical form. This study will lead us to the notion of "invariant" theme, which is an interesting subclass for which this uniqueness property is true. 

\subsection{Endomorphisms of a theme.}

We shall begin by the study of \ $\A-$linear  injections between two \ $[\lambda]-$primitive themes.

\begin{lemma}\label{inj.}
Let \ $E' \subset E$ \ be two \ $[\lambda]-$primitive themes of the same rank \ $k$. Let \ $\mu_1, \dots, \mu_k$ \ and \ $\lambda_1, \dots, \lambda_k$ \ be their respective fundamental invariants. Then we have
\begin{enumerate}
\item For all \ $j \in [1,k], \mu_j \geq \lambda_j $.
\item The dimension of the complex vector space \ $E\big/E'$ \ is \ $\sum_{j=1}^k \ \mu_j - \lambda_j $.
\end{enumerate}
\end{lemma}

\parag{Proof} We shall make an induction on the rank \ $k$. The rank 1 case is clear, so assume that the lemma is proved in rank \ $\leq k-1$. As \ $F'_{k-1} \subset F_{k-1}$ \ we have immediately the inequalities \ $\mu_j \geq \lambda_j$ \ for \ $j \in [1,k-1]$. Moreover, the restriction to \ $E'$ \ of the quotient map \ $ \pi_k : E \to E\big/F_{k-1} \simeq E_{\lambda_k}$ \ is not zero, because \ $E'$ \ is not contained in \ $F_{k-1}$. So \ $\pi_k(E')$ \ is a rank 1 quotient of \ $E'$ \ contained in \ $E_{\lambda_k}$. This gives \ $\mu_k \geq \lambda_k$ \ as the unique rank 1 quotient of \ $E'$ \ is \ $E_{\mu_k}$.\\
Remark that the kernel of the restriction to \ $E'$ \ of \ $\pi_k$ \ is then \ $F'_{k-1}$ \ so we have  \ $F_{k-1}\cap E' = F'_{k-1}$ \ and this gives  the exact sequence of vector spaces
$$ 0 \to F_{k-1}\big/F'_{k-1} \to E\big/E' \to E\big/(F_{k-1} + E') \to 0 .$$
Now the induction hypothesis gives \ $\dim_{\C}(F_{k-1}\big/F'_{k-1}) = \sum_{j=1}^{k-1} \ \mu_j - \lambda_j$. Moreover the equatity \ $\pi_k(E') = E_{\mu_k}$ \ implies \ $\dim_{\C}(E\big/(F_{k-1} + E')) = \mu_k - \lambda_k$. We conclude using the exact sequence above.$\hfill \blacksquare$\\
 
\begin{thm}\label{hom. 1}
Let \ $E$ \ and \ $E'$ \ be two \ $[\lambda]-$primitive themes of the same rank \ $k$. The vector subspace of  \ $Hom_{\A}(E',E)$ \ of \ $\A-$linear maps with rank \ $\leq k-1$ \ has codimension \ $\leq 1$.\\
Assume now that the respective fundamental invariants \ $\mu_1, \dots, \mu_k$ \ and \ $\lambda_1, \dots, \lambda_k$ \ satisfy the inequalities
$$ \mu_j - \lambda_j \geq k-1 \quad \forall j \in [1,k] .$$
Then there exists an \ $\A-$linear injection \ $i : E' \to E$.
\end{thm}

\parag{Proof} Let us prove the first assertion by induction on the rank \ $k$. The case \ $k =1 $ \ is clear, so let assume the assertion proved for rank \ $k-1\geq 1$ \ and consider \ $\varphi_1$ \ and \ $\varphi_2$ \ two injections of \ $E'$ \ in \ $E$. Their restrictions to \ $F'_{k-1}$ \ are injective with value in \ $F_{k-1}$ \ so the induction hypothesis gives \ $\alpha \in \C^*$ \ such that \ $\varphi_1 - \alpha.\varphi_2$ \ is not injective on \ $F'_{k-1}$, so a fortiori on \ $E'$.\\
We shall show by induction on \ $k \geq 1$ \  the following  assertion which is a precise version of the second statement of the theorem :
\begin{itemize}
\item Let \ $E \simeq \A\big/\A.P$ \ and \ $E' \simeq \A\big/\A.P'$ \ two rank \ $k$ \ \ $[\lambda]-$primitive themes in standard forms, where
\begin{align*}
& P : = (a - \lambda_1.b).S_1^{-1} \dots S_{k-1}^{-1}.(a - \lambda_k.b) \\
& P' : = (a - \mu_1.b).T_1^{-1} \dots T_{k-1}^{-1}.(a - \mu_k.b)
\end{align*}
satisfying the condition \ $\mu_j - \lambda_j \geq k-1 \quad \forall j \in [1,k]$. Then exists an element
\begin{equation*}
 x  = \sigma.b^{\mu_k-\lambda_k}.e_k + \sum_{j=1}^{k-1} \ W_j.e_j  \quad {\rm where} \quad W_j \in b^{j-1}.\C[[b]] \quad \forall  j \in [1,k-1]  \tag{@}
 \end{equation*}
in \ $E$ \ such that \ $P'.x = 0$, where \ $e_1, \dots, e_k$ \ is the standard basis of \ $E$ \ associated to the standard form \ $E \simeq \A\big/\A.P$ \  and where \ $\sigma \in \C^*$.
\end{itemize}
Of course, sending the standard generator \ $e'_k$ \ of \ $E'$ \ to \ $x$ \ defines an \ $\A-$linear injection of \ $E'$ \ in \ $E$.\\
As the case \ $k = 1$ \ is immediate, we may assume our assertion proved in rank \ $k-1$. Using for \ $F'_{k-1}$ \ and \ $b.F_{k-1}$ \ the standard forms deduced from these of \ $E$ \ and \ $E'$ \ above, the induction hypothesis\footnote{Note that \ $\mu_j - \lambda_j \geq k-1$ \ implies \ $\mu_j - (\lambda_j+1)\geq k-2$ \ so the induction hypothesis applies to \ $F'_{k-1}$ \ and \ $b.F_{k-1}$. Note also that if \ $\tilde{P}_j$ \ is obtained from \ $P_j$ \ by \ $\lambda_h \mapsto \lambda_h + 1$ \ for each \ $h$, then \ $\tilde{P}_j.b.e_k = b.P_j.e_k$ \ as \ $(a - (\lambda+1).b).b = b.(a -\lambda.b)$.} gives an element
$$ y = \rho.b^{\mu_{k-1} - \lambda_{k-1} - 1}.b.e_{k-1} + \sum_{h=1}^{k-2} \ V_h.b.e_h \quad V_h \in b^{h-1}.\C[[b]] $$
in \ $b.F_{k-1}$ \ such that \ $Q'.y = 0$ \ where we put \ $P' = (a - \mu_1.b).T_1^{-1}.Q'$.\\
So we look for \ $x \in E$ \ of the form  \ $(@)$ \ and satisfying
\begin{enumerate}[i)]
\item \ $\sigma \not= 0 $.
\item \ $W_j \in b^{j-1}.\C[[b]] \quad \forall j \in [1,k-1]$.
\item \ $(a - \mu_k.b).x = T_{k-1}.y $.
\end{enumerate}
Remark that because the sequence \ $\lambda_j+j$ \ is increasing, we have
\begin{equation*}
 \mu_k \geq \lambda_k + k - 1 \geq \lambda_j + j - 1\quad \forall j \in [1,k] \tag{@@}
 \end{equation*}
 The relation iii) gives the equations
 \begin{align*}
& b.W'_{k-1} - (\mu_k - \lambda_{k-1}).W_{k-1} = \rho.T_{k-1}.b^{\mu_{k-1}-\lambda_{k-1} -1} - \sigma.S_{k-1}.b^{\mu_k-\lambda_k-1} \\
 & b^2.W'_{h} - (\mu_k-\lambda_h).b.W_h = T_{k-1}.b.V_h - S_h.W_{h+1} \tag{$@@@_h$}
 \end{align*}
The first equation will have a solution in \ $\C[[b]]$ \ unique modulo \ $\C.b^{\mu_k-\lambda_{k-1}}$ \ as soon as the coefficient of \ $b^{\mu_k- \lambda_{k-1}}$ \ is zero in the right handside. If \ $\alpha' \not=0$ \ is the coefficient of \ $b^{p'_{k-1}}$ \ in \ $T_{k-1}$ \ and \ $\alpha \not= 0$ \ the one of \ $b^{p_{k-1}}$ \ in  \ $S_{k-1}$, it is enough to choose \ $\sigma = \rho.\alpha/\alpha'$ \ to have the existence of \ $W_{k-1} \in b^{\mu_k - \lambda_{k-1}}.\C[[b]] \subset b^{k-2}.\C[[b]]$, unique modulo \ $\C.b^{\mu_k-\lambda_{k-1}}$ \ using the relations: 
$$ \mu_k - \lambda_{k-1} = \mu_{k-1} + p'_{k-1} - 1 - \lambda_{k-1} = \mu_k - \lambda_k + p_{k-1} - 1.$$
Assume that we have proved the existence of \ $W_{h+1} \in b^h.\C[[b]]$, unique modulo \ $\C.b^{\mu_k - \lambda_{h+1}}$ \ for some \ $h \geq 1$. To solve the equation \ $(@@@_h)$ \ it is enough to insure that the coefficient of \ $b^{\mu_k-\lambda_h+1}$ \ in \ $T_{k-1}.b.V_h - S_h.W_{h+1}$ \ is zero. But again we may choose arbitrarily the coefficient of \ $b^{\mu_k-\lambda_{h+1}}$ \ in \ $W_{h+1}$ \ thanks to our assumption, and because the coefficient of \ $b^{p_h}$ \ in \ $S_h$ \ is not zero.  We conclude easily because \ $\lambda_{h+1} = \lambda_h +p_h - 1$ \ gives \ $\mu_k - \lambda_h + 1 = \mu_k - \lambda_{h+1} + p_h$, and because \ $V_h \in b^{h-1}.\C[[b]]$ \ and \ $W_{h+1} \in b^h.\C[[b]]$ \ implies \ $W_h \in b^{h-1}.\C[[b]]$. $\hfill \blacksquare$

\parag{Remarks}
 \begin{enumerate}
 \item A immediate consequence of the previous proof is that, when there exists an injection from \ $E'$ \ to \ $E$, it sends a standard generator \ $e'$ \ of \ $E'$ \ to  \ $\tau.b^{\mu_k-\lambda_k}.e \quad {\rm modulo} \  F_{k-1}$ \ where \ $e$ \ is a standard generator of \ $E$ \ and \ $\tau \in \C^*$.
 \item In the final remark before the Appendix we  give two \ $[\lambda]-$primitive themes  of rank \ $3$ \ $E' : = E\big/F_1$ \ and \ $F_3$ \ such that \ $\mu_j - \lambda_j \geq k-2 = 1$ \ and such that there is no  \ $\A-$linear injection of \ $E'$ \ in \ $E$. 
   \end{enumerate}
  
  \begin{cor}\label{hom. 2}
  Let \ $E$ \ be a rank \ $k$  \ $[\lambda]-$primitive theme and let \ $R_j \subset Hom_{\A}(E,E)$ \ the vector subspace of endomorphisms of \ $E$ \ of  rank \ $\leq j$. For any \ $j \in [0,k-1]$, the vector space \ $R_{j+1}\big/R_j $ \ has dimension \ $\leq 1$. \\
  So we always have \ $\dim_{\C}(Hom_{\A}(E,E)) \leq k$ \ with equality if and only if \ $(R_j)_{j \in [0,k]}$ \ is a full flag\footnote{this means that each \ $R_{j+1}\big/R_j $ \ has dimension \ $1$ \ for \ $j \in [0,k-1]$.} of \ $Hom_{\A}(E,E)$.
  \end{cor}
  
  \parag{Proof} Let \ $\varphi$ \ be an endomorphism of rank \ $j$ \ of \ $E$. Then its kernel is \ $F_{k-j}$ \ as it is a normal rank \ $k-j$ \ submodule. Moreover the normalization of the image of \ $\varphi$ \ is \ $F_j$. So we have the following factorisation of \ $\varphi$ :
  $$ E \to E\big/F_{k-j} \overset{f}{\rightarrow} F_j \hookrightarrow  E.$$
  The correspondance \ $\varphi \mapsto f$ \ induces a linear bijection between \ $R_j\big/R_{j-1}$ \ and \ $Hom_{\A}(E\big/F_{k-j}, F_j)\big/G$ \ where \ $G$ \ is the subspace of morphism of rank at most  \ $ j-1$ \ in \ $Hom(E\big/F_{k-j}, F_j)$. The first assertion of the previous theorem gives that \ $\dim_{\C}(Hom_{\A}(E\big/F_{k-j}, F_j)\big/G) \leq 1$. $\hfill \blacksquare$
  
  \begin{cor}\label{hom. 3}
   Let \ $E$ \ be a rank \ $k$  \ $[\lambda]-$primitive theme. A sufficient condition for the existence of an \ $\A-$linear injection of \ $E\big/F_j$ \ in \ $F_{k-j}$ is that for each \ $h \in [1,k-j]$ \ we have
   $$ p_h + \dots + p_{h+j-1} \geq k-1 .$$
   For instance, if we have \ $p_h \geq k-1$ \ for each \ $h \in [1,k-1]$, then \ $Hom_{\A}(E,E)$ \ is \ $k-$dimensional.
   \end{cor}
   
   \parag{Proof} As \ $E\big/F_j$ \ is a rank \ $k-j$ \ $[\lambda]-$primitive theme with fundamental invariants \ $\lambda_{j+1}, \dots, \lambda_k$ \ and \ $F_{k-j}$ \ is also a rank \ $k-j$ \ $[\lambda]-$primitive theme with fundamental invariants \ $\lambda_1, \dots , \lambda_{k-j}$, we may conclude using the previous theorem as soon as we have
   $$ \lambda_{j+h} - \lambda_h \geq k-j-1 \quad \forall h \in [1,k-j] .$$
   But \ $\lambda_{j+h} - \lambda_h = p_h + \dots + p_{j+h-1} - j$ \ so the corollary is a consequence of theorem \ref{hom. 1}. The last assertion is now a consequence of the previous corollary \ref{hom. 2}.
   $\hfill \blacksquare$\\
   
   It is easy to see that the  necessary condition to have an \ $\A-$linear injection of \ $E\big/F_j$ \ into \ $F_{k-j}$ \ given by the lemma \ref{inj.} corresponds to the inequalities
   $$ p_h + \dots + p_{h+j-1} \geq j \quad \forall h \in [1,k-j].$$
   So it is trivially satisfied as soon as \ $p_h \geq 1$ \ for any \ $h \in [1,k]$, that is to say when the sequence \ $\lambda_1, \dots, \lambda_k$ \ is increasing (large).

   \subsection{Invariant \ $[\lambda]-$primitive themes.}
   
   Let us begin by two simples remarks.
   
   \parag{Remarks}
   \begin{enumerate}
   \item When \ $E$ \ is a \ $[\lambda]-$primitive theme the vector space \ $Hom_{\A}(E, E_{\mu})$ \ has dimension at most \ $1$ \ for any \ $\mu$. As \ $E$ \ has an unique rank 1 quotient (a,b)-module which is \ $E\big/F_{k-1} \simeq E_{\lambda_k}$, a non zero \ $\A-$linear map to \ $E_{\mu}$ \ factorizes via an \ $\A-$linear injection  of \ $E_{\lambda_k}$ \ to \ $E_{\mu}$. So \ $Hom_{\A}(E, E_{\mu})$ \ is zero if \ $\mu \not\in \lambda_k - \mathbb{N}$ \ and 1-dimensional if \ $\lambda_k = \mu + q, q \in \mathbb{N}$.
   \item Let \ $E_1$ \ and \ $E_2$ \ be \ $[\lambda]-$primitive themes with respective rank \ $k_1$ \ and \ $k_2$, and let \ $i : E_1 \hookrightarrow E_2$ \ be an \ $\A-$linear injection. Then \ $k_2 \geq k_1$ \ and the normalization\footnote{the smallest normal submodule containing \ $i(E_1)$.} of \ $i(E_1)$ \ is the normal submodule \ $G_{k_1}$ \ of rank \ $k_1$ \ in \ $E_2$. We may then apply the lemma \ref{inj.} to compute the codimension of the vector space \ $G_{k_1}\big/i(E_1)$.
   \end{enumerate}
   
   \begin{prop}\label{invariant 0}
   Let \ $E$ \ be a rank \ $k\geq 1$ \ $[\lambda]-$primitive theme\footnote{The proposition \ref{inv. primitive} will show that the existence of a rank \ $k-1$ \ endomorphism for any rank \ $k$ \ theme \ $E$ \  implies that \ $E$ \ is \ $[\lambda]-$primitive, for some \ $\lambda$.} \ and assume that \ $\varphi_0$ \ is an \ $\A-$linear endomorphism of \ $E$ \ of rank \ $k-1$. Then we have the following properties:
   \begin{enumerate}[i)]
   \item For each \ $j \in [0,k]$ \ the rank of \ $\varphi_0^j$ \ is \ $ k-j$.
   \item A \ $\C-$basis of \ $End_{\A}(E)$ \ is given by \ $\id, \varphi_0, \dots, \varphi_0^{k-1}$. So \ $End_{\A}(E)$ \ is \ $k-$dimensional. Moreover this algebra is commutative  and isomorphic to \ $\C[x]\big/(x^k)$.
   \item For each \ $j \in [1,k-1]$ \ the restriction of \ $\varphi_0$ \ to \ $F_j$ \ induces an  endomorphism  of  rank 
    \ $j-1$ \ of \ $F_j$. So we have a surjective restriction \ $End_{\A}(E) \to End_{\A}(F_j)$.
   \item For each \ $j \in [1,k-1]$ \ $\varphi_0$ \ induces an endomorphism of \ $E\big/F_j$ \ of rank \ $k-j-1$. So we have also a surjective map  \ $End_{\A}(E) \to End_{\A}(E\big/F_j).$
   \end{enumerate}
   \end{prop}
   
   \parag{Proof} We shall show i) by induction on \ $j$. As the cases \ $j=0$ \ and \ $j=1$ \ are clear, assume the assertion proved for some  \ $j \in [1,k-1]$ \ and let us show that \ $\varphi_0^{j+1}$ \ has rank \ $k-j-1$.\\
   The kernel of \ $\varphi_o^j$ \ is normal with rank \ $j$, so it is \ $F_j$. So \ $\varphi_0^j$ \ induces an \ $\A-$linear  injection of \ $E\big/F_j$ \ in \ $F_{k-j}$ \ whose image \ $\Phi_j$ \ has finite \ $\C-$codimension in \ $F_{k-j}$ \ (see remark 2 above). Applying again \ $\varphi_0$ \ whose kernel \ $F_1$ \  meets \ $\Phi_j$ \ in a rank 1 submodule, because \ $F_1 \subset F_{k-j}$. We see that the kernel of the restriction of \ $\varphi_0$ \ to \ $\Phi_j$, has rank 1 and so its image, which is \ $\Phi_{j+1}$ \ has rank \ $k-j-1$. So i) is proved.\\
   Assume that we have in \ $End_{\A}(E)$ \  a linear relation  \ $ \sum_{j=0}^{k-1} \ \alpha_j.\varphi_0^j = 0 $ \ where \ $\alpha_0, \dots, \alpha_{k-1}$ \ are complex numbers. Let \ $j_0$ \ the first integer such that \ $\alpha_{j_0} \not= 0$. Then we would have \ $\Phi_{j_0} \subset \sum_{h=j_0+1}^{k-1}\Phi_h$. But for each \ $j$ \ we know that \ $\Phi_j \subset F_{k-j}$. So we conclude that \ $\Phi_{j_0} \subset F_{k-j_0-1}$, and this contradicts the fact that \ $\varphi_0^{j_0}$ \ has rank \ $k - j_0$. So all \ $\alpha_j$ \ are zero and \ $\id, \varphi_0, \dots, \varphi_0^{k-1}$ \ are \ $k$ \ linearly independant vectors in \ $End_{\A}(E)$. But we have seen that the dimension of this vector space in \ $\leq k$. So it is \ $k-$dimensional and \ $\id, \varphi_0, \dots, \varphi_0^{k-1}$ \ is a basis. The algebra structure is then obvious because \ $\varphi_0^k = 0$.\\
   Let us proved iii). As the kernel of \ $\varphi_0$ \ is \ $F_1 \subset F_j$, the rank of the restriction of \ $\varphi_0$ \ to \ $F_j$ \ is \ $j-1$. So \ $\varphi_0(F_j) \subset F_{j-1} \subset F_j$ \ and \ $\varphi_0$ \ induces an endomorphism of \ $F_j$.  The surjectivity of the restriction map \ \ $End_{\A}(E) \to End_{\A}(F_j)$ \ is consequence of ii) applied to \ $F_j$.\\
   We have seen that \ $\varphi_0(F_j) \subset F_j$ \ so \ $\varphi_0$ \ induced an endomorphism \ $\tilde{\varphi}_0$ \ of \ $E\big/F_j$. The image \ $\Phi_1$ \ of \ $\varphi_0$ \ has finite codimension in \ $F_{k-1}$, so the quotient \ $\Phi_1\big/F_j \cap \Phi_1$ \ which is the image of \ $\tilde{\varphi}_0$, has rank \ $k - j - 1$ \ as a quotient of a rank \ $k-1$ \ $\C[[b]]-$module by a submodule of rank \ $j$. \\
   The surjectivity of the map  \ $End_{\A}(E) \to End_{\A}(E\big/F_j)$ \ is consequence of ii) applied to \ $E\big/F_j$. $\hfill \blacksquare$
   
   \parag{Main Example} Let us define  the monodromy operator \ $\mathcal{T} :  \Xi_{\lambda}^{(k-1)} \to \Xi_{\lambda}^{(k-1)}$ \ by the following rules :
   \begin{enumerate}
   \item It is \ $\A-$linear.
   \item \ $\mathcal{T}[s^{\lambda-1}.(Log\,s)^j] = e^{2i\pi.\lambda}.s^{\lambda-1}.(Log\,s + 2i\pi)^j$.
   \end{enumerate}
   Then it corresponds to the usual action of the monodromy of multivalued holomorphic functions on the punctured disc.\\
    Let \ $\varphi \in \Xi_{\lambda}^{(k-1)} \setminus \Xi_{\lambda}^{(k-2)}$ \ and assume that the theme \ $E : = \A.\varphi$ \ is invariant by the monodromy \ $\mathcal{T}$ \ of \ $\Xi_{\lambda}^{(k-1)}$. Then \ $\mathcal{T} - e^{2i\pi.\lambda}.\id$ \ induces an \ $\A-$endomorphism of rank \ $k-1$ \ of \ $E$ :\\
    To prove this we may assume\footnote{see proposition \ref{theme 4}.} that \ $\varphi = s^{\lambda_k-1}.(Log\,s)^{k-1} + \psi$ \ where \ $\psi $ \ is in \ $\Xi_{\lambda}^{(k-2)}$. Then we have to show that the degree in \ $Log\,s$ \ of \ $(\mathcal{T} - e^{2i\pi.\lambda}.\id)[\varphi] $ \ is precisely \ $k-2$. This is consequence of the fact that \ $\mathcal{T} - e^{2i\pi.\lambda}.\id$ \ sends \ $\Xi_{\lambda}^{(N-1)}$ \ in \ $\Xi_{\lambda}^{(N-2)}$ \ for any positive integer \ $N$ \ because we have
      $$ (\mathcal{T} - e^{2i\pi.\lambda}.\id)[\varphi] = e^{2i\pi.\lambda}.s^{\lambda_k-1}.(Log\,s + 2i\pi)^{k-1} - e^{2i\pi.\lambda}.s^{\lambda_k-1}.(Log\,s)^{k-1} \quad {\rm modulo} \ \Xi_{\lambda}^{(k-3)}.$$
      
\bigskip

\begin{defn}\label{invariant 1}
We shall say that a rank \ $k$ \ $[\lambda]-$primitive theme \ $E$ \ is {\bf invariant} \ when it has a endomorphism of rank \ $k-1$.
\end{defn}

With this definition we may reformulate the results of the proposition \ref{invariant 0} as follows.

\begin{cor}\label{invariant 2}
Let \ $E$ \ be an invariant rank \ $k$ \ $[\lambda]-$primitive theme. Then any normal \ $F$ \ of \ $E$ \ is invariant. The quotient \ $E\big/F$ \ is also invariant.
\end{cor}

\begin{prop}\label{invariant 3}
Let \ $E$ \ be a rank \ $k$ \ $[\lambda]-$primitive theme. The following properties are equivalent :
\begin{enumerate}[i)]
\item $E$ \ is invariant.
\item \ $\dim_{\C}(End_{\A}(E)) = k $.
\item The image of any \ $\A-$linear injection of \ $E$ \ in \ $\Xi_{\lambda}$ \ is independant of the choosen injection.
\item There exists an \ $\A-$linear injection \ $j : E \to \Xi_{\lambda}$ \ such that \ $j(E)$ \ is (globally) invariant by the monodromy \ $\mathcal{T}$ \ of \ $\Xi_{\lambda}$.
\end{enumerate}
\end{prop}

\parag{Proof} The implication \ $ i) \Rightarrow ii)$ \ is proved in the  proposition \ref{invariant 0}. The implication \ $ii) \Rightarrow iii)$ \ is consequence of the fact that if \ $i : E \to \Xi_{\lambda}$ \ is injective, then the composition with \ $i$ \ gives an injective \ $\C-$linear map \ $\tilde{i} : End_{\A}(E) \to Hom_{\A}(E, \Xi_{\lambda})$. As these two vector spaces have the same dimension \ $k$, the first by the assumtion ii), the second thanks to the theorem 2.2.1 of [B.05], we obtain iii) because any \ $\A-$linear map from \ $E$ \ to \ $\Xi$ \ has its image in \ $\Xi_{\lambda}^{(k-1)}$ \ (see the proposition \ref{theme 4}).\\
The implication \ $iii) \Rightarrow iv)$ \ is easy because \ $\mathcal{T}\circ i$ \ is an \ $\A-$linear injection of \ $E$ \ in \ $\Xi_{\lambda}^{(k-1)}$ \ when  \ $i$ \ is.\\
The implication \ $iv) \Rightarrow i)$ \ is consequence of the "main example" above. $\hfill \blacksquare$\\

\begin{cor} 
Let \ $E$ \ be a rank \ $k$ \ $[\lambda]-$primitive theme. A sufficient condition in order that \ $E$ \ is invariant is that we have \ $p_j \geq k-1\quad \forall j\in [1,k-1]$.
\end{cor}

\parag{Proof} This is a trivial consequence of the corollary \ref{hom. 3} and of ii) in the previous proposition. $\hfill \blacksquare$\\

\begin{cor}\label{invariant 4}
Let \ $E$ \ be an invariant  \ $[\lambda]-$primitive theme with fundamental invariants \ $\lambda_1, \dots, \lambda_k$. Then for \ $\delta > \lambda_k + k - 1$ \ a rational number, the theme \ $E^*\otimes_{a,b} E_{\delta}$ \ is invariant.
\end{cor}

\parag{Proof} It is enough to prove that the vector space \ $End_{\A}(E^*\otimes_{a,b} E_{\delta})$ \ is \ $k-$dimensional. But this vector space is isomorphic to \ $End_{\A}(E^*)$ \ because the tensor product by \ $E_{\delta}$ \ is just the change of \ $a$ \ in \ $a + \delta.b$, which does not change the \ $\A-$linear endomorphisms. But the transposition is an isomorphism \ $\C-$linear between \ $End_{\A}(E)$ \ and \ $End_{\A}(E^*)$. $\hfill \blacksquare$

\parag{Example} The \ $[\lambda]-$primitive themes with rank 2 are all invariant except those with fundamental invariants \ $\lambda_1, p_1 = 0$. To prove this it is enough to produce a rank 1 endomorphism when we assume \ $p_1 \geq 1$ \ and to show that in the case \ $p_1 = 0$ \ there is no such endomorphism. But an endomorphism of rank 1 is just an injection of \ $E\big/F_1$ \ in \ $F_1$, that is to say an \ $\A-$linear injection of \ $E_{\lambda_2}$ \ in \ $E_{\lambda_1}$. This exists and is unique up to a multiplicative non zero constant if and only if we have \ $\lambda_2 \geq \lambda_1$. But, by definition \ $\lambda_2 = \lambda_1 + p_1 - 1$. $\hfill \blacksquare$

\begin{lemma}\label{invariant 5}
Let \ $E$ \ be a rank \ $k \geq 2$  \ $[\lambda]-$primitive theme, and assume that either \ $p_{k-1} = 0 $ \ or \ $k \geq 3$ \ and  \ $p_{k-1} = 1, p_{k-2} \geq 2$. Then \ $E$ \ is not invariant.
\end{lemma}

\parag{Proof} Let \ $ E \simeq \A\big/\A.P$ \ a standard form for \ $E$ \ and let \ $e_1, \dots,e_k$ \ the corresponding standard basis of \ $E$. It is enough to show that there is no \\
 $x \in F_{k-1} \setminus F_{k-2}$ \ such that \ $P.x = 0$. Thanks to proposition \ref{base standard 1} we may write
\begin{equation*}
 x = \rho.b^{\lambda_k-\lambda_{k-1}}.e_{k-1} + \sum_{j=1}^{k-2} \ U_j.e_j \tag{*}
 \end{equation*}
where \ $\rho \in \C^*$ \ and \ $U_j \in \C[[b]]$ \ for \ $j \in [1,k-2]$. If \ $p_{k-1} = 0$, so \ $\lambda_k - \lambda_{k-1} = -1$, such an \ $x$ \ can not exist. So assume \ $k \geq 3$ \ and \ $p_{k-1} = 1$, so \ $\lambda_k = \lambda_{k-1}$. Then  \begin{equation*}
(a - \lambda_k.b).x = S_{k-1}.y \tag{**}
\end{equation*}
 with \ $y \in F_{k-2}$ \ and if we write \ $ P : = Q.S_{k-1}^{-1}.(a - \lambda_k.b)$ \ we shall have \ $Q.y = 0$ \ so, using again  proposition \ref{base standard 1} we have, for some \ $\sigma \in \C^*$
$$ y - \sigma.b^{\lambda_{k-1}-\lambda_{k-2}}.e_{k-2} \in F_{k-3}.$$
Now substituing \ $(^*)$ \ in \ $(^{**})$ \ gives that \ $U_{k-2}$ \ must satisfy
$$ S_{k-2} + b^2.U'_{k-2} - (\lambda_k - \lambda_{k-2}).b.U_{k-2} = \rho.S_{k-1}.b^{\lambda_{k-1}-\lambda_{k-2}} .$$
As we assume that \ $\lambda_{k-1} > \lambda_{k-2}$ \ this equation has no solution in \ $\C[[b]]$ \ because \ $S_{k-2}(0) = 1$. $\hfill \blacksquare$\\

This lemma has the interesting following corollary.

\begin{cor}\label{invariant 6}
Let \ $E$ \ be an invariant  \ $[\lambda]-$primitive theme. Then either the sequence \ $\lambda_1, \dots, \lambda_k$ \ is strictly increasing, or constant.
\end{cor}

\parag{Proof} The case \ $k = 2$ \ is clear thanks to lemma \ref{rang 2}. So let assume that \ $k \geq 3$. We want first to show that there exists \ $j_0 \in [1,k]$ \ such that we have
$$ \lambda_1 = \dots = \lambda_{j_0} < \lambda_{j_0+1} < \dots < \lambda_k .$$
Let \ $E$ \ be a rank \ $k$ \  invariant \ $[\lambda]-$primitive theme with fundamental invariants \ $\lambda_1, p_1, \dots, p_{k-1}$. Then each \ $p_j$ \ is positive because \ $p_j = 0$ \ would imply that \ $F_{j+1}\big/F_{j-1}$ \ is isomorphic to \ $E_{\lambda_j,\lambda_j}$ \ and so is not invariant. This would contradict the corollary \ref{invariant 2}.\\
If all \ $p_j$ \ are at least \ $2$ \ then \ $j_0 = 1$ \ satisfies our requirement. If this is not the case, let \ $j_0$ \ be the maximal integer \ such we have  $p_h = 1$ \ for each \ $h \in [1,j_0-1]$. Then we have \ $\lambda_1 = \dots = \lambda_{j_0} < \lambda_{j_0+1}$. So the invariant theme \ $E\big/F_{j_0-1}$ \ as fundamental invariants equal to \ $\mu_1 = \lambda_{j_0} < \mu_2 = \lambda_{j_0+1} \leq \dots \leq \mu_{k-j_0-1} = \lambda_k$. Assume that we have \ $\mu_{h-1} = \mu_h$ \ for some \ $h \geq 3$. Then the  lemma \ref{invariant 5} applied to the rank 3 invariant theme \ $F_{j_0+h+1}\big/F_{j_0+h-2}$ \ gives a contradiction.\\
We conclude using the remark  which follows the proof of the theorem \ref{tw. dual}, because for \ $\delta \in \mathbb{Q}$ \ large enough, \ $E^* \otimes_{a,b}E_{\delta}$ \ is a theme and it is invariant thanks to corollary \ref{invariant 4}. Then we have either \ $j_0 = 1$ \ or \ $j_0 = k$. $\hfill \blacksquare$\\

We shall say that an invariant \ $[\lambda]-$primitive theme is {\bf special}  when its fundamental invariants is of the form \ $\lambda_1, \lambda_1, \dots, \lambda_1$ \ which means that \ $p_1 = p_2 = \dots = p_{k-1} = 1$. In rank 2 all such fundamental invariants corresponds to an invariant theme. This no longer true in higher ranks. Will shall see in section 5.2.3  that in rank 3 the  \ $[\lambda]-$primitive themes
$$ \A\big/\A.(a - \lambda_1.b)(1 + \beta.b)^{-1}(a - \lambda_1.b)(1 + \alpha.b)^{-1}(a - \lambda_1.b)$$
are invariant for \ $\alpha = \beta \in \C^*$ \ but not invariant for \ $(\alpha,\beta) \in (\C^*)^2, \alpha \not= \beta$.

The following easy lemma is the first step of the determination of all special invariant  $[\lambda]-$primitive themes.

\begin{lemma}\label{invariant special}
Let \ $E$ \ be an  \ $[\lambda]-$primitive theme $E$ \ with fundamental invariants \ $\lambda_1, p_1 = \dots = p_{k-1} = 1$. Then \ $E$ \ is (special) invariant  if and only if \ $E\big/F_1$ \ is isomorphic to \ $F_{k-1}$. In this case, for any standard form \ $E \simeq \A\big/\A.P$ \ the coefficients of \ $b$ \ in \ $S_j$ \ is independent of \ $j \in [1,k-1]$ \ (and of the choice of the standard form).
\end{lemma}

\parag{Proof} By definition  \ $E$ \ is invariant if  there exists a rank \ $k-1$ \ endomorphism \ $\varphi_0$ \ of \ $E$. Its kernel is \ $F_1$ \ and its image is contains in \ $F_{k-1}$. In fact, the codimension of \ $\varphi_0(E)$ \ in \ $F_{k-1}$ \ is given by the lemma \ref{inj.}. But it gives zero and \ $\varphi_0$ \ induces an isomorphism between \ $E\big/F_1$ \ and \ $F_{k-1}$. The converse is obvious.\\
Now remark that, as \ $p_1 = \dots = p_{k-1} = 1$, the coefficient of \ $b$ \ in \ $S_j$ \ is the parameter of the rank \ $2$ \  \ $[\lambda]-$primitive theme \ $F_{j+1}\big/F_{j-1}$. So it only depends on the isomorphism class of \ $E$. So the point is to prove that these numbers are independant of \ $j$.\\
We shall prove this  by induction on the rank of \ $E$. \\
As \ $F_{k-1}$ \ is again an invariant special  \ $[\lambda]-$primitive theme, so it is enough to prove that the coefficient of \ $b$ \ in \ $S_{k-2}$ \ and in \ $S_{k-1}$ \ are the same.\\
The isomorphism between \ $E\big/F_1$ \ and \ $F_{k-1}$ \ induces an isomorphism between \ $E\big/F_{k-2}$ \ and \ $F_{k-1}\big/F_{k-3}$. An the coefficient of \ $b$ \ in \ $S_{k-1}$ \ and \ $S_{k-2}$ \ are respectively the parameters of these rank 2  \ $[\lambda]-$primitive themes which are isomorphic. $\hfill \blacksquare$\\

Note that conversely, it is obvious that the rank \ $k$  \ $[\lambda]-$ primitive theme given by \ $\A\big/\A.P$ \ with
$$ P : = (a - \lambda_1.b).S_1^{-1}.(a - \lambda_1.b).S_1^{-1} \dots S_1^{-1}.(a - \lambda_1.b) $$
where \ $S_1$ \ is in \ $\C[b]$ \ and satisfies \ $S_1(0) = 1$ \ and \ $S'(0) \not= 0$, is an invariant special theme because we have an obvious isomorphism between \ $E\big/F_1$ \ and \ $F_{k-1}$ \ which gives a rank \ $k-1$ \ endomorphism for \ $E$. \\
This implies  the existence of special invariant \ $[\lambda]-$ primitive themes with any rank.\\
But remark that in general \ $E$ \ is not given in a canonical  form\footnote{when \ $\deg(S_1) \geq 2$.}, and given such an \ $E$ \ in a canonical  form, it may be not so easy to recognize that it is invariant.

\subsection{Uniqueness of the canonical form.} 

The main result of this paragraph is the following uniqueness theorem.

\begin{thm}\label{unique 1}
Let \ $E$ \ be a rank \ $k \geq 2$ \ invariant \ $[\lambda]-$primitive theme. Then the canonical form \ $E \simeq \A\big/\A.P$ \ is unique.
\end{thm}

The proof will be an easy consequence of the next proposition.

\begin{prop}\label{unique 2}
Let \ $E$ \ be a rank \ $k \geq 2$ \ $[\lambda]-$primitive theme. Assume that the restriction map \ $End_{\A}(E) \to End_{\A}(E\big/F_1)$ \ is surjective. Then let
$$ P_1 : = (a - \lambda_2.b).S_2^{-1}\dots S_{k-1}^{-1}.(a - \lambda_k.b) $$
and \ $P : = (a - \lambda_1.b).S_1^{-1}.P_1$ \ where \ $E \simeq \A\big/\A.P$ \ is a standard form. Let \ $e = e_k$ \ be the standard generator of \ $E$, and let \ $e' \in E$ \ such that
\begin{enumerate}[i)]
\item \ $P_1.e' = T_1.e_1$ \ where \ $e_1$ \ is the standard generator of \ $F_1$.
\item \ $e - e' \in F_{k-1}$.
\end{enumerate}
Then \ $(T_1-S_1).e_1 $ \ is in \ $P_1.F_1$.
\end{prop}

\parag{Proof} The images \ $[e]$ \ and \ $[e']$ \ of  \ $e$ \ and \ $e'$ \ in \ $E\big/F_1$ \ are generators with the same annihilator ideal \ $\A.P_1$. As \ $[e-{e'}]$ \ lies in \ $F_{k-1}\big/F_1$ \ the endomorphism \ $\psi$ \  of \ $E\big/F_1$ \ defined by \ $\psi([e]) = [e] - [e']$ \ has rank \ $\leq k-2$. Our assumption gives then an endomorphism \ $\varphi$ \ of \ $E$ \ which induces \ $\psi$. So it is not an isomorphism. Put \ $\varepsilon : = \varphi(e)$. We have \ $P.\varepsilon = 0$, and because \ $\varepsilon $ \ lies in \ $F_{k-1}$ \ (the rank of \ $\varphi$ \ is \ $\leq k-1$), we have \ $P_1.\varepsilon = 0$ \ thanks to proposition \ref{base standard 1}. But the fact that \ $\psi$ \ induces \ $\varphi$ \ implies the relation
$$ \varepsilon = e - e' + U.e_1 \quad {\rm with} \quad U \in \C[[b]].$$
Then \ $P_1.\varepsilon = 0 = S_1.e_1 - T_1.e_1 + P_1.U.e_1 $. $\hfill \blacksquare$

\parag{Proof of the theorem \ref{unique 1}} By induction on the rank \ $k$ \ of \ $E$ \ we may assume that the canonical form for \ $E\big/F_1$ \ is unique. So we have uniqueness for \ $P_1$. \\
As we assume that \ $E$ \ is invariant, the assumption of surjectivity for the restriction map \ $End_{\A}(E) \to End_{\A}(E\big/F_1)$ \ is satisfied thanks to proposition \ref{invariant 0}. Then the previous proposition gives the uniqueness of \ $S_1.e_{\lambda_1}$ \ modulo \ $P_1.F_1$ \ which implies the uniqueness of the canonical form for \ $E$. $\hfill \blacksquare$

\begin{defn}\label{U 1}
We shall say that a \ $[\lambda]-$primitive theme has the {\bf property U} when its canonical form is unique.
\end{defn}

The reader will find in the Appendix some examples of non invariant \ $[\lambda]-$primitive themes which have the property U and also example where the property U is not satisfied.\\
The  problem of the characterization of \ $[\lambda]-$primitive themes having the property U is rather tricky.\\
Let us begin by a result of  non uniqueness.

\begin{lemma}\label{U 2}
Let \ $E$ \ be a rank \ $k \geq 3$ \ $[\lambda]-$primitive theme which is not invariant, but such that \ $E\big/F_1$ \ is invariant. Let \ $E \simeq \A\big/\A.P$ \ a standard form for \ $E$ with standard generator \ $e = e_k$. Then there exists a generator \ $e'$ \ of \ $E$ such that its annihilator is equal to \ $(a - \lambda_1.b).T_1^{-1}.P_1$ \ with \ $(S_1 - T_1).e_1 \not\in P_1.F_1$.
\end{lemma}

\parag{Remark} So in the situation of the previous lemma \ $E$ \ has not the property \ $U$.

\parag{Proof} Let \ $\psi$ \ an endomorphism of \ $E\big/F_1$ \ of rank \ $k-2$, and put \ $\psi([e]) = [\eta]$. Then  \ $\eta \in F_{k-1} \setminus F_{k-2}$ \ and the relation \ $P_1.[e] = 0$ \ in \ $E\big/F_1$ \ implies \ $P_1.\eta \in F_1$. Put \ $P_1.\eta = Z.e_1$ \ where \ $Z \in \C[[b]]$. If there exists \ $V \in \C[[b]]$ \ such that \ $Z.e_1 = P_1.V.e_1$, this would implies that \ $P_1(\eta - V.e_1) = 0$, and a fortiori, \ $P.(\eta - V.e_1) = 0$. But as \ $k \geq 3$, \ $\eta - V.e_1$ \ is in \ $F_{k-1}\setminus F_{k-2}$ \ and sending \ $e$ \ to \ $\eta - V.e_1$ \ is then an endomorphism of \ $E$ \ of rank \ $k-1$, contradicting our assumption that \ $E$ \ is not invariant. \\
So \ $Z.e_1 \not\in P_1.F_1$. But \ $P_1.\eta$ \ is in \ $b.E$ \ because \ $\eta \in F_{k-1} \subset a.E + b.E$ \ and \ $a.P_1.E \subset a^k.E + b.E \subset b.E$. So we have \ $Z(0) = 0$ \ as \ $F_1$ \ is normal. Define now \ $e' : = e - \eta$. This is generator of \ $E$ \ and it satisfies
$$ P_1.e' = T_1.e_1 \quad {\rm with} \quad  T_1 : = (S_1 - Z) .$$
We have \ $T_1(0) = 1$ \ and  \ $(S_1 - T_1).e_1$ \ is not in \ $P_1.F_1$. $\hfill \blacksquare$\\

The next proposition shows that there are not so many \ $[\lambda]-$primitive themes which satisfy the property U  and which are not invariant.

\begin{prop}\label{U 3}
Let \ $E$ \ be a a rank \ $k \geq 3$ \ non invariant  \ $[\lambda]-$primitive theme with the property U. Then for each \ $j \in [1,k-2]$ \ the theme \ $E\big/F_j$ \ is not invariant with the property U. This implies that \ $p_{k-1} = 0$.\\
Conversely, if \ $E\big/F_1$ \ satisfies the property U and \ $End_{\A}(E\big/F_1) = \C.\id$, then \ $E$ \ satisfies the property U.
\end{prop}

\parag{Proof} The fact that \ $E\big/F_j$ \ satisfies the property U when \ $E$ \ does  is obvious. If \ $E\big/F_1$ \ were invariant, the lemma \ref{U 2} would show that (we assume \ $k \geq 3$) \ $E$ \ is invariant. By induction on \ $j \in [1,k-2]$, we deduce that the quotients \ $E\big/F_j$ \ are not invariant and satisfies the property U. For \ $j = k-2$ \ we obtain \ $p_{k-1} = 0$.\\
To show the converse consider two generators \ $e, e'$ \ of \ $E$ \ such their images \ $[e],[e']$ \ in \ $E\big/F_1$ \ have the same annihilator ideal \ $\A.P_1$, the difference \ $e - e'$ \ is in \ $F_{k-1}$, and put \ $P_1.e = S_1.e_1, P_1.e' = T_1.e_1$ \ where \ $e_1$ \ is a generator of \ $F_1$ \ such that \ $a.e_1 = \lambda_1.b.e_1$. The endomorphism of \ $E\big/F_1$ \ which is defined by sending \ $[e]$ \ to \ $[e] - [e']$ \ is not surjective, because \ $[e] - [e']$ \ lies in \ $F_{k-1}\big/F_1$. So it is zero and we conclude that \ $e - e' \in F_1$. This implies that \ $(S_1 - T_1).e_1 \in P_1.F_1$ \ $\hfill \blacksquare$\\

An easy corollary describes completely the situation in rank \ $3$.

\begin{cor}\label{U 4}
The rank \ $3$ \ $[\lambda]-$primitive themes which satisfy the property U are the invariant themes and the themes such that \ $p_2 = 0$.
\end{cor}

\parag{Proof} As any rank \ $2$ \ $[\lambda]-$primitive theme has the property U, a non invariant rank \ $3$ \ theme with the property U must satisfy \ $p_2 = 0$ \ from the previous proposition. \\
Conversely, if \ $p_2 = 0$, then \ $E\big/F_1$ \ is isomorphic to \ $E_{\lambda_2,\lambda_2}$ \ and it satisfies the condition \ $End_{\A}(E\big/F_1) = \C.\id$, so the previous proposition implies it satisfies the property U.$\hfill \blacksquare$\\

We end this section by showing that if a rank \ $k$ \  theme \ $E$ \   has an endomorphism of rank \ $k-1$ \ this forces \ $E$ \ to be \ $[\lambda]-$primitive, for some \ $\lambda$.

\begin{prop}\label{inv. primitive}
Let \ $E$ \ be a rank \ $k$ \ theme which admits an endomorphism of rank \ $k-1$. Then \ $E$ \ is \ $[\lambda]-$primitive \ for some \ $\lambda \in \mathbb{Q}$.
\end{prop}

\parag{Proof} We shall prove the proposition by induction on the rank \ $k$. The case \ $k = 1$ \ is clear, so let us assume \ $k\geq 2$ \ and the result prove for  rank \ $k-1$.\\
 Let \ $\varphi$ \ an endomorphism of \ $E$ \  of rang \ $k-1$ \ and denote \ $F$ \ its image. Its kernel has  rank 1 so is isomorphic to \ $E_{\lambda}$ \ for some \ $\lambda \in \mathbb{Q}^{*+}$. As \ $F \simeq E \big/E_{\lambda}$ \ it is a rank \ $k-1$ \ theme.\\
  Assume first that \ $E_{\lambda} \cap F  = \{0\}$ \ then the restriction \ $\psi$ \ of \ $\varphi$ \ to \ $F$ \ is an  injective endomorphism of \ $F$. It is an isomorphism, because of the proposition below which is proved in [K.09] ; denote \ $\chi$ \ its inverse in \ $End_{\A}(F)$. Put \ $\theta : = \id_E - \chi\circ \varphi : E \to E$. We have \ $\varphi\circ\theta(x) = \varphi(x) - \varphi(\chi(\varphi(x))  = 0 $ \ for any \ $x \in E$. So \ $\theta$ \ takes values in \ $E_{\lambda}$. This implies that we have an (a,b)-linear isomorphism
  $$ E \simeq E_{\lambda} \oplus F \quad {\rm corresponding \ to} \quad \id_E = \theta + \chi\circ \varphi .$$
  So \ $E$ \ cannot be monogenic.\\
  Then we have \ $G : = E_{\lambda} \cap F \not= \{0\}$, and \ $G$ \  has rank 1. But \ $G$ \ is the kernel of the restriction \ $\psi : F \to F$ \ of \ $\varphi$. So \ $\psi$ \ has rank \ $k-2$ \ on \ $F$ \ which is a theme of rank \ $k-1$. The induction hypothesis gives then that \ $F$ \ is \ $[\mu]-$primitive, for some \ $\mu$. But \ $G = E_{\lambda+q}$ \ for some integer \ $q$ \ and \ $G \subset F$ \ shows that \ $\lambda \in [\mu]$ \ and  \ $F$ \ is \ $[\lambda]-$primitive. The exact sequence \ $0 \to E_{\lambda} \to E \to F \to  0$ \ implies now that \ $E$ \ is \ $[\lambda]-$primitive. $\hfill \blacksquare$
  
  \begin{prop}([K.09] proposition 2.26 .)
  Let \ $E$ \ be a regular (a,b)-module and \ $\varphi$ \ an (a,b)-endomorphism of \ $E$. If \ $\varphi$ \ is injective, it is an isomorphism.\\
  \end{prop}

\section{Holomorphic families of \ $[\lambda]-$primitive themes.}

\subsection{Definitions and first examples.}

\subsubsection{Definitions.}

Let \ $X$ \ be a  complex space. We shall denote \ $\mathcal{O}_X[[b]]$ \ the sheaf of \ $\C-$algebras on \ $X$ \ associated to the presheaf
$$ U \mapsto \mathcal{O}_X(U)[[b]] .$$
It is a sheaf of \ $\mathcal{O}_X-$algebras. For \ $\mathcal{I} \subset (\mathcal{O}_X)^p$ \ a subsheaf of \ $\mathcal{O}_X-$modules (resp. $\mathcal{O}_X-$coherent), we shall denote \ $\mathcal{I}[[b]]$ \ the subsheaf of \ $\mathcal{O}_X[[b]]-$modules (resp.  \ $\mathcal{O}_X[[b]]-$coherent) of  \ $(\mathcal{O}_X[[b]])^p$ \ which is generated by \ $\mathcal{I}$.\\

\begin{defn}\label{X-(a,b) 1}
Let \ $X$ \ be a complex space. A sheaf of \ $\mathcal{O}_X-(a,b)-$modules \ $\mathbb{E}$ \ on \ $X$ \ is a  locally free sheaf of finite type of \ $\mathcal{O}_X[[b]]-$modules endowed of a sheaf map
$$ a : \mathbb{E} \to \mathbb{E}  $$
which is \ $\mathcal{O}_X-$linear and continuous for the \ $b-$adic topology of \ $\mathbb{E} $ \ and satisfies the commutation relation \ $a.b - b.a = b^2$.\\
A morphism between two sheaves of \ $\mathcal{O}_X-(a,b)-$modules  is a morphism of sheaves of \ $\mathcal{O}_X[[b]]-$modules which commutes with the respective actions of \ $a$.
\end{defn}

\parag{Example} Let \ $\lambda \in ]0,1]$ \ be a  rational number and let \ $N$ \ be a natural integer. Define
$$ \Xi_{X,\lambda}^{(N)} : = \oplus_{j=0}^N \ \mathcal{O}_X[[b]].e_{\lambda,j} $$
and define \ $a$ \ by induction from the relations\footnote{the reader may think that \ $e_{\lambda,j} = s^{\lambda-1}.(Log\,s)^j\big/j! $.}
\begin{equation*}
  a.e_{\lambda,0} = \lambda.b.e_{\lambda,0} \quad {\rm and} \quad
  a.e_{\lambda,j} = \lambda.b.e_{\lambda,j} + b.e_{\lambda,j-1} \quad {\rm for} \ j \geq 1 
\end{equation*}

\parag{Notation} We shall denote, for \ $\Lambda$ \ a finite subset of \ $]0,1] \cap \mathbb{Q}$,
 \ $$\Xi_{X, \Lambda}^{(N)} : = \oplus_{\lambda \in \Lambda} \quad  \Xi^{(N)}_{X,\lambda}.$$
 An holomorphic map of a complex space \ $X$ \ in \ $\Xi_{\Lambda}^{(N)}$ \ will be, by definition,  a global section on \ $X$ \ of the sheaf \ $\Xi_{X, \Lambda}^{(N)}$.\\
 Remark that we have \ $a.\Xi_{X,\Lambda}^{(N)} \subset b.\Xi_{X, \Lambda}^{(N)}$ \ so they are "simple poles" \ $\mathcal{O}_X-(a,b)-$modules.

\bigskip

Let \ $x \in X$. We have an evaluation map \ $\mathcal{O}_X \to \mathcal{O}_X\big/\mathcal{M}_x \simeq \C_x$ \ where \ $\mathcal{M}_x \subset \mathcal{O}_X$ \ is the subsheaf of holomorphic germs which vanish at \ $x$. When \ $\mathbb{E}$ \ is a sheaf of \ $\mathcal{O}_X-(a,b)-$modules on \ $X$, we shall have, in an analoguous way, an evaluation map at \ $x$
$$ \mathbb{E} \to E(x) : = \mathbb{E}\big/\mathcal{M}_x[[b]].\mathbb{E}.$$
Then \ $E(x)$ \ is the fiber at \ $x$ \ of the sheaf \ $\mathbb{E}$. We shall consider a sheaf of \ $\mathcal{O}_X-(a,b)-$modules on \ $X$ \ as a family of (a,b)-modules parametrized by \ $X$.

\begin{defn}\label{holom. 0}
An holomorphic map \ $\varphi : X \to \Xi_{ \Lambda}^{(N)}$ \ is  {\bf $k-$thematic} when the following condition is satisfied :
\begin{itemize}
\item The \ $\mathcal{O}_X[[b]]-$submodule \ $\mathbb{E}_{\varphi}$ \ of \ $\Xi_{X, \Lambda}^{(N)}$ \ generated by the \ $a^{\nu}.\varphi, \nu \in \mathbb{N}$ \ is free of rank \ $k$ \ with basis \ $\varphi, a.\varphi, \dots, a^{k-1}.\varphi$.
\end{itemize}
\end{defn}

For each \ $x \in X$ \ we shall denote \ $E(\varphi(x))$ \ the rank \ $k$ \ theme given as 
$$\mathbb{E}_{\varphi}\big/\mathcal{M}_x[[b]].\mathbb{E}_{\varphi} \simeq \A.\varphi(x) \subset \Xi_{\Lambda}^{(N)}.$$

Note that in this definition we may always assume that \ $N = k-1$.

\begin{lemma}\label{holom. 1}
Let \ $X$ \ be a reduced complex space and \ $\varphi : X \to \Xi_{\Lambda}^{(N)}$ \ a \ $k-$thematic holomorphic map ; then the Bernstein polynomial \ $B_{\varphi(x)}$ \ of \ $E(\varphi(x))$ \ is locally constant on \ $X$. Moreover, if \ $\Lambda = \{\lambda\}$ \ the fundamental invariants of the \ $[\lambda]-$primitive themes \ $E(\varphi(x))$ \ are also  locally constant.
\end{lemma}

\parag{Proof} We may write, by assumption, 
$$ a^k.\varphi = \sum_{j=0}^{k-1} \ S_j.a^j.\varphi $$
where \ $S_1, \dots , S_{k-1}$ \ are global sections on \ $X$ \ of the sheaf \ $\mathcal{O}_X[[b]]$. As for each \ $x \in X$ \ the theme \ $\A.\varphi(x)$ \ has rank \ $k$, its Bernstein element is given by
$$ a^k - \sum_{j=0}^{k-1} \ \sigma_j(x).a^j  $$
where \ $\sigma_j(x)$ \ is the coefficient of \ $b^{k-j}$ \ in \ $S_j(x)$. Remark that, when it is not zero,  \ $\sigma_j(x).b^{k-j}$ \ is the initial form of \ $S_j(x)$. But the holomorphic function \ $x \to \sigma_j(x)$ \ takes only rational values, so it is locally constant.\\
In the \ $[\lambda]-$primitive case, the sequence \ $\lambda_j + j$ \ is increasing, so \ $\lambda_1, \dots, \lambda_k$ \ are determined by the Berstein element of \ $\A.\varphi(x)$. $\hfill \blacksquare$ 

\parag{Remark} Given an holomorphic map \ $\varphi : X \to \Xi_{\Lambda}^{(N)}$ \ it is not enough to check that for each \ $x \in X$ \ the (a,b)-module \ $E(\varphi(x))$ \ is a rank \ $k$ \ theme to have, even locally,  a \ $k-$thematic map, as it is shown by the following example :\\
Let \ $\lambda > 1$ \ a rational number and put for \ $z \in \C$ 
$$ \varphi(z) : = s^{\lambda-1}.Log\,s + (z + b).s^{\lambda-2} = s^{\lambda-1}.Log\,s + z.s^{\lambda-2} + \frac{1}{\lambda-1}.s^{\lambda-1} .$$
Then the Bernstein element of \ $\A.\varphi(z)$ \  is \ $(a - \lambda.b)(a - \lambda.b)$ \ for \ $z \not= 0$, but for \ $z = 0$ \ the Bernstein element of \ $\A.\varphi(0)$ \ is \ $(a - (\lambda+1).b)(a - \lambda.b)$. We conclude using  the previous lemma.\\

\parag{Example} Let \ $X$ \ be any reduced and irreducible complex space and let \\
 $\varphi : X \to \Xi_{\lambda}^{(k-1)}$ \ be an holomorphic map such that the coefficient of \ $e_{\lambda,k-1}$ \ is equal to \ $b^n.S$ \ where \ $S$ \ is an invertible element of the algebra \ $\mathcal{O}(X)[[b]]$\footnote{This is equivalent to say that the constant term (constant in \ $b$) of \ $S$ \ is invertible in \ $\mathcal{O}(X)$.}, and such that the \ $b-$valuation of \ $\varphi - b^n.S$ \ is strictely bigger than \ $n$. Then the sheaf \ $\mathbb{E}_{\varphi} : = \sum_{j=0}^{k-1} \ \mathcal{O}_X[[b]].a^i.\varphi$ \ is free of rank \ $k$ \ on \ $\mathcal{O}_X[[b]]$ \ and stable by \ $a$ :\\
It does not reduce the generality to assume that \ $S = 1$, and in this case \\ 
$\psi : = (a - (\lambda+n).b).\varphi$ \ satisfies the same hypothesis that \ $\varphi$ \ with \ $k$ \ replace by \ $k-1$ \ and \ $n$ \ replace by \ $n+1$. An easy induction allows to conclude.\\
Note that in this kind of example we have \ $p_1 = \dots = p_{k-1} =0 $. So this method constructs only very specific \ $[\lambda]-$primitive themes.

The reader will find a general and systematic  method to produce holomorphic \ $k-$thematic maps in the Appendix (see corollary \ref{exist. k-thm.}).

\begin{defn}\label{holom. 2}
Let \ $X$ \ be a reduced complex space and let \ $\mathbb{E}$ \ be a sheaf of \ $\mathcal{O}_X-(a,b)-$modules on \ $X$. We shall say that \ $\mathbb{E}$ \ is an holomorphic family of rank \ $k$ \  themes parametrized by \ $X$ \ when the following condition is satisfied :
\begin{itemize}
\item There exists an open covering \ $(\mathcal{U}_{\alpha})_{\alpha \in A}$ \ of \ $X$ \ and for each \ $\alpha \in A$ \ a finite set \ $\Lambda_{\alpha} \subset ]0,1] \cap \mathbb{Q}$ \  a \ $k-$thematic holomorphic map
$$ \varphi_{\alpha} : \mathcal{U}_{\alpha} \to \Xi_{\Lambda_{\alpha}}^{(k-1)} $$
with an isomorphism of  sheaves of \ $\mathcal{O}_X-(a,b)-$modules \ $\mathbb{E}_{\vert \mathcal{U}_{\alpha}} \simeq \mathbb{E}_{\varphi_{\alpha}}$ \ on \ $\mathcal{U}_{\alpha}$.
\end{itemize}
\end{defn}

\parag{Pull back} Of course when we have an holomorphic map \ $f : Y \to X$ \ between two reduced complex space and \ $\mathbb{E}$ \ an holomorphic family of \ $\Lambda-$primitive themes parametrized by \ $X$ \ we may take the pull back of this family by setting for \ $y \in Y$ \ 
 $E(y) : = E(f(y))$.  It is easy to see that the pull back of an holomorphic family is again an holomorphic family, the corresponding sheaf of \ $\mathcal{O}_Y-(a,b)-$modules on \ $Y$ \ being the "analytic" pull back \ $f^*(\mathbb{E}) : = \mathcal{O}_Y \otimes f^{-1}(\mathbb{E})$,  simply because when \ $\varphi : X \to \Xi_{\Lambda}^{(N)}$ \ is an holomorphic and  \ $k-$thematic map, the composition \ $\varphi\circ f : Y \to \Xi_{\Lambda}^{(N)}$ \ is again holomorphic and \ $k-$thematic.

\subsubsection{Examples}

\parag{Example : the rank 1 case.}
Let \ $X$ \ be a connected reduced complex space, and let \ $\varphi : X \to \Xi_{\lambda}^{(0)}$ \ be an holomorphic \ $1-$thematic map. Then there exists \ $S \in \Gamma(X, \mathcal{O}_X[[b]])$ \ which is invertible\footnote{Note that if for some \ $x \in X, S(x)$ \ is not invertible, the Bernstein element jumps, so the map is not \ $1-$thematic.} and an integer\ $n$ \ such that \ $\varphi = S.s^{\lambda+n-1}$. So the sheaf \ $\mathbb{E}_{\varphi}$ \ is isomorphic to the sheaf \ $\mathbb{E}_{\psi}$ \ where \ $\psi : X \to  \Xi_{\lambda}^{(0)}$ \ is the constant map with value \ $s^{\lambda+n-1}$.\\

We study the rank \ $2$ \ case in the next two results.

\begin{prop}\label{holom. rank 2}
Fix \ $\lambda_1 > 1 $ \ a rational number in \ $[\lambda] \in \C\big/\mathbb{Z}$ \ and a positive integer \ $p$. Let \ $X$ \ be a complex space, and let \ $\varphi : X \to \Xi_{\lambda}^{(1)}$ \ be an holomorphic \ $2-$thematic map, such that the fundamental invariants of the themes \ $E(\varphi(x), x \in X$ \ are given by \ $\lambda_1, p_1 = p \geq 1$. Then there exists an holomorphic map \ $\alpha : X \to \C^*$ \ such that we have
\begin{enumerate}
\item The map \ $\tilde{\psi} : X \to  \Xi_{\lambda}^{(1)}$ \ defined by
\begin{equation*}
\tilde{\psi}(x) : =  \alpha(x).s^{\lambda_1+p-2}.Log\,s + c(\lambda_1,p).s^{\lambda_1-1}   \tag{@}
 \end{equation*}
 where \ $c(\lambda_1,p) = \frac{-1}{p}.(\lambda_1-1).\lambda_1\dots (\lambda_1+p-2) $.
\item \ $\tilde{\psi}$ \  is \ $2-$thematic 
\item The sheaves of  \ $\mathcal{O}_X-(a,b)-$modules \ $\mathbb{E}_{\tilde{\psi}}$ \ and \ $  \mathbb{E}_{\varphi}$ \ co{\"i}ncide.
 \end{enumerate}
 \end{prop}
 
 The proof will use the following lemma.
 
 \begin{lemma}\label{rank 2}
 Fix \ $\lambda_1 > 1 $ \ a rational number in \ $[\lambda] \in \C\big/\mathbb{Z}$ \ and a positive integer \ $p$. If \ $x \in \C^*$ \ and \ $Z \in \C[[b]]$ \ are such that
 $$ \chi : = x.s^{\lambda_1+p-2}.Log\,s + Z.s^{\lambda_1-2} $$
 satisfies \ $(a -(\lambda_1+p-1).b).\chi = (1 + \alpha.b^p).s^{\lambda_1-1}$ \ for some given \ $\alpha \in \C^*$ \ then we have
 $$ x =  \alpha/\gamma \quad {\rm and} \quad Z = \frac{-1}{p} + z.b^p $$
 for some \ $z \in \C$, where \ $\gamma : = (\lambda_1-1). \lambda_1.(\lambda_1+1) \dots (\lambda_1+p-2)$.
 \end{lemma}
 
 The proof is an exercice left to the reader.\\
 
 \parag{Remark} As \ $\A.\chi$ \ contains \ $\C[[b]].s^{\lambda_1-1}$ \ it contains \ $\chi_z : = \chi -z.b^p.s^{\lambda_1-2}$ \ and both \ $\chi$ \ and \ $\chi_z$ \ have the same annihilator which is the ideal \\
  $$\A.(a - \lambda_1.b).(1 + \alpha.b^p)^{-1}.(a -(\lambda_1+p-1).b),$$
  we have for each \ $z \in \C$ \ an automorphism of \ $\A.\chi$ \ which sends \ $\chi$ \ to \ $\chi_z$.
 
 \parag{Proof of the proposition} Without lost of generality we may assume that
   \begin{equation*}
 \varphi(x) = s^{\lambda_1 + p -2}.Log\,s + \Sigma(x).s^{\lambda_1-2} \tag{1}
 \end{equation*}
 where \ $\Sigma \in \Gamma(X, \mathcal{O}_X[[b]])$ \  is invertible in \ $\C[[b]]$ \ for each \ $x \in X$. This is consequence of the fact that we must have for each \ $x \in X$ \ a rank \ $2$ \ theme with fundamental invariants \ $\lambda_1, p_1 = p \geq 1$. Then we have
 \begin{align*}
 & (a - (\lambda_1 + p-1).b).\varphi(x) = \frac{s^{\lambda_1+p-1}}{\lambda_1+p-1} + \Sigma(x).s^{\lambda_1-1} + \\
 & \qquad \qquad  + b^2.\Sigma'(x).s^{\lambda_1-2} - (\lambda_1+p-1).b.\Sigma(x).s^{\lambda_1-2} \\
 & \qquad = (b.\Sigma(x)' - p.\Sigma(x) + \gamma.b^{p}).\frac{s^{\lambda_1-1}}{\lambda_1-1}
 \end{align*}
 where \ $\gamma : = (\lambda_1-1)\lambda_1\dots (\lambda_1+p-2) $ \ and where \ $\Sigma(x)'$ \ denote the derivative in \ $b$ \ of \ $\Sigma(x) \in \C[[b]]$. Note that \ $\gamma$ \ gives the identity
 $$ \frac{s^{\lambda_1+p-1}}{\lambda_1+p-1} = \gamma.b^p.\frac{s^{\lambda_1-1}}{\lambda_1-1}.$$
 So we have \ $(a - (\lambda_1 + p-1).b).\varphi(x) = S(x).s^{\lambda_1-1}$ \ with 
  $$(\lambda_1 - 1).S(x) : = b.\Sigma(x)' - p.\Sigma(x) + \gamma.b^p.$$
 Remark that the constant term in \ $b$ \ in \ $S(x)$ \ is equal to \ $(-p/(\lambda_1-1))-$times the constant term of \ $\Sigma(x)$ \ which not zero. Put \ $ S(x) : = S_0(x) + S_p(x).b^p + b.\tilde{S}(x) $ \ where \ $\tilde{S}(x)$ \ has no term in \ $b^{p-1}$. From our previous remark we see that \ $S_0$ \ is an invertible holomorphic function on \ $X$.\\
 Now let \ $T \in \Gamma(X, \mathcal{O}_X[[b]])$ \ be a solution of the equation
 $$ b.T(x)' - (p-1).T(x) = \tilde{S}(x) \quad \forall x \in X.$$
 Such a \ $T$ \ exists  because \ $\tilde{S}$ \ has no term in \ $b^{p-1}$. Define now \ $\psi : X \to \Xi_{\lambda}^{(1)}$ \ as
 $$ \psi(x) : = \varphi(x) - T(x).s^{\lambda_1-1} .$$
 As \ $\mathbb{E}_{\varphi}$ \ contains \ $\mathcal{O}_X[[b]].s^{\lambda_1-1}$, thanks to the invertibility of \ $S$ \ in \ $\mathcal{O}_X[[b]]$, the equality \ $\mathbb{E}_{\psi} = \mathbb{E}_{\varphi}$ \ is clear. So \ $\psi$ \ is \ $2-$thematic. \\
 Now compute \ $(a - (\lambda_1 + p-1).b).\psi$:
 \begin{align*}
 & (a - (\lambda_1 + p-1).b).\psi(x) = (a - (\lambda_1 + p-1).b).(\varphi(x) - T(x).s^{\lambda_1-1}) \\
 & \qquad  = S(x).s^{\lambda_1-1} - \big[T(x).a + b^2.T(x)'\big].s^{\lambda_1-1} \\
 & \qquad =  (S_0(x) + S_p(x).b^p).s^{\lambda_1-1}
 \end{align*}
 This shows that \ $\tilde{\psi} : = S_0^{-1}.\psi$ \ satisfies the relation
 $$ (a - (\lambda_1 + p-1).b).\tilde{\psi}(x) = (1 + \alpha(x).b^p).s^{\lambda_1-1} $$
 for each \ $x \in X$ \ where \ $\alpha(x) : = S_0^{-1}(x).S_p(x)$. And, of course, \ $\mathbb{E}_{\tilde{\psi}} = \mathbb{E}_{\psi} = \mathbb{E}_{\varphi}$ \ so conditions  2. and 3.  are satisfied. The lemma shows that there exists an holomorphic function \ $z : X \to \C$ \ such that \ $\tilde{\psi}(x) - z(x).b^{p-1}.s^{\lambda_1-1}$ \ has the desired form. But sending \ $\tilde{\psi}$ \ to \ $\tilde{\psi} - z.b^{p-1}.s^{\lambda_1}$ \ is an automorphism of the \ $\mathcal{O}_X-(a,b)-$module \ $\mathbb{E}_{\varphi}$. $\hfill \blacksquare$\\
 
 The case \ $p_1 = p = 0$ \ is given by the following lemma ; as it is a simple variant of the previous proposition, we let the proof to  the reader.
 
 \begin{lemma}\label{rank 2 bis}
 Fix \ $\lambda_1 > 1 $ \ a rational number in \ $[\lambda] \in \C\big/\mathbb{Z}$. Let \ $X$ \ be a complex space, and let \ $\varphi : X \to \Xi_{\lambda}^{(1)}$ \ be an holomorphic \ $2-$thematic map, such that the fundamental invariants of the themes \ $E(\varphi(x), x \in X$ \ are given by \ $\lambda_1, p_1 = 0$. Then if \ $\psi : = (\lambda_1-1). s^{\lambda_1-2}.Log\,s $ \ there is an isomorphism  of \ $\mathcal{O}_X-(a,b)-$modules between the sheaf \ $\mathbb{E}_{\varphi}$ \ and the sheaf \ $\mathbb{E}_{\psi}$ \ corresponding to the constant map with value \ $\psi$.
 \end{lemma}
 
 Note that the sheaf \ $\mathbb{E}_{\psi}$ \ is given as \ $\mathcal{O}_X[[b]].e_1 \oplus \mathcal{O}_X[[b]].e_2$ \ with \ $a$ \ defined by
 $$ a.e_1 = \lambda_1.b.e_1 \quad {\rm and} \quad a.e_2 = (\lambda_1-1).b.e_2 + e_1 .$$
 
 \bigskip
 
 \begin{defn}\label{parametre 1}
 When \ $E$ \ is a rank \ $2$  \ $[\lambda]-$primitive theme with fundamental invariants \ $\lambda_1, p_1 \geq 1$ \ we shall call the {\bf parameter} of \ $E$ \ the  number\ $\alpha \in \C^*$ \ such that \ $E$ \ is isomorphic to \ $\A\big/\A.(a - \lambda_1.b)(1 + \alpha.b^p)^{-1}.(a - (\lambda_1+p_1-1).b) $.
 \end{defn}
 
 The following lemma shows that this number is quite easy to detect.
 
 \begin{lemma}\label{parametre 0}
 Let $E$ \ be a rank \ $2$  \ $[\lambda]-$primitive theme with fundamental invariants \ $\lambda_1, p_1 \geq 1$. Let \ $S \in \C[[b]]$ \ such that \ $S(0) = 1$ \ and such that 
  $$E \simeq  \A\big/\A.(a - \lambda_1.b).S^{-1}.(a - (\lambda_1+p_1-1).b).$$
   Then the coefficient of \ $b^{p_1}$ \ in \ $S$ \ is the parameter \ $\alpha$ \ of \ $E$.
 \end{lemma}
 
 \parag{Proof} Let \ $e$ \ be a generator of \ $E$ \ with annihilator the ideal 
  $$\A.(a - \lambda_1.b).S^{-1}.(a - (\lambda_1+p_1-1).b) ,$$
   and put \ $e_2 : = e$ \ and \ $e_1 : = S^{-1}.(a - (\lambda_1+p_1-1).b).e$. Then \ $e_1, e_2$ \ is a \ $\C[[b]]-$base of \ $E$ \ and so we look for a generator \ $\varepsilon$ \  of \ $E$, write
 $$ \varepsilon : = U.e_2 + V.e_1 $$ 
 with \ $U, V \in \C[[b]]$, which is annihilated by \ $(a - \lambda_1.b)(1 + \alpha.b^p)^{-1}.(a - (\lambda_1+p_1-1).b)$. Remark that we know " a priori" that such an \ $\varepsilon$ \ exists thanks to the proposition \ref{holom. rank 2} with \ $X = \{pt\}$. Then compute 
 \begin{align*}
 & (a - (\lambda_1+p_1-1).b).\varepsilon = b^2.U'.e_2 + U.S.e_1 + b^2.V' - (p_1-1).b.V.e_1.
 \end{align*}
 As we are in a theme the kernel of \ $a - \lambda_1.b$ \ is \ $F_1 = \C[[b]].e_1$. So this implies that \ $U'  = 0$ \ and \ $U = U(0)$. We obtain also the equation
 $$ U(0).S + b^2.V' - (p_1-1).b.V = \rho.(1 + \alpha.b^{p_1}) $$
 for some \ $\rho \in \C$. As we assume \ $S(0) = 1$ \ this implies \ $U(0) = \rho$. The fact that \ $\varepsilon$ \ is a generator of \ $E$ \ insure that \ $U(0) \not= 0$, and then
 $$ b^2.V' - (p_1-1).b.V = \rho(1 - S) + \rho.\alpha.b^p $$
 and this equation implies that the coefficient of \ $b^{p_1}$ \ in \ $S$ \ is equal to \ $\alpha$. $\hfill \blacksquare$
 
 \parag{Exercice} Let \ $E$ \ be a rank \ $2$  \ $[\lambda]-$primitive theme with fundamental invariants \ $\lambda_1, p_1= p  \geq 1$. Show that, for \ $\delta \in \mathbb{Q}$ \ such that \ $\lambda_1 + \delta >1$ \  the parameter of \ $E \otimes_{a,b} E_{\delta}$ \ is the same than the parameter of \ $E$.\\
 Deduce then that for \ $-\lambda_1 - p +  \delta + 1 > 1$ \ the parameter of the theme \ $E^* \otimes_{a,b} E_{\delta}$ \ is also  the same than the parameter of \ $E$.
 
 \subsubsection{Holomorphy criterion.}
 
 Let us go back to holomorphic families of \ $[\lambda]-$primitive themes of any rank.
 
 \begin{prop}\label{holom. 2}
 Let \ $X$ \ be a connected reduced complex space and let \ $\mathbb{E}$ \ be an holomorphic family of rank  \ $k$ \ $[\lambda]-$primitive themes parametrized by \ $X$. Denote \ $\lambda_1, p_1, \dots, p_{k-1}$ \ the corresponding  fundamental invariants. For each \ $j \in [0,k]$ \ there exists an unique  holomorphic family \ $\mathbb{F}_j$ \ of rank \ $j$ \ \ $[\lambda]-$primitive themes parametrized by \ $X$ \ with the following properties :
 \begin{enumerate}[i)]
 \item \ $\mathbb{F}_j \subset \mathbb{F}_{j+1}$, $\mathbb{F}_k = \mathbb{E}$ ;
 \item for each \ $x \in X$ \ the theme \ $F_j(x)$ \ is the normal submodule of rank \ $j$ \ of \ $E(x)$;
 \item the family \ $\mathbb{E}\big/\mathbb{F}_j$ \ is holomorphic for each \ $j \in [0,k]$.
 \end{enumerate}
 \end{prop}
 
 \parag{Proof} The statement is local on \ $X$ \ and we may assume that \ $\mathbb{E} = \mathbb{E}_{\varphi}$ \\
  where \ $\varphi : X \to \Xi_{\lambda}^{(k-1)}$ \ is a \ $k-$thematic holomorphic map. Up to change \ $\varphi$ \ by the action of an invertible section of the sheaf \ $\mathcal{O}_X[[b]]$ \ we may assume that \\
   $\varphi = s^{\lambda_k-1}.(Log\,s)^{k-1} + \theta $ \ where \ $\theta$ \ is a section of the sheaf \ 
    $\Xi_{X, \lambda}^{(k-2)}$.\\
     Put \ $\psi : = (a - \lambda_k.b).\varphi$. Then \ $\psi$ \ is \ $(k-1)-$thematic, because \ $\psi, a.\psi, \dots, a^{k-2}.\psi$ \ is \ $\mathcal{O}(X)[[b]]-$free and generates \ $\mathbb{E}_{\psi}$ :\\
 If we have \ $\sum_{j=0}^{k-2} \ U_j.a^j.\psi = 0 $ \ with \ $U_j \in \mathcal{O}(X)[[b]]$ \ this implies
 $$ \sum_{j=0}^{k-2} \ U_j.a^{j+1}.\varphi -  \lambda_k.\sum_{j=0}^{k-2} \ U_j.a^j.b.\varphi = 0. $$
 As, by hypothesis \ $\varphi, a.\varphi, \dots, a^{k-1}.\varphi$ \ is \ $\mathcal{O}(X)[[b]]-$free, this gives successively \ $U_{k-2} = 0, U_{k-3} = 0, \dots, U_1 = 0$.\\
 Then it is easy to see that \ $\mathbb{F}_{k-1} : = E_{\psi}$ \ induces \ $F_{k-1}(x)$ \ for each \ $x \in X$. This proves i) and ii).\\
To prove iii) it is enough, by an easy induction, to prove it for \ $j = 1$. In this case, it is sufficient to prove that the composition \ $\theta : = f_{\lambda}\circ \varphi : X \to \Xi_{\lambda}^{(k-2)}$, where \ $f_{\lambda} :  \Xi_{\lambda}^{(k-1)} \to  \Xi_{\lambda}^{(k-2)}$ \ is the quotient by \ $\Xi_{\lambda}^{(0)}$, is \ $(k-1)-$thematic, as we know from the proof of the proposition \ref{theme 3} that for \ $x \in X$ \ we have 
 $$F_1(x) = Ker\,f_{\lambda} \cap \A.\varphi(x).$$
  So we want to prove that if we have \ $\sum_{j=0}^{k-2} \ S_j(x).a^j.\varphi(x) \in \Xi_{\lambda}^{(0)}$, then \ $S_j(x) = 0 \quad \forall j \in [0,k-2]$. If all \ $S_j(x)$ \ are not zero, there exists \ $q \in \mathbb{N}$ \  and \ $T \in \C[[b]]$ \ invertible, such that
 $$ \sum_{j=0}^{k-2} \ S_j(x).a^j.\varphi(x)  = T.s^{\lambda+q-1} .$$
 But then, \ $(a - (\lambda+q).b).T^{-1}.\big(\sum_{j=0}^{k-2} \ S_j(x).a^j\big) $ \ is a polynomial in \ $a$ \ with coefficients in \ $\C[[b]]$, of degree \ $\leq k-1$, which annihilated \ $\varphi(x)$, contradicting the fact that \ $E(x) = \A.\varphi(x)$ \ has rank \ $k$. Now \ $\theta$ \ is \ $(k-1)-$thematic and induces the family \ $E(x)\big/F_1(x), x \in X$. $ \hfill \blacksquare$\\

 \begin{thm}\label{crt. hol.}
 Let \ $E(x)_{x \in X}$ \ be a family of rank \ $k \geq 2$ \ \ $[\lambda]-$primitive themes with fixed fundamental invariants \ $\lambda_1, p_1, \dots, p_{k-1}$ \ defined by a \ $\mathcal{O}_X-(a,b)-$module \ $\mathbb{E}$. Let \ $\alpha : X \to \C^*$ \ the function which associates to \ $x \in X$ \ the parameter of the rank \ $2$ \ $[\lambda]-$primitive theme \ $E(x)\big/F_{k-2}(x)$. Then the family is holomorphic if and only if the following two conditions are satisfied:
 \begin{enumerate}[i)]
 \item The family \ $F_{k-1}(x)$ \ is holomorphic, that is say there exists a \ $\mathcal{O}_X-(a,b)-$submodule \ $\mathbb{F}_{k-1}$ \ of \ $\mathbb{E}$ \ which is an holomorphic family and take the value \ $F_{k-1}(x)$ \ for each \ $x \in X$.
 \item The function \ $\alpha$ \ is holomorphic on \ $X$.
 \end{enumerate}
 \end{thm}
 
 \parag{Remark} Thanks to proposition \ref{holom. rank 2} the condition ii) of the previous theorem is equivalent to the fact that the family \ $E(x)\big/F_{k-2}(x), x \in X$ \ is holomorphic. So this result allows, using induction, to reduce the problem of the holomorphy of a family to the rank \ $2$ \ case.\\
 
 The proof will use the following lemma.
 
 \begin{lemma}\label{holom. 3}
 Let \ $j, q $ \ be natural integers and \ $\lambda$ \ be a rational number in \ $]0,1]$. Denote \ $H(j,q)$ \ the hyperplan of \ $\Xi_{\lambda}^{(j)}$ \ corresponding to the annulation of the coefficient \footnote{we use here the notations of the example given at the begining of section 5.1.}of \ $b^q.e_{\lambda, j}$. Then the \ $\C-$linear map
 $$ (a - (\lambda+q).b) : H(j,q)\oplus \C.b^q.e_{\lambda, j+1} \to b.\Xi_{\lambda}^{(j)} $$
 is an isomorphism of  Frechet spaces. So its inverse is linear and continuous.
 \end{lemma}
 
 \parag{Proof} The following equality is easy, for any \ $(h,m) \in \mathbb{N}$ \ with the convention \ $e_{\lambda,-1} = 0$
 $$ (a - (\lambda+q).b).b^m.e_{\lambda,h} = (m-q).b^{m+1}.e_{\lambda, h} + b^{m+1}.e_{\lambda, h-1} .$$
 So the image of \ $\Xi_{\lambda}^{(j)}$ \ by \ $ (a - (\lambda+q).b)$ \ is the hyperplane of \ $b.\Xi_{\lambda}^{(j)}$ \  given by the annulation of the coefficient of \ $b^{q+1}.e_{\lambda,j}$ \ and its kernel is \ $\C.b^q.e_{\lambda, 0}$. The conclusion is then easy. $\hfill \blacksquare$
 
 \parag{Remark} If we begin with an element in \ $b.\Xi_{\lambda}^{(j)}$ \ for which the coefficient of \ $b^q.e_{\lambda, j}$ \ is \ $\rho$, then the coefficient of \ $b^q.e_{\lambda, j+1}$ \ in its image by the inverse map will be also \ $\rho$. In particular, it will be non zero for \ $\rho \not= 0$.
 
 \parag{Proof of the theorem \ref{crt. hol.}} The statement is local on \ $X$ \ so we may assume that we have a \ $(k-1)-$thematic holomorphic map \ $\psi : X \to \Xi_{\lambda}^{(k-2)}$ \ such that \ $\mathbb{E}_{\psi}$ \ defines the family \ $F_{k-1}(x), x \in X$. We may also assume that the map
 $$ \psi - b^{\lambda_{k-1} - \lambda}.e_{k-2} $$
 takes its values in \ $\Xi_{\lambda}^{k-3}$, up to multiply \ $\psi$ \ by an invertible element of \ $\mathcal{O}(X)[[b]]$. Put \ $q : = \lambda_k - \lambda, \ S_{k-1} : = 1 + \alpha.b^{p_{k-1}} $ \ and define
 $$ \varphi : X \to \Xi_{\lambda}^{(k-1)} $$
 as the composition of \ $S_{k-1}.\psi$ \ with the inverse map constructed in the previous lemma, with \ $j : = k-2$.  and the obvious inclusion of \ $ H(k-2, \lambda_k - \lambda) \oplus \C.b^q.e_{\lambda, k-1}$ \ in \ $\Xi_{\lambda}^{(k-1)}$. Note that we have \ $\psi(x) \in b.\Xi_{\lambda}^{(k-2)}$ \ for each \ $x \in X$ \ because the inclusion  \ $F_{k-1}(x) \subset a.E(x) + b.E(x)$, \ $\Xi_{\lambda}$ \ has a simple pole  and any \ $\A-$linear map from \ $F_{k-1}(x)$ \ to \ $\Xi_{\lambda}$ \ may be extended\footnote{see [B.05].} to \ $E(x)$.\\
 Remark that the coefficient of \ $b^q.e_{\lambda, k-2}$ \ in \ $S_{k-1}.\psi$ \ co{\"i}ncides with the coefficient of \ $b^{p_{k-1}}$ \ in \ $S_{k-1}$, so it is given by \ $\alpha$. As it is non zero, the coefficient of \ $ b^q.e_{\lambda,k-1}$ \ in \ $\varphi$ \ does not vanish, thanks to the relation \ $\lambda_k - \lambda + 1 = \lambda_{k-1} - \lambda + p_{k-1}$. This is of course necessary in order that \ $\A.\varphi(x)$ \ will be a rank \ $k$ \ theme. Now it is easy to see that the holomorphic \ $k-$thematic map \ $\varphi$ \ is such that \ $\mathbb{E}_{\varphi}$ \ induces the family \ $E(x), x \in X$. $\hfill \blacksquare$
 
 \subsubsection{Holomorphy and duality.} 
 
 The last theorem of this paragraph shows that the twisted duality preserves the holomorphy of a family of \ $[\lambda]-$primitive themes.
 
 \begin{thm}\label{holom. dual}
 Let \ $E(x), x \in X$ \ be an holomorphic family of  \ $[\lambda]-$primitive themes parametrized by a reduced complex space. Let \ $\delta \in \mathbb{Q}$ \ such that for each \\
  $x \in X$ \ the (a,b)-module \ $E(x)^* \otimes_{a,b} E_{\delta}$ \ is a theme. Then the family 
   $$E(x)^* \otimes_{a,b} E_{\delta}, x \in X \quad {\rm  is \  holomorphic}.$$
 \end{thm}
 
 The proof will use the next lemma, which is an easy exercice left to the reader
 
 \begin{lemma}\label{holom. decal.}
 Let \ $E(x), x \in X$ \ be an holomorphic family of  \ $[\lambda]-$primitive themes parametrized by a reduced complex space. Let \ $\delta \in \mathbb{Q}$ \ such that for each \\
  $x \in X$ \ the (a,b)-module \ $E(x) \otimes_{a,b} E_{\delta}$ \ is a theme. \\
  Then the family \ $E(x) \otimes_{a,b} E_{\delta}, x \in X$ \ is holomorphic.
 \end{lemma}
 
 \parag{Proof of the theorem \ref{holom. dual}} We make an induction on the rank \ $k$. As the case \ $k = 1$ \ is trivial and the case \ $k = 2$ \ is already known (see the exercice following the definition \ref{parametre 1}), assume that the result is proved for \ $k-1 \geq 2$ \ and consider the situation in rank \ $k$. Let \ $F_1(x), x \in X$ \ be the family of normal subthemes of rank \ $1$ \ of the themes \ $E(x), x \in X$. The proposition \ref{holom. 2} implies the holomorphy of the family \ $E(x)\big/F_1(x)$. This proposition gives also the holomorphy of the family \ $F_2(x), x \in X$. The induction hypothesis implies the holomorphy of the family \ $(E(x)\big/F_1(x))^*\otimes_{a,b} E_{\delta}$ \ and also of the family \ $(F_2(x))^*\otimes_{a,b} E_{\delta}$. But now we may apply the theorem \ref{crt. hol.} to the family \ $(E(x))^*\otimes_{a,b} E_{\delta}$ \ because the corresponding family of subthemes of rang \ $k-1$ \ is the holomorphic family \ $(E(x)\big/F_1(x))^*\otimes_{a,b} E_{\delta}$ \  and the corresponding  family of rank \ $2$ \ quotients is the holomorphic family \ $(F_2(x))^*\otimes_{a,b} E_{\delta}$. $\hfill \blacksquare$

 \subsection{The canonical family.}
 
 \subsubsection{Definition and holomorphy.}
 
 We shall fix in this paragraph an integer \ $k \geq 1$ \  a rational number \ $\lambda_1 > k-1$ \ and natural integers \ $p_1, \dots, p_{k-1}$. For \ $j \in [1,k-1]$ \ we define the affine open set \ $W_j \subset V_j$ \ of  an affine hyperplane of the  complex  vector space \ $V_j$ (see the its definition  before  the proposition \ref{canonical form 1})
 $$ W_j : = \{ S_j \in V_j \ / \  S_j(0) = 1 \quad {\rm and \ the \ coefficient \ of} \quad b^{p_j} \quad {\rm is} \quad \not= 0 \} .$$
 Then put
 \begin{equation*}
 \mathcal{S}(\lambda_1, p_1, \dots, p_{k-1}) : = \{(S_1, \dots, S_{k-1}) \in \C[[b]]^{k-1} \ / \  S_j \in W_j \quad \forall j \in [1,k-1] \}. \tag{*}
 \end{equation*}
 Note that \ $ \mathcal{S}(\lambda_1, p_1, \dots, p_{k-1})$ \ is always a connected complex manifold.
 
 \begin{defn}\label{fam.can. 1}
 For \ $\sigma \in \mathcal{S}(\lambda_1, p_1, \dots, p_{k-1}) $ \ we define \ $E(\sigma)$ \ as the \ $[\lambda]-$primitive theme (with fundamental invariants \ $\lambda_1, p_1, \dots, p_{k-1}$) \ $\A\big/\A.P(\sigma)$ \ where 
 \begin{equation*}
 P(\sigma) : = (a - \lambda_1.b)S_1^{-1} \dots S_{k-1}^{-1}.(a - \lambda_k.b) \tag{**}
 \end{equation*}
  with  \ $\lambda_{j+1} = \lambda_j + p_j - 1$ \ for \ $j \in [1,k-1]$.
 \end{defn}
 
 \parag{Examples} For \ $k = 1$ \ and \ $\lambda_1$ \ given, the canonical family reduces to the theme \ $E_{\lambda_1}$.
 For \ $k = 2$ \ and \ $\lambda_1, p_1$ \ given, we have
 \begin{enumerate}
 \item For \ $p_1 = 0$ \ $\mathcal{S}(\lambda_1, 0) = \{1\}$ \  and the corresponding value of \ $P$ \ is
  $$P : = (a - \lambda_1.b).(a - (\lambda_1-1).b).$$
  
 \item  For \ $p_1 \geq 1$ \ we have \ $\mathcal{S}(\lambda_1,p_1) = \{ 1 + \alpha.b^{p_1}, \alpha \in \C^* \} \simeq \C^*$ \ and the \ $P$ \ corresponding to \ $\alpha \in \C^*$ \ is given by
  $$ P(\alpha) = (a - \lambda_1.b)(1 + \alpha.b^{p_1})^{-1}(a - (\lambda_1 + p_1 - 1).b) .$$
  \end{enumerate}
  
  \begin{thm}\label{fam. can. 2}
  For any given \ $k \geq 1, \lambda_1 > k-1 \quad {\rm rational} \quad {\rm and} \quad p_1, \dots, p_{k-1}$ \ natural integers, the canonical family is holomorphic.
  \end{thm}
  
  \parag{Proof} We shall deduce this result using theorem \ref{crt. hol.} and an induction on \ $k$. First remark that  the family \ $E(\sigma)\big/F_1(\sigma)$ \ is the pull back of  the canonical family parametrized by \ $\mathcal{S}(\lambda_2, p_2, \dots, p_{k-1})$ \ via the obvious projection (recall that \ $\lambda_2 : = \lambda_1 + p_1 - 1> k-2$ \ as \ $\lambda_1 > k-1$)
  $$ \mathcal{S}(\lambda_1, p_1, \dots, p_{k-1}) \to \mathcal{S}(\lambda_2, p_2, \dots, p_{k-1}) $$
  and the induction hypothesis; so  it is an holomorphic family. Moreover the family \ $F_2(\sigma)$ \ is also the pull back  via the obvious projection 
   $$ \mathcal{S}(\lambda_1, p_1, \dots, p_{k-1}) \to \mathcal{S}(\lambda_1, p_1)$$ 
    and the example above \ gives the holomorphy. So the theorem \ref{holom. dual} allows now to use the theorem \ref{crt. hol.} and this  gives the conclusion . $\hfill \blacksquare$
 
 \subsubsection{Versality and universality.}
 
 Again we shall fix \ $k \geq 1$, $ \lambda_1 > k-1$ \ a rational and natural integers \ $ p_1, \dots , p_{k-1}$ \ in this paragraph.
 
 \begin{defn}\label{versal. univ. 1}
 Let \ $X$ \ be a reduced complex space, \ $x_0$ \ a point in \ $X$, and let \ $\mathbb{E}$ \ an holomorphic family of rank \ $k$  \ $[\lambda]-$primitive themes with fundamtenal invariants \ $\lambda_1, p_1, \dots, p_{k-1}$. We shall say that this family is {\bf versal} near \ $x_0$ \  when for any reduced complex space \ $Y$ \ with a base point \ $y_0$ \  and any holomorphic family  \ $\mathbb{H}$ \ of rank \ $k$  \ $[\lambda]-$primitive themes with fundamental invariants \ $\lambda_1, p_1, \dots, p_{k-1}$, parametrized by \ $Y$, such that \ $H(y_0)$ \ is isomorphic to \ $E(x_0)$,   there exists an holomorphic map \ $ f : U \to X$ \ of an open neighbourhood \ $U$ \ of \ $y_0$ \ in \ $Y$ \ and an isomorphism of \ $\mathcal{O}_U-(a,b)-$modules \ $\theta : f^*(\mathbb{E}) \simeq  \mathbb{H}_{\vert U}$.\\
 We shall say that \ $\mathbb{E}$ \ is {\bf universal} near \ $x_0$  \  when  the germ at \ $y_0$ \ of  such an holomorphic map \ $f$ \  is always unique.
 The family will be called {\bf  locally versal} (resp. {\bf universal}) when it is versal (resp. universal) in a neighbourhood of each point.
 \end{defn}
 
 Of course a versal family contains all isomorphy class of rank \ $k$  \ $[\lambda]-$primitive themes with fundamental invariants \ $\lambda_1, p_1, \dots, p_{k-1}$. Any  isomorphy class appears exactely one time in an universal family. The existence of an universal family is the same problem that the representability of the functor which associates to a reduced complex space \ $Y$ \ the set of holomorphic families of rank \ $k$  \ $[\lambda]-$primitive themes with fundamental invariants \ $\lambda_1, p_1, \dots, p_{k-1}$, parametrized by \ $Y$. \\
 The existence of a versal family implies the finiteness of the number of complex parameters in order to determine an isomorphy class of such a theme.
 
 \begin{thm}\label{versal univ. 2}
 For any choice of \ $k \geq 1$, $ \lambda_1 > k-1$ \ a rational and natural integers \ $ p_1, \dots , p_{k-1}$ \ the canonical family parametrized by \ $\mathcal{S}(\lambda_1, p_1, \dots, p_{k-1})$ \ is locally  versal.\\
 If any theme corresponding to  a choice of \ $k \geq 1$, $ \lambda_1 > k-1$ \ a rational and natural integers \ $ p_1, \dots , p_{k-1}$ \ has an unique canonical form, then this canonical family is (globally) universal.
 This is  the case, for instance, when \ $p_j \geq k-1 \quad \forall j \in [1,k-1]$.
 \end{thm}
 
 \parag{Proof} First we shall show the versality by induction on the rank \ $k \geq 1$. The cases \ $k = 1$ \ and \ $k = 2$ \ are already known (see above section 5.1.), so we shall assume that \ $k \geq 3$ \ and the case \ $k-1$ \ known. \\
 Let \ $\mathbb{E}$ \ be an holomorphic family of rank \ $k$  \ $[\lambda]-$primitive themes with fundamental invariants \ $\lambda_1, p_1, \dots, p_{k-1}$ \ parametrized by the reduced complex space \ $X$. Fix \ $x_0$ \ in \ $X$ \ and denote \ $\mathbb{F}_1 \subset \mathbb{E}$ \ the subsheaf of \ $\mathcal{O}_X-(a,b)-$modules \ defining the family of normal subthemes of rank 1 in \ $\mathbb{E}$. Then \ $\mathbb{E}\big/\mathbb{F}_1$ \ is an holomorphic family of rank \ $k-1$  \ \ $[\lambda]-$primitive themes with fundamental invariants \ $\lambda_2, p_2, \dots, p_{k-1}$, with \ $\lambda_2 : = \lambda_1 + p_1 -1$. So thanks to the induction hypothesis, there exists an open neighbourhood\ $U$ \  of \ $x_0$ \ in \ $X$ \ and an holomorphic map \ $g : U \to \mathcal{S}(\lambda_2, p_2, \dots, p_{k-1})$ such that \ $g^*(\mathbb{S}) $ \ is isomorphic to \ $( \mathbb{E}\big/\mathbb{F}_1)_{\vert U}$, where \ $\mathbb{S}$ \ is the canonical family parametrized by \ $\mathcal{S}(\lambda_2, p_2, \dots, p_{k-1})$.\\
 As the result is local near \ $x_0$, we may assume, up to thrink \ $U$ \ around \ $x_0$,  that \ $\mathbb{E}$ \ is given by an holomorphic \ $k-$thematic map \ $\varphi : U \to \Xi_{\lambda}^{(k-1)}$, which satisfies
 $$ \varphi(x) = s^{\lambda_k-1}.(Log\,s)^{k-1} + \psi(x)  $$
 where \ $\psi : U \to  \Xi_{\lambda}^{(k-2)}$ \ is holomorphic and \ $(k-1)-$thematic. The holomorphic map \ $g$ \ gives  in fact holomorphic maps \ $S_2, \dots, S_{k-1} : U \to 1 +b.\C[[b]]$ \ such that if we define
 $$P_1 : = (a -\lambda_2.b).S_2^{-1}.(a - \lambda_3.b).S_3^{-1} \dots S_{k-1}^{-1}.(a - \lambda_k.b) $$
 \ $\A.P_1$ \ will be the annihilator of the canonical generator \ $e$ \ of the canonical family parametrized by \ $\mathcal{S}(\lambda_2, p_2, \dots, p_k)$. So the generator  \ $g^*(e)$ \ of \ $\mathbb{E}\big/\mathbb{F}_1$ \ will satisfy also \ $P_1.g^*(e) = 0$. If we identify \ $\mathbb{E}_{\vert U}$ \ to \ $\mathbb{E}_{\varphi} \subset \Xi_{\lambda}^{(k-1)}$, we identify \ $\mathbb{E}\big/\mathbb{F}_1$ \ to a sub-sheaf of the quotient  sheaf  
 $$ \Xi_{X, \lambda}^{(k-1)}\Big/ \Xi_{X, \lambda}^{(0)}\  \simeq  \  \Xi_{X, \lambda}^{(k-2)}. $$
 So we may find \ $T_0, \dots, T_{k-1}$ \ local sections of \ $\mathcal{O}_X[[b]]$ \ such that the image of the section
 $$ \gamma : = \sum_{j=0}^{k-1} \ T_j.a^j.\varphi $$
 in \ $\mathbb{E}\big/\mathbb{F}_1$ \  co{\"i}ncides with \ $g^*(e)$. As \ $g^*(e)$ \ generates \ $\mathbb{E}\big/\mathbb{F}_1$ \ the section \ $T_0$ \ of \ $\mathcal{O}_X[[b]]$ \ must be  invertible near \ $x_0$. So \ $\gamma$ \ generates \ $\mathbb{E}$ \ near \ $x_0$ \ and satifies \ $P_1.\gamma \in \mathbb{F}_1$. As \ $\mathbb{F}_1 = \mathcal{O}_x[[b]].s^{\lambda_1-1}$ \ we may write
 $$ P_1.\gamma = \Theta.s^{\lambda_1-1} $$
 where \ $\Theta$ \ is a section of \ $ \mathcal{O}_x[[b]]$ \ in on an open neighbourhood of \ $x_0$. The decomposition \ $E_{\lambda_1} = P_1.E_{\lambda_1} \oplus V_1$ \ allows to write
 $$\Theta.s^{\lambda_1-1} = (P_1.\beta + S_1).s^{\lambda_1} $$
 where \ $\beta$ \ and \ $S_1$ \ are local sections respectiveley of \ $\mathcal{O}_X[[b]]$ \ and \ $\mathcal{O}_X\otimes V_1$. Moreover, as \ $\Theta$ \ is invertible, which means invertibility in \ $\mathcal{O}_X$ \ of its constant term (in \ $b$), so is \ $S_1$. So, up to multiplication by an holomorphic invertible function \ $I$ \   on an open neighbourhood of \ $x_0$, we may assume that the constant term of \ $S_1$ \ is \ $1$. Then \ $\tau : = I.(\gamma - \beta.s^{\lambda_1-1})$ \ is still a generator of \ $\mathbb{E}$ and satisfies
 $$ P_1\tau = S_1.s^{\lambda_1-1} \quad {\rm with} \quad S_1 \in \mathcal{O}_X \otimes V_1, S_1(0) = 1 .$$
 this gives \ $(a - \lambda_1.b)S_1^{-1}.P_1.\tau = 0 $. So we see that \ $\mathbb{E}$ \ is isomorphic to the pull back of the canonical family by the holomorphic map \ $f  : U(x_0) \to \mathcal{S}(\lambda_1, p_1, \dots, p_{k-1})$ \ given in the open neighbourhood \ $U(x_0)$ \ of \ $x_0$ \ by \ $S_1, \dots, S_k$, by sending the local generator \ $\tau$ \ on \ $ f^*(e)$ \ where \ $e$ \ is the standard generator of the canonical family parametrized by \ $ \mathcal{S}(\lambda_1, p_1, \dots, p_{k-1})$. $\hfill \blacksquare$

 \subsubsection{Examples.}
 
We shall give an example of rank \ $3$ \ \ $[\lambda]-$primitive themes for which there does not exists, even locally, an univeral family.\\

Let  the fundamental invariants be \ $\lambda_1, p_1 = p_2 = 1$;  so we have \ $q_1 = p_1 + p_2 =  2 $ \  and \  $q_2 = p_2 = 1$. Remark that the corresponding themes are special.

\begin{prop}\label{ex. 1}
There does not exist an universal family near each invariant special rank \ $3$ \ $[\lambda]-$primitive  theme.
\end{prop}

The proof will use  the following three lemmas.

\begin{lemma}\label{isom.}
Let  \ $\alpha,\beta, \gamma \in \C, \alpha.\beta \not= 0 $ \ and consider the rank \ $3$ \ themes defined as follows :
\begin{align*}
&  (a - \lambda.b).e_3 = (1 + \alpha.b).e_2 \\
&  (a - \lambda.b).e_2 = (1 + \beta.b + \gamma.b^2).e_1 \\
& (a - \lambda.b).e_1 = 0 .
\end{align*}
For \ $\beta \not= \alpha$, \ $E_{\alpha,\beta,\gamma}$ \ is isomorphic to \ $E_{\alpha,\beta, 0}$ \ for each \ $\gamma$.\\
For  \ $\beta = \alpha$, the themes \ $E_{\alpha,\alpha,\gamma}$ \ and \ $E_{\alpha,\alpha,\gamma'}$ \ are isomorphic if and only if \ $\gamma = \gamma'$.
\end{lemma}

\parag{Proof} We look for a \ $\C[[b]]-$basis \ $\varepsilon_3, \varepsilon_2, \varepsilon_1$ \ of \ $E_{\alpha,\beta,\gamma}$ \ satisfying the following conditions :
\begin{align*}
& \varepsilon_3 = e_3 + U.e_2 + V.e_1, \quad {\rm with} \quad U,V \in \C[[b]] \tag{0}\\
& (a - \lambda.b).\varepsilon_3 = (1 + \alpha.b).\varepsilon_2 \tag{1} \\
& (a - \lambda.b).\varepsilon_2 = (1 + \beta.b + \gamma'.b^2).\varepsilon_1 \tag{2}\\
& (a - \lambda.b).\varepsilon_1 = 0. \tag{3}
\end{align*}
We know that \ $\alpha$ \ and \ $\beta$ \ are determined by the isomorphism class of the theme \ $E(\alpha,\beta,\gamma)$:  as we have \ $p_1 = p_2 = 1$ \ they are the parameters of the rank \ $2$ \ themes \ $E\big/F_1$ \ and  \ $F_2$ \ respectively .\\
Remark also that \ $(3)$ \ implies that \ $\varepsilon_1 = \rho.e_1$ \ with \ $\rho \in \C^*$.\\
Let us compute the conditions for \ $U$ \ and \ $V$ :
\begin{align*}
& (a - \lambda.b).\varepsilon_3 =   (1 + \alpha.b).e_2 + b^2.U'.e_2 + U. (1 + \beta.b + \gamma.b^2).e_1 + b^2.V'.e_1 \\
&\qquad\qquad   = (1+\alpha.b).\varepsilon_2
\end{align*}
and so, we have \ $ \varepsilon_2 = Z.e_2 + T.e_1$ \ with \ $ Z = (1 +\alpha.b)^{-1}.(1 + \alpha.b + b^2.U') $ , and 
\begin{equation*}
 (1 + \alpha.b).T = U.(1 + \beta.b + \gamma.b^2) + b^2.V'. \tag{4}
\end{equation*}
Then we get
\begin{align*}
& (a - \lambda.b).\varepsilon_2 = Z.(1 + \beta.b + \gamma.b^2).e_1 + b^2.Z'.e_2 + b^2.T'.e_1 \\
& \qquad \qquad = (1 + \beta.b + \gamma'.b^2).\rho.e_1
\end{align*}
and this implies \ $Z' = 0 $ \ and as  \ $Z = 1 + (1 + \beta.b)^{-1}.b^2.U'$ \ we must have \ $U \in \mathbb{C}$, and \ $Z = 1$.  The relation \ $(2)$ \ implies now, as  \ $\varepsilon_2 = e_2 + T.e_1$
$$  (1 + \beta.b + \gamma.b^2).e_1 + b^2.T'.e_1 =  (1 + \beta.b + \gamma'.b^2).\rho.e_1 .$$
Then we have \ $\rho = 1$ \ and \ $ T' = \gamma' - \gamma $. \\
Looking at the constant terms in \ $(4)$ \ we obtain  \ $T = U + (\gamma' - \gamma).b$.\\
But \ $(4)$ \ implies also
\begin{equation*}
 \alpha.U + \gamma' - \gamma = U.\beta \qquad \qquad {\rm and}  \quad
 U.\gamma + V' = \alpha.(\gamma - \gamma') \tag{5} 
\end{equation*}
So, for \ $\alpha \not= \beta $ \ we will have
$$U = \frac{\gamma - \gamma'}{\alpha - \beta} \qquad {\rm and} \qquad  V = V_0 + \frac{\gamma'- \gamma}{\beta - \alpha}.\big(\alpha.(\beta-\alpha) - \gamma\big).b .$$
If \ $\beta = \alpha$, the relation \ $(5)$ \ imposes  \ $\gamma = \gamma'$. $\hfill \blacksquare$

\bigskip

For  \ $ (\alpha,\beta) \in (\C^*)^2, \alpha \not= \beta $ \ denote \ $E(\alpha, \beta)$ \ the rank \ $3$ \ theme defined by 
 $$E(\alpha,\beta) : = \A\big/\A.(a -\lambda.b)(1+\beta.b)^{-1}(a - \lambda.b)(1 + \alpha.b)^{-1}(a - \lambda.b).$$

\begin{lemma}\label{non invariant}
There is no rank \ $2$ \ endomorphism of \ $E(\alpha,\beta)$  \ for \ $\alpha \not= \beta$.
\end{lemma}

\parag{Proof} It is enough to prove that there exists no element \ $x : = e_2 + U.e_1$ \ in \ $E(\alpha,\beta)$ \ such that \ $(a-\lambda.b)(1+\alpha.b)^{-1}(a -\lambda.b).x = 0 $, where \ $U \in \C[[b]]$. As each element in  \ $E(\alpha,\beta)$ \ which is annihilated by  \ $(a - \lambda.b)$ \ is in \ $\C.e_1$, such an \ $x$ must satisfies 
$$ (a - \lambda.b)x = \rho.(1 + \alpha.b).e_1 $$
which implies that \ $U$ \ is a solution of the equation
$$ (1 + \beta.b) + b^2.U' = \rho.(1 + \alpha.b) .$$
So we conclude that we must have \ $\rho = 1$ \ and so \ $\alpha = \beta$. $\hfill \blacksquare$

\bigskip

But if  \ $\alpha = \beta \not= 0$ \ we shall have an invariant theme for any \ $\gamma$.

\begin{lemma}\label{Inv. alpha = beta}
For \ $\alpha \not= 0$ \ the (a,b)-module \  $E_{\alpha,\alpha,\gamma}$ \ is an invariant rank \ $3$ \ theme.
\end{lemma}

\parag{Proof} It is enough to find  \ $x : = e_2 + U.e_1$  \ such that
$$(a-\lambda.b)(1+\alpha.b)^{-1}(a -\lambda.b).x = 0 , $$
where \ $U \in \C[[b]]$. \\
As \ $F_2$ \ is a theme, the kernel of \ $(a - \lambda.b)$ \ in \ $F_2$ \ is \ $\C.e_1$. So \ $x$ \ must satisfies
$$ (a - \lambda.b)x = \rho.(1 + \alpha.b).e_1 $$
and this implies that \ $U$ \ is solution of the equation
$$ (1 + \alpha.b+ \gamma.b^2) + b^2.U' = \rho.(1 + \alpha.b) .$$
So we conclude that we must have \ $\rho = 1$ \ and  \ $U = -\gamma.b + cste$. So we obtain a solution \ $ x : = e_2 -\gamma.b.e_1$. $\hfill \blacksquare$

\parag{Proof of the proposition \ref{ex. 1}} The fact that the invariant themes in this canonical family are exactely the \ $E_{\alpha,\alpha,\gamma}$ \ is proved in lemmas \ref{non invariant} and \ref{Inv. alpha = beta}.\\
The canonical family  \ $(E_{\alpha,\beta,\gamma})_{(\alpha,\beta,\gamma) \in S(\lambda_1,  p_1 = p_2 =1)}$ \ is holomorphic and versal near any point thanks to the theorem \ref{versal univ. 2}. Let assume that we have find an holomorphic family  \ $(E_y)_{y \in Y}$ \  which is locally universal around a theme \ $ E(\alpha_0,\alpha_0,\gamma_0) \simeq E_{y_0}$, where \ $Y$ \ is a reduced complex space that we may assume to be embedded in \ $\C^N$ \ near \ $y_0$. Let  $\varphi : \Omega \to Y \hookrightarrow \C^N$ \ the classifying map for the canonical family on an open set  \ $\Omega$ \ de \ $(\alpha_0,\alpha_0,\gamma_0) \in (\C^*)^2\times \C $. As for \ $\alpha \not= \beta$ \ the isomorphy class of the theme $E_{\alpha,\beta,\gamma}$ \  does not depend on \ $\gamma$, thanks to the lemma \ref{isom.}, we shall have \ $\frac{\partial \varphi}{\partial \gamma} \equiv 0$ \ on the open set \ $\{\alpha \not= \beta\}$ \ of \ $\Omega$. This implies that \ $\varphi$ \  does not depends on \ $\gamma$ \ for all  \ $\alpha$ \ near enough to \ $\alpha_0$, and for all \ $\gamma,\gamma'$ \ near enough from \ $\gamma_0$. This contradicts the lemma \ref{isom.}. $\hfill \blacksquare$

\begin{cor}\label{Univers.}
The family  \ $E(\alpha,\beta)_{(\alpha,\beta)\in X}$ \ is universal near any point in \\
 $X : = (\C^*)^2 \setminus \{\alpha = \beta\} $.
\end{cor}

\parag{Proof} Let \ $\mathbb{E}$ \ be the sheaf on \ $X$ \ of \ $\mathcal{O}_X-(a,b)-$modules associated to the holomorphic family \ $E(\alpha,\beta)$. It is enough to prove that the holomorphic map
$$\pi :  X \times \C \to X $$
defined by \ $\pi(\alpha,\beta, \gamma) = (\alpha,\beta)$ \ is such that the sheaf \ $\pi^*(\mathbb{E})$ \ is isomorphic as a sheaf of \ $\mathcal{O}_X-(a,b)-$modules to the sheaf associated to the canonical family restricted to \ $X \times \C$. But the inverse of the isomorphism we are looking for is given by the computation in the lemma \ref{isom.} which gives, in the case \ $\gamma' = 0$, where \ $(\alpha,\beta,\gamma)$ \ is taken as an holomorphic parameter in \ $X \times \C$, holomorphic sections \ $U, V$ \ of the sheaf \ $\mathcal{O}_{X\times \C}[[b]]$. The inverse isomorphism is obtained by sending the generator \ $e_3$ \ of the canonical family to  \ $\varepsilon_3(\gamma' = 0) : = e_3 + U.e_2 + V.e_1$, which is the generator of the family \ $\pi^*(\mathbb{E})$. $\hfill \blacksquare$

 \parag{Another example without universal family}
 
 We give here an example of {\bf non special} rank \ $4$ \ $[\lambda]-$primitive themes without universal family.
 Fix the fundamental invariants to be \ $\lambda_1 > 3$ \ and \ $p_1 = 3, p_2 = p_3 = 2$. So the canonical family is defined by the following relations :
 
 \begin{align*}
 & (a - (\lambda_1+4).b).e_4 = S_3.e_3 \quad S_3 : = 1 + \alpha.b^2  \quad\quad (q_3 = p_3 = 2) \\
 & (a - (\lambda_1+3).b).e_3 = S_2.e_2 \quad  S_2 : = 1 + \beta.b + \gamma.b^2 \quad\quad (q_2 = p_2 = 2) \\
 & (a - (\lambda_1+2).b).e_2 = S_1.e_1 \quad S_1 : = 1 + \delta.b + \varepsilon.b^2 + \theta.b^3 \quad\quad  (q_1 = p_1 = 3)\\
 & (a - \lambda_1.b).e_1 = 0.
 \end{align*}
 with the condition  \ $ \alpha.\gamma.\theta \not= 0 $.\\
 
 We look for invariant themes in this family. So we look for an element \ $x \in F_3 \setminus F_2$ \ such that \ $P_0.x = 0$. The existence of such an \ $x$ \  is equivalent to the existence of a rank \ $3$ \ endomorphism.\\
 Using proposition \ref{base standard 1} the condition \ $P_0.x = 0$ \ is in fact equivalent to the condition \ $P_1.x = 0$ \ because we ask \ $x$ \ to be in \ $F_3$. This gives the equality \ $(a - \lambda_2.b).S_2^{-1}.P_2.x = 0$.
 The kernel in \ $E$ \ of \ $(a - \lambda_1.b)$ \ is \ $\C.e_1$, so we obtain the equation
 \begin{equation*}
 P_2.x = \rho.S_2.b^2.e_1 \tag{1} 
 \end{equation*}
 It gives \ $(a - \lambda_3.b).S_3^{-1}.P_3.x = \rho.S_2.b^2.e_1$ \ which implies that the class of \  \ $S_3^{-1}.P_3.x$ \ in \ $E\big/F_1$ \  is in the kernel of \ $(a - \lambda_3.b)$ \ which is  \ $\C.b^{\lambda_3-\lambda_2}.e_2 \quad modulo \  F_1$. So we may put :
 \begin{equation*}
 P_3.x = \sigma.S_3.b.e_2 + S_3.T.e_1 \tag{2}
 \end{equation*}
and the equation \ $(1)$ \ implies 
 \begin{align*}
 &   P_2.x = (a - \lambda_3.b)( \sigma.b.e_2 + T.e_1) = \rho.S_2.b^2.e_1 \\
 & \qquad = \sigma.b.S_1.e_1 + (b^2.T' - 3b.T) =  \rho.S_2.b^2.e_1. 
 \end{align*}
So we get \ $ b^2.T' - 3.b.T = \rho.S_2.b^2 - \sigma.S_1.b $ \ and after simplification by \ $b$ \ this gives 
 \begin{equation*}
 b.T' - 3.T = \rho.S_2.b - \sigma.S_1. \tag{3}
 \end{equation*}
 Then we conclude that we must have \ $\rho.\gamma = \sigma.\theta$. But as \ $\gamma.\theta \not= 0$, the condition \ $\rho \not= 0 $ \ is equivalent to \ $\sigma \not= 0$.\\
 
 We shall have a solution \ $T$ \ by putting \ $\sigma : = \rho.\gamma/\theta$, and it satisfies   \ $-3T(0) = -\sigma$, that is to say \ $T(0) = \sigma/3 =  \rho.\gamma/3\theta$.\\
 Define now  \ $x : = X.e_3 + U.e_2 + V.e_1 $. The relation \ $(2)$ \ gives 
 $$ (a - \lambda_4.b).x = (b^2.X' - b.X).e_3 \quad {\rm modulo} \ F_2 $$
 and so  \ $X = \tau.b$ \ with \ $\tau \in \C^*$ \ as we assume \ $x \not\in F_2$. Then we get
 \begin{align*}
& (a - \lambda_4.b).x = (\tau.b.S_2 + b^2.U' - 2b.U).e_2 + (b^2.V' - 4b.V + U.S_1).e_1 \\
& \qquad    = \sigma.S_3.b.e_2 + S_3.T.e_1 
\end{align*}
and so the equations
\begin{align*}
& b.U' - 2U = \sigma.S_3 - \tau.S_2  \tag{5} \\
& b^2.V' - 4b.V = S_3.T - U.S_1 \tag{6}
\end{align*}
The first one implies that \ $\tau.\gamma = \sigma.\alpha$ \ and forces \ $U(0) = (\tau - \sigma)/2= \sigma.(\alpha/\gamma -1)/2  $, the second one imposes
$$ T(0) = U(0) .$$
Then we obtain
$$ T(0) = \sigma/3  = \sigma(\alpha/\gamma - 1)/2 . $$
So, if \ $ \alpha/\gamma \not= 1 + 2/3 = 5/3 $ \ we must have  \ $\sigma = 0$ \ and this imposes \ $\tau = 0$ \ and this is incompatible with our assumtion \ $x \in F_3 \setminus F_2$.\\
Then for  \ $3\alpha \not= 5\gamma $ \ the theme is not invariant, and as  \ $E/F_1$ \ is invariant thanks to the theorem \ref{unique 1} (here the rank is \ $3$ \ and each \ $p_i$ \ is at least \ $2$), the canonical form will not be unique thanks to lemma \ref{U 2}.\\
Now for \ $3\alpha = 5\gamma $ \ we find a solution \ $U$ \ of the equation \ $(5)$ \ and then a solution \ $V$ \ for \ $(6)$ \ when we choose a solution \ $T$ \ of \ $(3)$. So along the hyperplane  \ $3\alpha = 5\gamma $ \ the themes are invariant.  The situation is similar to the previous example.\\

Let us show, to conclude the proof, that when \ $3\alpha = 5\gamma $ \ the coefficient  \ $\varepsilon$ \ is not relevant to determine the isomorphism class of \ $E$. As we know that \ $\theta$ \ is determined by the isomorphy class of \ $E$, it is enough to prove that \ $\delta$ \ is also determined by the isomorphy class of \ $E$.\\
 Because, up to an homothetie, any automorphism of \ $E$ \ is obtained  by sending \ $e_4$ \ to \ $e_4 + y$ \ where \ $y \in F_3$  \ and satisfies \ $P_1.y \in F_1$, to prove that \ $\delta$ \ is determined by the isomorphy class of \ $E$, it is enough to prove that any such \ $y$ \ satisfies \ $P_1.y \in b^2.F_1$.\\
But  \ $P_1.y \in F_1$ \  shows that if we define \ $\varphi(e) = y$ \ this produces an \ $\A-$linear map \ $E\big/F_1 \to F_3\big/F_1$. As the rank is \ $\leq 2$, $\varphi(F_2\big/F_1) = 0$ \ and so  \ $P_2.e \in F_2$ \ implies \ $P_2.y \in F_1$. So we have 
$$ (a - \lambda_3.b).S_3^{-1}.P_2.y \in F_1 .$$
The kernel of \ $(a - \lambda_3.b)$ \ acting on  \ $F_3\big/F_1$ \ is \ $\C.b.e_2 + F_1$ \ as \ $\lambda_3 - \lambda_2 = 1$. Then we may write
$$ P_3.y = \rho b.S_3.e_2 + S_3.T.e_1 $$
which gives 
$$ P_2.y = \big[ \rho b.S_1 + b^2.T' - 3b.T \big].e_1 \in b.F_1 .$$
So  \ $P_1.y \in b^2.F_1$ \ is proved. \\
It is not necessary to prove that \ $\varepsilon$ \ can move by some isomorphism because  this is consequence of the lemma \ref{U 2}. $\blacksquare$

\parag{Conclusion} So in this example we have the following properties, analoguous of the properties in the previous rank \ $3$  example:
\begin{enumerate}
\item For \ $3\alpha \not= 5\gamma$ \ the themes are not invariant.
\item For \ $3\alpha = 5\gamma$ \ each theme is invariant but there is no local universal family around sucn an isomorphy class.
\item For each theme \ $E$ \  such that \ $3\alpha \not= 5\gamma$ \ there exists a local universal family obtained by the restriction of the canonical family to an open neighbourghood of \ $E$ \ in the hyperplane \ $\{\varepsilon = \varepsilon_0\}$.\\
\end{enumerate}

\parag{Remark} Let \ $E$ \ be a theme of the previous family with \ $3\alpha \not= 5\gamma$. Then  there is no injective map from \ $E\big/F_1$ \ to \ $F_3$ \ because \ $E$ \ is not invariant. But the fundamental invariants of these \ $[\lambda]-$primitive rank \ $3$ \ themes are respectively
\begin{align*}
& \mu_1 = \lambda_1 + 2, \  \mu_2 = \lambda_1+ 3, \  \mu_3 = \lambda_1 + 4 \quad {\rm for} \quad E\big/F_1\\
& {\rm and } \quad \lambda_1, \  \lambda_1+ 2, \ \lambda_1 + 3 \quad {\rm for} \quad F_3.
\end{align*}
So we have \ $\mu_1 - \lambda_1 = 2, \mu_2 - \lambda_2 = 1, \mu_3 - \lambda_3 = 1$; this shows that  the condition \ $\mu_j - \lambda_j \geq k-1\quad \forall j \in [1,k]$ \ of the theorem \ref{hom. 1} is sharp.

 \section{Appendix.}
 
 \subsection{A lemma.}
 
 The following result  from [B.95] plays a key role in the construction of the canonical form and so in the construction of the canonical families. We give here the main lines of the proof for the convenience of the reader.

\begin{lemma}\label{Ext}
Let \ $E$ \ and \ $F$ \ be two regular (a,b)-modules. Then we have :
$$ \dim_{\C}(Ext^1_{\A}(E,F)) - \dim_{\C}(Ext^0_{\A}(E,F)) = rg(E).rg(F).$$
\end{lemma}

\parag{Proof} First we show the case where \ $rg(E) = 1$.\\
Then we have \ $E \simeq \A\big/\A.(a - \lambda.b)$ \ so \ $Ext^0_{\A}(E,F)$ \ and \ $Ext^1_{\A}(E,F)$ \ are respectiveley the kernel and cokernel of \ $a - \lambda.b : F \to F $. We shall prove the formula by induction on the rank of \ $F$. For the rank \ $1$ \ we have \ $F \simeq \A\big/\A.(a - \mu.b)$,  and the computation is easy :
\begin{enumerate}
\item For \ $\lambda \not\in \mu + \mathbb{N}$ \  we have  \ $Ext^0_{\A}(E,F) = \{0\}$ \ and \ $Ext^1_{\A}(E,F) \simeq \C.e_{\mu}$.
\item Then \ $Ext^0_{\A}(E,F) = \C.b^n.e_{\mu}$ \ and \
$Ext^1_{\A}(E,F) \simeq \C.e_{\mu} \oplus \C.b^{n+1}.e_{\mu}$ \ for \ $\lambda = \mu + n $.
\end{enumerate}
and the formula follows.\\
Assume that \ $rk(F) \geq 2$ \ and the formula proved for lower ranks. There exists an  exact sequence
$$ 0 \to G \to F \to E_{\mu} \to 0 $$
and \ $rk(G) = rk(F) - 1$. Now we have a long exact sequence of finite dimensional complex vector spaces (see [B.95] theorem 1)
\begin{align*}
& 0 \to Ext^0_{\A}(E,G) \to Ext^0_{\A}(E,F)\to Ext^0_{\A}(E,E_{\mu})\to \\
& \to Ext^1_{\A}(E,G) \to Ext^1_{\A}(E,F) \to Ext^1_{\A}(E,E_{\mu})\to 0 
\end{align*}
The alternating sum of the dimension is zero and we obtain
\begin{align*}
&  \dim(Ext^1_{\A}(E,G)) -  \dim(Ext^0_{\A}(E,G)) + \dim( Ext^1_{\A}(E,E_{\mu})) - \dim( Ext^0_{\A}(E,E_{\mu})) = \\
& \dim(Ext^1_{\A}(E,F)) - \dim(Ext^1_{\A}(E,F)) = (rg(F) - 1) + 1 = rg(F)
\end{align*}
and we conclude thanks to the induction hypothesis.\\
The general case for \ $E$ \ with \ $rk(F) = 1$ \  is obtained in the same manner. \\
An induction on \ $rk(E) + rk(F)$ \ gives the general case by an analoguous \\
computation. $\hfill \blacksquare$

\subsection{Existence of holomorphic k-thematic maps.}\label{App.2}

Our purpose here is to show that if we have a proper holomorphic map \\
 $F : X \to D\times T$ \ between two complex manifolds which is a submersion over \ $D^*\times T$, then to the following data 
\begin{enumerate}[i)]
\item A smooth \ $T-relative$  \ $(p+1)-$differential form \ $\omega$ \ on \ $X$ \ such that \ $d_{/T}\omega = 0 = df \wedge \omega$ \ where \ $f : X \to D$ \ is the composition of \ $F$ \ with the projection on \ $D$. \item  A vanishing cycle \ $\gamma$ \ in the generic fiber of \ $F$ \ which is in the generalized eigenspace of the monodromy of  the restriction of \ $F$ \ over \ $D \times \{t_0\}$.
\end{enumerate}
we may associated locally around the generic point of \ $T$ \ an \ $k-$thematic holomorphic map which defines the family of themes associated to the family parametrized by \ $T$  \ of vanishing periods.\\

We shall begin by an easy  lemma of algebraic geometry over the algebra  \ $Z : = \C[[b]]$.

 \begin{lemma}\label{rg. sem. cont.1}
 Let \ $E$ \ be a regular rank \ $k$ \  (a,b)-module. Fix a \ $\C[[b]]-$basis  \ $e_1, \dots, e_k$ \  of \ $E$ \ and consider \ $E$ \ as the affine space  \ $Z^k$ \ over \ $Z$. Then for each natural integer \ $p$ the subset  \ $X_p \subset E = Z^k$ \  defined by
 $$ X_p : = \{ x \in E \ / \ rg( \A.x) \leq p \} $$
 is an algebaric subset of \ $E$, that is to say that there exists finitely many polynomials \ $P_1, \dots, P_N$ \ in  \ $Z[x_1, \dots, x_k]$ \ such that we have
 $$ X_p = \{ x \in Z^k \ / \ P_j(x) = 0 \quad \forall j \in [1,N] \}.$$
 \end{lemma}
 
 \parag{Proof} As \ $E$ \ is regular with rank \ $k$, for each \ $x \in E$ \ the sub-(a,b)-module \ $\A.x$ \ is monogenic and regular with rank \ $\leq k$. So it is generated as \ $\C[[b]]-$module by \ $x, a.x, \dots a^{k-1}.x$. To write that the rank of \ $\A.x$ \ is \ $\leq p$ \ it enough to write that all  \ $(q,q)$ \ minors of the matrix of these \ $k$ \ vectors in the basis  \ $e_1, \dots, e_k$ \ are zero for all \ $q \in [p+1,k]$. This gives the polynomiales \ $P_1, \dots, P_N$. $\hfill \blacksquare$\\
 
Here is an immediate consequence.
 
 \begin{cor}\label{rg. sem. cont.2}
 Let \ $X$ \ be a reduced complex space and \ $E$ \ be a rank \ $k$ \ regular (a,b)-module. Let \ $f : X \to  E$ \ be an holomorphic map\footnote{by fixing a basis \ $e_1, \dots, e_k$ \ of \ $E$ \ on \ $\C[[b]]$ \ this is a global section of the sheaf \ $\mathcal{O}_X[[b]]^k$.}. Then we have a finite  stratification
  $$ X_0 \subset X_1 \subset \dots \subset X_k = X $$
  by closed analytic subsets, such that, for any \ $q \in [1,k]$ \ the subset \ $X_q\setminus X_{q-1}$ \ is exactely the set of \ $x \in X$ \ where the rank of \ $\A.x$ \ is equal to \ $q$.
  \end{cor}
  
  \parag{Remark} The quotient of two holomorphic functions \ $f : D \to \C[[b]]$ \ and \\
  $g : D \to \C[[b]]$ \ with \ $g(0) \not= 0$ \ may be well defined for each value of \  \ $z \in D$, and the function \ $f\big/ g$ \ may not be holomorphic on \ $D$. For instance this is the case for \ $z \to \frac{z+b^2}{z+b} $, because a relation like
  $$ z + b^2 = (z+b).(\sum_{j=0}^{\infty} \ a_j(z).b^j) $$
  implies that \  $a_0 \equiv 1$ \ and  \ $a_1 = \frac{-1}{z}$ ! $\hfill \square$
  
   \begin{lemma}\label{mero}
   Let \ $f , g : X \to \C[[b]]$ \ be two holomorphic functions on a irreducible  reduced complex space \ $X$. Assume that \ $g $ \ does not vanish  and that, for each \ $x \in X$ \ the quotient \ $f(x)\big/g(x)$ \ is in \ $\C[[b]]$. Then there exists a Zariski open dense  set \ $X'$ \ in \ $X$ \ such that the map defined by \ $f\big/g$ \ is holomorphic on \ $X'$.
   \end{lemma}
 
 \parag{Proof} We may assume that \ $f \not\equiv 0$ \ on \ $X$. So there exists two  Zariski open dense  sets \ $X_1$ \ and \ $X_2$ \ such the valuation in \ $b$ \ of \ $f(x)$ \ (resp. of \ $g(x)$) is constant  \ $= k$ \ on \ $X_1$ (resp. is constant \ $= l$ \ on \ $X_2$). Then our assumption implies \ $k \geq l$ \ and it is clear \ that on the Zariski open dense set   \ $X' = X_1 \cap X_2$  \ the map \ $f\big/g$ \ is holomorphic. $\hfill \blacksquare$\\
 
 Remark that the set where \ $g$ \ vanishes is a closed analytic set. So when \ $g$ \ is not identically zero we may again find an open dense set to apply the lemma.\\
 
 \begin{cor}\label{exist. k-thm.}
 Let  \ $f : X \to E$ \ be a non identically zero  holomorphic map of an irreducible reduced complex space \ $X$ \ with value in a  regular (a,b)-module \ $E$. Then there exists a dense open set in \ $X$ \ on which the restriction of \ $f$ \ induces  a \ $k-$thematic holomorphic map via \ $x \mapsto \A.x$ \ where \ $k \leq rank(E)$.
 \end{cor}
 
 \parag{Proof} The point is that we may find an open dense set \ $X'$ \ on which the rank of \ $\A.x$ \ is maximal, thanks to the first lemma above. Then we solve a Cramer system with parameter on \ $X'$ \ to find the functions  \ $x \mapsto S_j(x)\in \C[[b]]$ \ which give the relation
  $$ a^k.f(x) = \sum_{j=0}^{k-1} \ S_j(x).a^j.f(x) .$$
  But these functions are "meromorphic". The second lemma above gives then a  open dense set  in \ $X'$ \ on which\ $f$ \  holomorphic and \ $k-$thematic. $\hfill \blacksquare$

\newpage

\section{Bibliography}

\begin{itemize}

\item{[A-G-V]} Arnold, V. \  Goussein-Zad\'e, S. \ Varchenko,A. {\it Singularit\'es des applications diff\'erentiables}, volume 2,  \'edition MIR, Moscou 1985.

\item{[B.93]} Barlet, D. {\it Theory of (a,b)-modules I}, in Complex Analysis and Geometry, Plenum Press New York (1993), p. 1-43.

\item{[B.95]} Barlet, D. {\it Theorie des (a,b)-modules II. Extensions}, in Complex Analysis and Geometry, Pitman Research Notes in Mathematics, Series 399 Longman (1997), p. 19-59.

\item{[B. 05]}  Barlet, D. \textit{Module de Brieskorn et formes hermitiennes pour une singularit\'e isol\'ee d'hypersurface}, Revue Inst. E. Cartan (Nancy) vol. 18 (2005), p. 19-46.

\item{[B.09b]} Barlet, D.  {\it Sur les fonctions \`a lieu singulier de dimension 1}, Bull. Soc. math. France 137 (4), (2009), p.587-612.

\item{[B.08]} Barlet, D. {\it Two finiteness theorems for (a,b)-modules}, preprint Inst. E. Cartan (Nancy) \ $n^0 5$ \ (2008)  p. 1-38.

\item{[B. 09a]} Barlet, D. {\it P\'eriodes \'evanescentes et (a,b)-modules monog\`enes},  Bulletino U.M.I. (9) II (2009), p. 651-697.

 \item{[Br.70]} Brieskorn, E. {\it Die Monodromie der Isolierten Singularit{\"a}ten von Hyperfl{\"a}chen}, Manuscripta Math. 2 (1970), p. 103-161.

\item{[K.09]} Karwasz, P. {\it Self dual (a,b)-modules and hermitian forms}, Th\`ese de doctorat de l'Universit\'e H. Poincar\'e (Nancy 1) soutenue le 10/12/09, p.1 -57.

\item{[K.76]} Kashiwara, M.  {\it b- function and holonomic systems}, Inv. Math. 38 (1976) p. 33-53.

\item{[M.74]} Malgrange, B.   {\it Int\'egrales asymptotiques et monodromie}, Ann. Sci. Ec. Norm. Sup. 4 (1974) p. 405-430.

  \item{[M. 75]} Malgrange, B. {\it Le polyn\^ome de Bernstein d'une singularit\'e isol\'ee}, in Lect. Notes in Math. 459, Springer (1975), p.98-119.

\item{[S. 89]} Saito, M. {\it On the structure of Brieskorn lattices}, Ann. Inst. Fourier 39 (1989), p.27-72.

\end{itemize}

\end{document}